\documentclass[10pt,a4paper,oneside,onecolumn,fleqn]{article}
\usepackage{mathrsfs}
\usepackage{amsfonts}
\usepackage{amssymb}
\usepackage{latexsym}
\usepackage{amsmath,amsthm}
\usepackage{epsf}
\usepackage{graphicx}
\usepackage{indentfirst}
\usepackage{cite}
\usepackage[pagewise]{lineno}
\usepackage{setspace}
\usepackage{color}
\usepackage{caption}

\theoremstyle{plain}

\theoremstyle{definition}

\theoremstyle{remark}

\title{\bf Further Results on the Pseudo-$L_{g}(s)$ Association Scheme with $g\geq 3$, $s\geq g+2$
 }
\author{  Congwei Wang, Shanqi Pang, Guangzhou Chen$^*$ \\
\\
\small {\it Henan Engineering Laboratory for Big Data Statistical Analysis and Optimal Control, }\\
\small {\it School of Mathematics
and Information Science, Henan Normal University, Xinxiang, 453007, P.R.China}\\
}
\date{}
\begin{document}
\maketitle
\begin{center}
{\bf Abstract}
\end{center}

\vspace{0.2cm} \ \ \ \ It is inevitable that the $L_{g}(s)$ association
scheme with $g\geq 3, s\geq g+2$ is a pseudo-$L_{g}(s)$ association
scheme. On the contrary, although $s^2$ treatments of the
pseudo-$L_{g}(s)$ association scheme can form one $L_{g}(s)$
association scheme, it is not always an $L_{g}(s)$ association
scheme. Mainly because the set of cardinality $s$, which contains two
first-associates treatments of the pseudo-$L_{g}(s)$ association scheme, is non-unique.
Whether the order $s$ of a Latin square $\mathbf{L}$ is a
prime power or not, the paper proposes two new conditions in order to
extend a $POL(s,w)$ containing $\mathbf{L}$. It has been known that a
$POL(s,w)$ can be extended to a $POL(s,s-1)$
so long as Bruck's \cite{brh} condition $s\geq
\frac{(s-1-w)^4-2(s-1-w)^3+2(s-1-w)^2+(s-1-w)}{2}$ is satisfied,
Bruck's condition will be completely improved through utilizing six properties
of the $L_{w+2}(s)$ association scheme in this paper. Several
examples are given to elucidate the application of our results.

\textbf{Keywords:} Latin square; transversal; net; $L_{g}(s)$ association scheme;
pseudo-$L_{g}(s)$ association scheme


\vskip 0.5cm

 {  \begingroup\makeatletter  \let\@makefnmark\relax  \footnotetext{ $^*$Corresponding author:  chenguangzhou0808@163.com(G. Chen)  }  \makeatother\endgroup}

\section{Introduction}

The construction of pairwise or mutually orthogonal Latin squares has
fascinated researchers for many years. The known results are well documented in the books
by D$\acute{e}$nes and Keedwell \cite{dk, dk1} and Laywine and Mullen \cite{lm},
or the article by Jungnickel \cite{j}. For further information about Latin squares also see
the works by Beth, Jungnickel and Lenz \cite{bjl}, Dinitz and Stinson \cite{ds},
Raghavarao \cite{r}, Street and Street \cite{s}.
Now we wonder under what conditions the $POL(s,w)$ with $w\geq 1$ can be extended to a $POL(s,s-1)$.

It has already been known from Bose
and Shimamoto \cite{bs} that the existence of an $L_g(s)$ association
scheme presupposes the existence of a $POL(s, g-2)$, which requires
the structure of pairwise orthogonal Latin squares.
We may verify that an $L_g(s)$ association scheme is a
pseudo-$L_{g}(s)$ association scheme.

With respect to given positive integers $g$ and $s$ with $g\leq s$, many scholars have made some
progress. For $g=1$, Bose and Connor \cite{bc} have already proved that a
pseudo-$L_1(s)$ association scheme becomes unique, moreover, it is an
$L_1(s)$ association scheme. For $g=2$, Shrikhande \cite{ss} has shown
that all the pseudo-$L_2(s)$ association schemes are the $L_2(s)$
association schemes except when $s=4$. He also points out that
there are two distinct pseudo-$L_2(4)$ association schemes, one of
which is exactly an $L_2(4)$ association scheme, the other is not.
Moreover, if $s>4$, he \cite{ss1} showed that two Latin squares can
be added to any $POL(s,s-3)$ to obtain a $POL(s,s-1)$.

For $g=3$, Liu \cite{lwr} has demonstrated all the pseudo-$L_3(s)$
association schemes with $s=3, 4$ or $s\geq 24$ are
the $L_3(s)$ association schemes, and given a pseudo-$L_3(5)$ association scheme
which is not the $L_3(5)$ association scheme. Next Hsu \cite{xbl}
illustrates that there exists a pseudo-$L_3(6)$ association scheme that is not the
$L_3(6)$ association scheme.

Up to present we have not known whether there exists a $POL(10,3)$.
If it is assumed that the $POL(10,3)$ exists, Liu \cite{lzw} has given a
pseudo-$L_5(10)$ association scheme which is not an $L_5(10)$
association scheme.

Finally, if the condition $s\geq \frac{g^4-2g^3+2g^2+g}{2}$ is satisfied,
Bruck \cite{brh} has concluded
that any pseudo-$L_g(s)$ association scheme is uniquely an $L_g(s)$
association scheme. Let $g=s-1-w$, no matter how a $POL(s,w)$ with $w\geq 1$
was constructed, in this case it is verified that the a $POL(s,w)$ can be
extended to a $POL(s,s-1)$.

In this paper we say $s^2$ treatments of the pseudo-$L_{g}(s)$ association scheme can form
an $L_g(s)$ association scheme, which has two hiding implications.
One of these two ones is that it is an $L_g(s)$ association scheme by itself, the other is although
its $s^2$ treatments may form at least one $L_g(s)$ association scheme,
but it is not an $L_g(s)$ association scheme.
Now we conversely consider when $s^2$ treatments of the
pseudo-$L_{g}(s)$ association scheme can from an $L_g(s)$ association scheme.

In order that we further illustrate the $POL(s,w)$ with $w\geq 1$ can be extended to a $POL(s,s-1)$
while filling an $s\times s$ square with $s^2$ treatments $1, 2,
\cdots, s^2$ of the pseudo-$L^{*}_{[s-1-w]}(s)$ association scheme,
we principally study under what conditions
$s^2$ treatments of the pseudo-$L^{*}_{[s-1-w]}(s)$ association scheme
can form an $L_{s-1-w}(s)$ association scheme
if $$w+4\leq s< \frac{(s-1-w)^4-2(s-1-w)^3+2(s-1-w)^2+(s-1-w)}{2}$$
is satisfied. Section 2 contains some basic concepts and three lemmas.
According to six properties of an $L_{g}(s)$ association
scheme, Section 3 presents some primary results, so that it can be identified that $s^2$ treatments of
the pseudo-$L_{g}(s)$ association scheme form an $L_{g}(s)$ association scheme,
however, the pseudo-$L_{g}(s)$ association scheme is not always an $L_{g}(s)$ association scheme.
In Section 4, Examples 4.1 and 4.2 point out: it is false that all the pseudo-$L_3(s)$
association schemes with $s=3, 4$ are the $L_3(s)$ association schemes \cite{lwr};
the remaining examples are exhibited to illustrate the application of our theorems.
An open question is put forward in Section 5.

\section{Basic concepts and main lemmas}

We first introduce some concepts and four lemmas.

{\bf Definition 2.1} Two $s\times s$ matrices with entries from a
set $S$ of cardinality $s$ are said to be orthogonal to each other
if when one is superimposed on the other the ordered pairs $(i,j)$
of corresponding entries consist of all possible $s^2$ pairs.

{\bf Definition 2.2 \cite{hss} } A Latin square of order $s$ is an $s\times s$
array with entries from a set $S$ of cardinality $s$ such
that each element of $S$ appears once in every row and every column.

{\bf Definition 2.3 \cite{hss} } Two Latin squares of order $s$ are said
to be orthogonal to each other if when one is superimposed on the
other the ordered pairs $(i,j)$ of corresponding entries consist of
all possible $s^2$ pairs.

{\bf Definition 2.4 \cite{hss} } A collection of $w$ Latin squares of
order $s$, any pair of which are orthogonal, is called a set of
pairwise orthogonal Latin squares, and denoted by $POL(s,w)$.

We allow w=1: any Latin square of order $s$ is a $POL(s,1)$.
A $POL(s,s-1)$ will be referred to as a complete set of pairwise orthogonal Latin
squares of order $s$. It follows that $w\leq s-1$ for any $POL(s,w)$.

{\bf Definition 2.5 \cite{hss} } A transversal in a Latin square
of order $s$ is defined to be a set of $s$ row-column pairs such that
each row and column occurs once and the entries
corresponding to these row-column pairs are distinct.

A partial transversal of a Latin square
of order $s$ denotes a set of row-column pairs as Definition 2.5, whose number
is not more than $s$.

{\bf Definition 2.6 \cite{hss} } Two or more transversals of a Latin
square are said to be parallel if the corresponding sets of
row-column pairs are disjoint.

If we label the rows and columns of a Latin square of order $s$ by
$0,1,\cdots,s-1$ and denote its entries by
$l_{ij}\in\{0,1,\cdots,s-1\}$. In order to state the following
arguments intuitively and clearly, from now on, for a partial
transversal of this Latin square, after its row-column pair
$(i,j)$ is replaced by the treatment $t_{ij}=is+j+1$, $t_{ij}$
corresponds to $l_{ij}$.

Set $POL(s, w)=\{\mathbf{L_{1}, \cdots, L_{w}}\}$,
and let $\mathbf{T}=(t_{ij})$ and $\mathbf{L_{k}}=(l^{k}_{ij})$, i.e.,
\begin{center}
$\mathbf{T}=\left(
\begin{smallmatrix}
1 & 2 & \cdots & s \\
s+1 &s+2& \cdots & 2s \\
\vdots& \vdots & \vdots &\vdots \\
s^{2}-s+1&s^{2}-s+2&\cdots &s^{2} \\
\end{smallmatrix}
\right),\ \ \mathbf{L_{k}}=\left(
\begin{smallmatrix}
l^{k}_{00}&l^{k}_{01}& \cdots &l^{k}_{0(s-1)}\\
l^{k}_{10}&l^{k}_{11}& \cdots &l^{k}_{1(s-1)}\\
 \vdots& \vdots & \vdots &\vdots \\
l^{k}_{(s-1)0}&l^{k}_{(s-1)1}& \cdots &l^{k}_{(s-1)(s-1)}\\
\end{smallmatrix}
\right)$.
\end{center}
Here $t_{ij}=is+j+1, k=1,\cdots,w$; $i,j=0,1,\cdots,s-1$.

It is obviously seen that
$t_{ij}$ corresponds to each $l^{k}_{ij}, k=1,\cdots,w$. Taking any
treatment $t_1$ from $\{1,2,\cdots, s^{2}\}$ of $\mathbf{T}$, we
have the following

{\bf Definition 2.7 } A deletion operation $\mathbf{\Delta^{T; L_{1},
\cdots, L_{w}}_{t_1}}$ on $\mathbf{T}$ is defined to delete
three-part treatments concerning $t_1$: Part I and $t_1$ lie on the
same row of $\mathbf{T}$, Part II and $t_1$ lie on the same column
of $\mathbf{T}$, Part III correspond to the same symbol as $t_1$ in
each $\mathbf{L_{k}}, k=1,\cdots,w$.

Select $t_1, t_2, \cdots, t_j$ and $u_1, u_2, \cdots, u_m$ from
$\{1,2,\cdots, s^{2}\}$, respectively; for any two treatments of
either of $t_1, t_2, \cdots, t_j$ and $u_1, u_2, \cdots, u_m$, it is satisfied that these
two treatments appear neither in the same row nor in the same column
of $\mathbf{T}$, and that they correspond to two distinct symbols of
each $\mathbf{L_{k}}, k=1,\cdots,w$. We therefore have the following

{\bf Definition 2.8 } $\mathbf{\Delta^{T; L_{1}, \cdots,
L_{w}}_{t_1--t_2--\cdots--t_j}}$ on $\mathbf{T}$ will denote that we
delete three-part treatments concerning $t_i$ in turn,
$i=1,2,\cdots,j$, here $1\leq j\leq s$.

When $j< s$, if it happened that there is no treatment in some row
or column of $\mathbf{\Delta^{T; L_{1}, \cdots,
L_{w}}_{t_1--t_2--\cdots--t_{j}}}$, then the deletion operation
pauses. It follows that $\{t_1-t_2-\cdots-t_j\}$ is one common
partial transversal of $\mathbf{L_{1}, \cdots, L_{w}}$.

When $j=s$, it must be shown that there is exactly one treatment in
each row and each column of $\mathbf{\Delta^{T; L_{1}, \cdots,
L_{w}}_{t_1--t_2--\cdots--t_{s}}}$, moreover, $t_1, t_2, \cdots,
t_{s}$ correspond to $s$ distinct symbols in each $\mathbf{L_{k}},
k=1,\cdots,w$. It follows that $\{t_1-t_2-\cdots-t_{s}\}$ is one
common transversal of $\mathbf{L_{1}, \cdots, L_{w}}$.

{\bf Definition 2.9 } $\mathbf{C^{T; L_{1}, \cdots,
L_{w}}_{u_1--u_2--\cdots--u_m}}$ will denote the collection of all
the common transversals of $\mathbf{L_{1}, \cdots, L_{w}}$ which
contain $m$ distinct treatments $u_1, u_2, \cdots, u_m$, here $1\leq
m\leq s$.

If there exists no common transversal of $\mathbf{L_{1}, \cdots,
L_{w}}$ which contain $u_1, u_2, \cdots, u_m$, then we define
$\mathbf{C^{T; L_{1}, \cdots,
L_{w}}_{u_1--u_2--\cdots--u_m}}=\emptyset$, here $\emptyset$ is the
vacuous set.

If $\{u_1-u_2-\cdots-u_{s}\}$ is one common transversal of
$\mathbf{L_{1}, \cdots, L_{w}}$, we therefore have
\begin{center}
$\mathbf{C^{T; L_{1}, \cdots, L_{w}}_{u_1--u_2--\cdots--u_{s}}}=\{u_1-u_2-\cdots-u_{s}\}.$
\end{center}

We notice: it is possible that the
cardinality of $\mathbf{C^{T; L_{1}}_{u_1--u_2--\cdots--u_m}}$ is
tremendous. However, the more the number $w$ is, the smaller the
cardinality of $\mathbf{C^{T; L_{1}, \cdots,
L_{w}}_{u_1--u_2--\cdots--u_m}}$ is.

With respect to any two treatments of either of $t_1-t_2-\cdots-t_j$
in $\mathbf{\Delta^{T; L_{1}, \cdots,
L_{w}}_{t_1--t_2--\cdots--t_j}}$ and $u_1-u_2-\cdots-u_{m}$ in
$\mathbf{C^{T; L_{1}, \cdots, L_{w}}_{u_1--u_2--\cdots--u_m}}$, it
is easily verified that the commutative law holds.

{\bf Definition 2.10 \cite{bs, bn} } Given $v$ distinct symbols $1, 2,
\cdots, v$, an association scheme with two associate classes is
referred to an associate relation between any two symbols satisfying
three following conditions:

(i) Any two symbols are either first or second associates, the
relation is symmetric.

(ii) Each symbol has $n_i$ $i$-th associates, $i=1,2$,
where the number $n_i$ is independent of this symbol selected.

(iii) With respect to any two symbols $\alpha, \beta$, which are
$i$-th associates, the number of symbols that are commonly $j$-th
associates of $\alpha$ and $k$-th associates of $\beta$ is $p_{jk}^{i}, j, k=1,2$,
and it is independent of the pair of $\alpha$ and $\beta$.

The numbers $v, n_{i}$ and $p_{jk}^{i}$ are called the parameters of
the association scheme. Show $p_{jk}^{i}=p_{kj}^{i}$ easily,
it follows that intersection matrix $P_{i}=(p_{jk}^{i})$ is symmetrical.

{\bf Definition 2.11 \cite{bs, lwr, lzw} } Presuppose that there exists a $POL(s,
g-2)$ $(g\leq s)$. An association scheme with two associate classes
is called an $L_{g}(s)$ association scheme, if there are $v=s^{2}$
treatments that may be set forth in an $s\times s$ array, such that
any two distinct treatments are first associates, which either occur
together in the same row or column of the array, or correspond to the same symbol
of one Latin square of the $POL(s, g-2)$; and second associates
otherwise.

{\bf Definition 2.12 \cite{brh} } A system $N$ satisfying the following
statements (I)-(IV) we shall call a net of \emph{order} $n$,
\emph{degree} $k$:

(I) Each line of $N$ contains exactly $n$ distinct points, where
$n$$\geq$$1$.

(II) Each point of $N$ lies on exactly $k$ distinct lines, where
$k$$\geq$$1$.

(III) $N$ has exactly $kn$ distinct lines. These fall into $k$
parallel classes of $n$ lines each. Distinct lines of the same
parallel class have no common points. Two lines of different classes
have exactly one common point.

(IV) $N$ has exactly $n^2$ distinct points.

{\bf Definition 2.13 \cite{brh} } For two distinct points P, Q of
$N$ in Definition 2.12, we say that P, Q are $joined$ in $N$
if there exists a line PQ of $N$ (necessarily unique) which contains
both P and Q; if the line PQ does not exist, then P, Q are
$not$ $joined$ in $N$.

{\bf Definition 2.14 \cite{brh} } A partial transversal of $N$ in Definition 2.12
denotes a nonempty set S of points of $N$ such that every two
distinct points in S are not joined in $N$. A transversal of $N$
denotes a partial transversal with exactly $n$ distinct points
(where $n$ is the order of $N$).

{\bf Definition 2.15 \cite{lwr, lzw} } An association scheme with two
associate classes having the following parameters is called a
pseudo-$L_{g}(s)$ association scheme:
\begin{center}
$v=s^2, n_{1}=g(s-1), p_{11}^{1}=(s-2)+(g-1)(g-2), p_{11}^{2}=g(g-1).$
\end{center}

It is easily seen that an $L_{g}(s)$ association scheme is a
pseudo-$L_{g}(s)$ association scheme \cite{lzw}. It is also verified that $n^2$ distinct
points of $N$ in Definition 2.12 form a pseudo-$L_{k}(n)$ association scheme; when we define two
distinct points of $N$ are first associates if and only if they are joined in $N$, and second
associates otherwise \cite{brh}.

{\bf Definition 2.16} A pseudo-$L^{*}_{[s+1-g]}(s)$ association scheme
is induced by an $L_{g}(s)$ association scheme if and only if its first
associates are exactly second associates of the $L_{g}(s)$
association scheme, and second associates otherwise.

It is obviously seen that an $L_{g}(s)$ association scheme and its inducing
pseudo-$L^{*}_{[s+1-g]}(s)$ association scheme have the same $s^2$
treatments. For two distinct treatments of the
pseudo-$L^{*}_{[s+1-g]}(s)$ association scheme, they are first
associates if and only if there exists a common partial transversal
of a $POL(s,g-2)$ (possibly nonunique) which contains them; and
second associates otherwise. That is to say, its two treatments are
first associates if and only if they appear neither in the same row
nor in the same column of an $s\times s$ array of the $L_{g}(s)$
association scheme, further, they correspond to two distinct symbols
of each of the $POL(s,g-2)$; and second associates otherwise.

{\bf Definition 2.17} A pseudo-net-$N^{*}$ of \emph{order} $s$,
\emph{degree} $s+1-g$ is induced by a net $N$ of \emph{order} $s$,
\emph{degree} $g$ if the following statements are true:

For two distinct points P, Q of pseudo-net-$N^{*}$, we say that
P, Q are $joined$ in pseudo-net-$N^{*}$ if and only if they are
$not$ $joined$ in $N$; P, Q are $not$ $joined$ in pseudo-net-$N^{*}$
if and only if they are $joined$ in $N$.

It is obviously known that a net $N$ and its inducing pseudo-net-$N^{*}$
have the same $s^2$ points. For two distinct points P, Q of pseudo-net-$N^{*}$,
if there exists a partial transversal of $N$ (possibly nonunique)
which contains both P and Q, then P, Q are joined in pseudo-net-$N^{*}$; if the partial
transversal does not exist, then they are not joined.

{\bf Definition 2.18 } A net of \emph{order}
$s$, \emph{degree} $s+1-g$ is called be a complementary net of a net
$N$ of \emph{order} $s$, \emph{degree} $g$ if and only if
its points are identical with those of $N$ and its
lines are a suitably selected set of transversals of $N$.

So long as $s(s+1-g)$ distinct transversals could be selected from all the ones of $N$,
these following conditions are satisfied: those $s(s+1-g)$ transversals fall into
$s+1-g$ parallel classes of $s$ ones each, distinct transversals of the same
parallel class have no common points, two transversals of different classes
have exactly one common point. In this case we say $s^2$ points of the pseudo-net-$N^{*}$
in Definition 2.17 are used to arrange
a new net of order $s$, degree $s+1-g$. This has
two hiding implications, one is that the pseudo-net-$N^{*}$ is exactly the new net, i.e.,
a complementary net of $N$, the other is although we are able to find at least
$s(s+1-g)$ distinct lines of pseudo-net-$N^{*}$ (transversals of $N$)
satisfying these above conditions, but the pseudo-net-$N^{*}$ is not the new net.
Either implication, once $s(s+1-g)$ distinct transversals
satisfying these conditions could be selected,
we know the new net and $N$ have the same $s^2$ points, moreover, its line is
the suitably selected transversal of $N$, its transversal is uniquely the
line of $N$. Thus we can imbed them in an affine plane of order $s$ \cite{brh}.

{\bf Definition 2.19 \cite{lwr} } A pseudo-$L^{*}_{s+1-g}(s)$ association
scheme is induced by a pseudo-$L_{g}(s)$ association scheme if and only if its
first associates are exactly second associates of the
pseudo-$L_{g}(s)$ association scheme, and second associates otherwise.

{\bf Lemma 2.20 \cite{lwr} } If the parameters of a pseudo-$L_{g}(s)$
association scheme are $v=s^2, n_{1}=g(s-1), p_{11}^{1}=g^{2}-3g+s$
and $p_{11}^{2}=g(g-1)$, then the parameters of a
pseudo-$L^{*}_{s+1-g}(s)$ association scheme are $v^{*}=s^2,
n_{1}^{*}=(s+1-g)(s-1), p_{11}^{1*}=(s+1-g)^{2}-3(s+1-g)+s$ and
$p_{11}^{2*}=(s+1-g)(s-g)$.

{\bf Proof: } It is known from \cite{bn} that \begin{center}
$v=n_{1}+n_{2}+1, \;\;\;\;\;\;\;\; n_{1}p_{12}^{1}=n_{2}p_{11}^{2}$, \\
$1+p_{11}^{1}+p_{12}^{1}=p_{11}^{2}+p_{12}^{2}=n_{1},\;\;\;\;\;\;\;\;
p_{21}^{1}+p_{22}^{1}=1+p_{21}^{2}+p_{22}^{2}=n_{2}$.
\end{center}
Hence,we have \begin{center}$n_{1}^{*}=n_{2}=v-n_{1}-1=s^2-g(s-1)-1=(s+1-g)(s-1)$,\\
$p_{11}^{1*}=p_{22}^{2}=n_{2}-1-n_{1}+p_{11}^{2}=(s+1-g)(s-1)-1-g(s-1)+g(g-1)
=(s+1-g)^{2}-3(s+1-g)+s$,\\
$p_{11}^{2*}=p_{22}^{1}=n_{2}-n_{1}+1+p_{11}^{1}=(s+1-g)(s-1)-g(s-1)+1+(g^{2}-3g+s)
=(s+1-g)(s-g)$.
\end{center} The proof is complete.
$ \blacksquare $

Since an $L_{g}(s)$ association scheme is a pseudo-$L_{g}(s)$ association scheme,
by definitions 2.16 and 2.19, a pseudo-$L^{*}_{[s+1-g]}(s)$ association scheme is
a pseudo-$L^{*}_{s+1-g}(s)$ association scheme. It is known
from Lemma 2.20 that a pseudo-$L^{*}_{s+1-g}(s)$ association scheme
is a pseudo-$L_{s+1-g}(s)$ association scheme.
Hence we acquire that a pseudo-$L^{*}_{[s+1-g]}(s)$ association
scheme is a pseudo-$L_{s+1-g}(s)$ association scheme, and that $s^2$
points of pseudo-net-$N^{*}$ in Definition 2.17 form a
pseudo-$L_{s+1-g}(s)$ association scheme; when we define that two
distinct points of pseudo-net-$N^{*}$ are first associates if and
only if they are joined in pseudo-net-$N^{*}$, and second
associates otherwise.

{\bf Lemma 2.21 \cite{lwr} } The pseudo-$L_{3}(s)$ association scheme must
be an $L_{3}(s)$ association scheme if and only if the following two
conditions hold ($s>4$):

$\mathbf{(I)}$ Take any treatment $x$, all the treatments that are
first associates of $x$ can be divided into three pairwise disjoint
sets as follows: $Y=\{y_{1},y_{2},\cdots,y_{s-1}\}$,
$Z=\{z_{1},z_{2},\cdots,z_{s-1}\}$ and
$W=\{w_{1},w_{2},\cdots,w_{s-1}\}$, so that any two different
treatments are first associates in each of three sets, such as $y_i,
y_j$ in $Y$, $z_i, z_j$ in $Z$ and $w_i, w_j$ in $W$, respectively,
where $i, j=1,2,\cdots,s-1, i\neq j$.

$\mathbf{(I^{'})}$ For every $y_i\in Y$, there is only one treatment
in either of $Z$ and $W$ such that it is first associates of $y_i$;
for every $z_i\in Z$, there is only one treatment in either of $Y$
and $W$ such that it is first associates of $z_i$; for every $w_i\in
W$, there is only one treatment in either of $Y$ and $Z$ such that
it is first associates of $w_i$.

Now we are wondering when $s^2$ treatments of the
pseudo-$L^{*}_{[s+1-g]}(s)$ association scheme can form an $L_{s+1-g}(s)$
association scheme, and when $s^2$ points of the pseudo-net-$N^{*}$ can be
arranged a net of order $s$, degree $s+1-g$. These will be elaborated in Section 3.

{\bf Lemma 2.22 \cite{hss} } A $POL(s,w)$ can be extended to a
$POL(s,w+1)$ if and only if the $w$ Latin squares in the $POL(s,w)$
possess $s$ common parallel transversals.

\section{ Main results }

In this section $\mathbf{T}$ is
identical with that of Definition 2.7, the $POL(s,w)$ will
always be written as $\mathbf{L_1, \cdots, L_{w}}, w\leq s-1$,
the corresponding relation of $\mathbf{T}$ and
$\mathbf{L_{k}}$ is as in Definition 2.7, here each $k\in\{1, \cdots, w\}$.

Let $g=w+2$ and $g\leq s$, we research the existence of an $L_g(s)$ association
scheme always presupposing the existence of a $POL(s,w)$.
With respect to two distinct treatments that are first
associates of a pseudo-$L_{g}(s)$ association scheme, it probably
appears that the set of cardinality $s$ containing them, in which
arbitrary two distinct treatments are first associates of the
pseudo-$L_{g}(s)$ association scheme, is non-unique.
Due to the converse-negative proposition of an $L_{g}(s)$ association
scheme having six properties,
even if $s^2$ treatments of the pseudo-$L_{g}(s)$ association scheme
can form an $L_{g}(s)$ association scheme,
but the pseudo-$L_{g}(s)$ association scheme is not always
an $L_{g}(s)$ association scheme.\\

Let $C$=$\left(
\begin{smallmatrix}
0&0&\cdots&0\\
1&1&\cdots&1\\
\vdots&\vdots&\vdots&\vdots\\
s-1&s-1&\cdots&s-1\\
\end{smallmatrix}
\right)$\ \ and\
$D$=$\left(
\begin{smallmatrix}
0 & 1 & \cdots&  s-1 \\
0 & 1 & \cdots&  s-1 \\
\vdots&\vdots &\vdots& \vdots \\
0 & 1 & \cdots &  s-1 \\
\end{smallmatrix}
\right)$ be two $s\times s$ matrices. We have the following.\\

{\bf Theorem 3.1 } An $s\times s$ matrix $L$ with entries from
$\{0,1,\cdots, s-1\}$ is a Latin square if and only if $C, D$ and
$L$ are orthogonal, respectively.

{\bf Proof.} Since each row of $L$ is the permutation of
$0,1,\cdots, s-1$ if and only if $L$ and $C$ are orthogonal
according to Definition 2.1. Similarly, also since each column of $L$
is the permutation of $0,1,\cdots, s-1$ if and only if $L$ and $D$
are orthogonal. Hence the theorem holds. $ \blacksquare $

{\bf Theorem 3.2 } Let $e_1, e_2, \cdots, e_s$ denote $s$ treatments
of $\mathbf{T}$ together corresponding to the same symbol of a Latin square
of order $s$. Then any two of $e_1, e_2, \cdots, e_s$ appear neither
in the same row nor in the same column of $\mathbf{T}$.

{\bf Proof: } It follows from Definition 2.2. $ \blacksquare $

{\bf Theorem 3.3 } If $e_1, e_2, \cdots, e_s$ and $e^{'}_1, e^{'}_2,
\cdots, e^{'}_s$ correspond to one symbol and another of
a Latin square of order $s$, respectively, then
$\{e_1, e_2, \cdots, e_s\}\cap \{e^{'}_1, e^{'}_2, \cdots,
e^{'}_s\}=\emptyset$ holds.

{\bf Proof: } It follows from Definition 2.2. $ \blacksquare $

When $s\geq g=3,4,\cdots$, it has been known from Definition 2.11 that $s^2$ treatments
and $g-2$ pairwise orthogonal Latin squares of order $s$ are
used to construct an $L_{g}(s)$ association scheme. We therefore have

{\bf Theorem 3.4 } For $s\geq g\geq 3$, an $L_{g}(s)$ association
scheme has the following six properties:

$\mathbf{(i)}$ There are only $g$ distinct paralleling classifications such that
its $s^2$ treatments can be exactly divided into $s$
pairwise disjoint sets $B^i_1, B^i_2, \cdots, B^i_s$ inside the
$i$-th classification, $i=1,2,\cdots,g$.

$\mathbf{(ii)}$ With respect to any two different classifications
$B^i_1, B^i_2, \cdots, B^i_s$ and $B^j_1, B^j_2, \cdots, B^j_s$, $i,
j\in \{1,2,\cdots,g\}$; for any $m, n\in\{1,2,\cdots,s\}$, $B^i_m$
and $B^j_n$ have exactly one common treatment.

$\mathbf{(iii)}$ For $g(s-1)$ treatments, which are first associates
of any treatment $t$, there is a unique partition of them such that
they can be equally divided into $g$ pairwise
disjoint sets $A_1, A_2, \cdots, A_g$.

$\mathbf{(iv)}$ With respect to $B^i_{m_i}$ containing the treatment $t$,
where $B^i_{m_i}\subset\{B^i_1, B^i_2, \cdots, B^i_s\}$,
$i=1,2,\cdots,g$; there exactly exists a one to one mapping of
$B^i_{m_i}$ onto $A_{k_i}$ satisfying $B^i_{m_i}=\{t\}\cup A_{k_i}, k_i\in\{1,2,\cdots,g\}$.

$\mathbf{(v)}$ With respect to arbitrary two distinct treatments $t_1$ and $t_2$,
they are first associates. It is assumed that $g(s-1)$
treatments, which are first associates of $t_v$, can uniquely be equally
divided into $g$ pairwise disjoint sets $A^{v}_1, A^{v}_2, \cdots,
A^{v}_g$, $v=1,2$; we have
$$\left\{
\begin{array}{lll}
(a). \ \  \ \{t_1\}\cup A^{1}_{i_1}=\{t_2\}\cup
A^{2}_{j_1}, \ \ \ for \ \ \ t_1\in A^{2}_{j_1},\ t_2\in A^{1}_{i_1}.\\
(b).\ \  For\ any\ i_2\in
\{1,2,\cdots,g\}\setminus\{i_1\},\\
\ \ \ \ \ \ (\{t_1\}\cup A^{1}_{i_2})\cap(\{t_2\}\cup
A^{2}_{j_2})=\emptyset,\ \ \
a \ single\ j_2\in \{1,2,\cdots,g\}\setminus\{j_1\}.\\
\ \ \ \ \ \ (\{t_1\}\cup A^{1}_{i_2})\cap(\{t_2\}\cup
A^{2}_{j'})=\{t_{1j'}\}, \ \ each\ \ j'\in \{1,2,\cdots,g\}\setminus\{j_1, j_2\}.\\
(c).\ \ For\ any\ j_3\in\{1,2,\cdots,g\}\setminus\{j_1\},\\
\ \ \ \ \ \ (\{t_2\}\cup A^{2}_{j_3})\cap(\{t_1\}\cup A^{1}_{i_3})=\emptyset,
\ \ \ a\ single\ i_3\in \{1,2,\cdots,g\}\setminus\{i_1\}.\\
\ \ \ \ \ \ (\{t_2\}\cup A^{2}_{j_3})\cap(\{t_1\}\cup A^{1}_{i'})=\{t_{2i'}\},
\ \ each\ i'\in \{1,2,\cdots,g\}\setminus\{i_1, i_3\}.\\
Here,\ t_{1j'}\ and\ t_{2i'}\ are\ commonly\ first\ associates\ of\ both\ t_1\ and\ t_2.
\end{array}
\right. $$

$\mathbf{(vi)}$ For every $a_i\in A_i$, $i=1,2,\cdots,g$; there are
exactly $(g-2)$ treatments in $A_j$ such that they are first
associates of $a_i$, any $j \in \{1,2,\cdots,g\}$$\setminus$$\{i\}$.

Where arbitrary two distinct treatments are first associates of the $L_{g}(s)$ association
scheme in each $\{t\}\cup A_{i}$, $\{t_1\}\cup A^{1}_{i}$, $\{t_2\}\cup A^{2}_{i},
i=1,2,\cdots,g$; each $B^i_m, B^j_n, m, n=1,2,\cdots,s$,
respectively.

{\bf Proof: } It is known from Definition 2.11 that there is a $POL(s,g-2)$=$\{\mathbf{L_1,
L_2, \cdots, L_{g-2}}\}$, moreover, $s^{2}$ treatments of the $L_{g}(s)$
association scheme can be filled an $s\times s$ array.

Let each $B^1_m$ and each $B^2_n$ be every row and every column of the array,
respectively; these obtain Classifications 1 and 2. We define each $B^{k+2}_m$ is the set of $s$
treatments corresponding to the same symbol of $\mathbf{L_k}$, this
obtains Classification $(k+2)$, $k=1,2,\cdots,g-2$. It follows that
$\mathbf{(i)}$ is verified.

According to Definition 2.2, because of pairwise orthogonality of
$\mathbf{L_1, L_2, \cdots, L_{g-2}}$, hence it may be seen from
$\mathbf{(i)}$ that $\mathbf{(ii)}$ holds.

$\mathbf{(i)}$ and $\mathbf{(ii)} \Longrightarrow \mathbf{(iii)}$
and $\mathbf{(iv)}$

With respect to any treatment $t$, it is known from $\mathbf{(i)}$
that there is exactly one sole $B^i_{m_i}$ such that $t\in
B^i_{m_i}$ holds, where $B^i_{m_i}\subset\{B^i_1, B^i_2, \cdots,
B^i_s\}$, $i=1,2,\cdots,g$. Let $B^i_{m_i}=\{t\}\cup A_{i}$, by
$\mathbf{(ii)}$, for any $i_1, i_2\in\{1,2,\cdots,g\}, i_1\neq i_2$,
we have $A_{i_1}\cap A_{i_2}=\emptyset$. It follows that
$\mathbf{(iv)}$ holds.

Next we consider $g(s-1)$ treatments that are first associates of
$t$. Due to $\mathbf{(iv)}$, we exhibit $\mathbf{(iii)}$ is verified.

$\mathbf{(i)}$ and $\mathbf{(ii)} \Longrightarrow \mathbf{(v)}$

With respect to arbitrary two distinct treatments $t_1$ and $t_2$, which are
first associates. It is very obvious from $\mathbf{(iii)}$ that $g(s-1)$
treatments, which are first associates of $t_v$, can uniquely be equally
divided into $g$ pairwise disjoint sets $A^{v}_1, A^{v}_2, \cdots,
A^{v}_g$, $v=1,2$.

Let $B^i_{m_i}=\{t_1\}\cup A^{1}_{i}$ and $B^j_{n_j}=\{t_2\}\cup
A^{2}_{j}$, $i,j=1,2,\cdots,g$. Since $t_1$ and $t_2$ are first
associates, there exist some $A^{1}_{i_1}$ and some $A^{2}_{j_1}$
such that $t_2\in A^{1}_{i_1}$ and $t_1\in A^{2}_{j_1}$ are satisfied. Thus we have
$B^{i_1}_{m_{i_1}}=\{t_1\}\cup A^{1}_{i_1}$ and
$B^{j_1}_{n_{j_1}}=\{t_2\}\cup A^{2}_{j_1}$.

$(a)$. According to $\mathbf{(i)}$ and $\mathbf{(ii)}$, we learn
$B^{i_1}_{m_{i_1}}$ and $B^{j_1}_{n_{j_1}}$ have exactly one
or no common treatments, which implies, it is impossible that they have two treatments in common.
Due to $\{t_1, t_2\}\subset B^{i_1}_{m_{i_1}}\cap B^{j_1}_{n_{j_1}}$, hence we must
have $\{t_1\}\cup A^{1}_{i_1}=B^{i_1}_{m_{i_1}}=B^{j_1}_{n_{j_1}}=\{t_2\}\cup A^{2}_{j_1}$.

$(b)$. Taking any $i_2\in \{1,2,\cdots,g\}$$\setminus$$\{i_1\}$, we
write $B^{i_2}_{m_{i_2}}=\{t_1\}\cup A^{1}_{i_2}$. By
$\mathbf{(i)}$, except for the $j_1$-th paralleling classification, there exists one
unique $B^{j_2}_{n_{j_2}}=\{t_2\}\cup A^{2}_{j_2}$ inside the
remaining $g-1$ paralleling classifications such that $B^{i_2}_{m_{i_2}}$ and
$B^{j_2}_{n_{j_2}}$ belong to the same paralleling classification, i.e.,
$B^{i_2}_{m_{i_2}}\cap B^{j_2}_{n_{j_2}}=(\{t_1\}\cup A^{1}_{i_2})\cap(\{t_2\}\cup
A^{2}_{j_2})=\emptyset$, here a single
$j_2\in \{1,2,\cdots,g\}$$\setminus$$\{j_1\}$. By $\mathbf{(ii)}$, for each $j'\in
\{1,2,\cdots,g\}$$\setminus$$\{j_1, j_2\}$,
we deduce $B^{i_2}_{m_{i_2}}$ and $B^{j'}_{n_{j'}}=\{t_2\}\cup
A^{2}_{j'}$ belong separately to distinct classifications, thus they
exactly have one common treatment $t_{1j'}$, i.e.,
$B^{i_2}_{m_{i_2}}\cap B^{j'}_{n_{j'}}=(\{t_1\}\cup A^{1}_{i_2})\cap(\{t_2\}\cup
A^{2}_{j'})=\{t_{1j'}\}$.

$(c)$. We may derive those conclusions in a similar manner.

$\mathbf{(i)}$ and $\mathbf{(ii)} \Longrightarrow \mathbf{(vi)}$

For $t_2\in A^{1}_{i_1}$ and $t_1\in A^{2}_{j_1}$, it is seen from $\mathbf{(v)}$ that
$(\{t_1\}\cup A^{1}_{i_2})\cap(\{t_2\}\cup A^{2}_{j_2})=\emptyset$ and
$(\{t_1\}\cup A^{1}_{i_2})\cap(\{t_2\}\cup A^{2}_{j'})=\{t_{1j'}\}$
hold for each $j'\in \{1,2,\cdots,g\}$$\setminus$$\{j_1, j_2\}$.
This means there are exactly $(g-2)$ treatments in $A^{1}_{i_2}$
such that they are all first associates of $t_2$ for any $i_2\in
\{1,2,\cdots,g\}$$\setminus$$\{i_1\}$.

For every $a_i\in A_i, i=1,2,\cdots,g$ and $A_j$, any $j \in
\{1,2,\cdots,g\}$$\setminus$$\{i\}$, after $t_2, A^{1}_{i_1},
A^{1}_{i_2}$ are replaced by $a_i, A_i, A_j$, respectively; it
follows that $\mathbf{(vi)}$ is verified. The theorem is therefore completely
proved. $ \blacksquare $

Since an $L_g(s)$ association scheme has six properties of Theorem 3.4,
it has already been known that an $L_g(s)$ association scheme is a
pseudo-$L_{g}(s)$ association scheme. However, the pseudo-$L_{g}(s)$
association scheme may or may not have these six properties. Also
since the original proposition is always logically equivalent to its
converse-negative proposition. Therefore, so long as
any property of these six ones is not
satisfied in the pseudo-$L_{g}(s)$ association scheme, in this case
the pseudo-$L_{g}(s)$ association scheme must not be an $L_{g}(s)$
association scheme. For example, if there exists no set of cardinality $s$
containing the fixed treatment at all,
in which arbitrary two distinct treatments are first associates of the
pseudo-$L_{g}(s)$ association scheme, or if there is some
treatment in the pseudo-$L_{g}(s)$ association scheme such that
$g(s-1)$ treatments which are first associates of it, cannot be
equally divided into $g$ pairwise disjoint sets,
then $s^2$ treatments of the pseudo-$L_{g}(s)$ association scheme cannot form
an $L_{g}(s)$ association scheme in any case, that is, the pseudo-$L_{g}(s)$
association scheme must not be an $L_g(s)$ association scheme.

Hereafter in the following discuss, it is always assumed that
there exists the set of cardinality $s$ such that arbitrary
two distinct treatments within it are first associates of a
pseudo-$L_{g}(s)$ association scheme.
For $s\geq g\geq 3$, the following five conditions possibly hold in the
pseudo-$L_{g}(s)$ association scheme:

$\mathbf{(I)}$ There are more than $g$ distinct paralleling classifications such that
its $s^2$ treatments can be exactly divided into $s$
pairwise disjoint sets $B^i_1, B^i_2, \cdots, B^i_s$ inside the
$i$-th classification, $i=1,2,\cdots,g,\cdots$.

$\mathbf{(II)}$ With respect to two different classifications
$B^i_1, B^i_2, \cdots, B^i_s$ and $B^j_1, B^j_2, \cdots, B^j_s$, $i,
j\in \{1,2,\cdots,g,\cdots\}$; $B^i_{m_1}$
and $B^j_{n_1}$ have greater than one and less
than $s$ common treatments, here $m_1, n_1\in\{1,2,\cdots,s\}$.

$\mathbf{(III)}$ For $g(s-1)$ treatments, which are first associates
of some treatment $t_1$, there are more than one partitions of them such that
they can be equally divided into $g$ pairwise
disjoint sets $A^{1}_1, A^{1}_2, \cdots, A^{1}_g$.

$\mathbf{(IV)}$ With respect to $B^i_{m_i}$ containing some treatment $t_1$,
where $B^i_{m_i}\subset\{B^i_1, B^i_2, \cdots, B^i_s\}$,
$i\in \{1,2,\cdots,g,\cdots\}$; $B^i_{m_i}\neq \{t_1\}\cup A^{1}_{k}$ holds, here
any $k\in\{1,2,\cdots,g\}$.

$\mathbf{(V)}$ There exist two distinct treatments $t_1$ and $t_2$,
which are first associates of the pseudo-$L_{g}(s)$ association scheme;
it is assumed that $\{t_1\}\cup A^{1}_{i_1}, \{t_1\}\cup A^{1}_{i_2}$ and
$\{t_2\}\cup A^{2}_{j_1}, \{t_2\}\cup A^{2}_{j_2}, \{t_2\}\cup A^{2}_{j_3}$
are five sets of cardinality $s$, in which arbitrary two distinct treatments
are first associates of the pseudo-$L_{g}(s)$ association scheme,
$i_1\neq i_2, j_1\neq j_2, j_1\neq j_3, j_2\neq j_3$.
If
$$\left\{
\begin{array}{lll}
(1)\ \ (\{t_1\}\cup A^{1}_{i_1})\cap(\{t_2\}\cup A^{2}_{j_1})=\emptyset,
(\{t_1\}\cup A^{1}_{i_1})\cap(\{t_2\}\cup A^{2}_{j_2})=\emptyset,\\
(2)\ \ \{t_1\}\cup A^{1}_{i_2} \
and\ \{t_2\}\cup A^{2}_{j_3}\ have\ greater\ than\ one\ and\ less\ than\ $s$ \
common\ treatments.\\ In \ especial,\ the\ set\ of\ cardinality\ $s$\ containing\
t_1\ and \ t_2\ is\ non$-$unique,\ i.e.,\\
\ \ \ \ \ \ \ \ \ \ \ \ \ \ \ \ \ \ \ \ \ \
\{t_1\}\cup A^{1}_{i_2}\neq \{t_2\}\cup A^{2}_{j_3},\ for \ t_1\in A^{2}_{j_3}, t_2\in A^{1}_{i_2}.\\
\end{array}
\right. $$
are satisfied.

Where arbitrary two distinct treatments are first associates of the
pseudo-$L_{g}(s)$ association scheme within each
$\{t_1\}\cup A^{1}_{i}$, $\{t_2\}\cup A^{2}_{i},
i=1,2,\cdots,g,\cdots$; each $B^i_m, B^j_n, m, n=1,2,\cdots,s$,
respectively.

{\bf Theorem 3.5 } So long as any one of $\mathbf{(I), (II), (III), (IV), (V)}$ holds,
the pseudo-$L_{g}(s)$ association scheme must not be an $L_{g}(s)$ association scheme.
However, $s^2$ treatments of the pseudo-$L_{g}(s)$ association scheme probably form
an $L_{g}(s)$ association scheme.

{\bf Proof: } Since the original proposition is always logically equivalent to its
converse-negative proposition. Due to Theorem 3.4, the pseudo-$L_{g}(s)$ association
scheme must not be an $L_{g}(s)$ association scheme. However, we are able to select
$g$ distinct classifications satisfying $\mathbf{(i)}$ and $\mathbf{(ii)}$ of Theorem 3.6 from the
pseudo-$L_{g}(s)$ association scheme, so that its $s^2$ treatments
can form an $L_{g}(s)$ association scheme. Hence the theorem holds. $ \blacksquare $

Next what we will need to do: how to select which properties
of those six ones in Theorem 3.6 so that they are satisfied,
it follows that $s^2$ treatments of the pseudo-$L_{g}(s)$
association scheme can form an $L_g(s)$ association scheme.

{\bf Theorem 3.6 } $\mathbf{(i)}$ $s^2$ treatments of the pseudo-$L_{g}(s)$
association scheme can be exactly divided into $s$
pairwise disjoint sets $B^i_1, B^i_2, \cdots, B^i_s$ inside the
$i$-th paralleling classification, $i=1,2,\cdots,g$.

$\mathbf{(ii)}$ With respect to any two different classifications
$B^i_1, B^i_2, \cdots, B^i_s$ and $B^j_1, B^j_2, \cdots, B^j_s$, $i,
j\in \{1,2,\cdots,g\}$; for any $m, n\in\{1,2,\cdots,s\}$, $B^i_m$
and $B^j_n$ have exactly one common treatment.

$\mathbf{(iii)}$ For $g(s-1)$ treatments, which are first associates
of any treatment $t$ in the pseudo-$L_{g}(s)$
association scheme, they can be equally divided into $g$ pairwise
disjoint sets $A_1, A_2, \cdots, A_g$.

$\mathbf{(iv)}$ With respect to $B^i_{m_i}$ containing the treatment $t$,
where $B^i_{m_i}\subset\{B^i_1, B^i_2, \cdots, B^i_s\}$,
$i=1,2,\cdots,g$; there exactly exists a one to one mapping of
$B^i_{m_i}$ onto $A_{k_i}$ satisfying $B^i_{m_i}=\{t\}\cup A_{k_i}, k_i\in\{1,2,\cdots,g\}$.

$\mathbf{(v)}$ With respect to arbitrary two distinct treatments $t_1$ and $t_2$,
they are first associates of the pseudo-$L_{g}(s)$ association scheme, it is assumed that $g(s-1)$
treatments, which are first associates of $t_v$, can be equally
divided into $g$ pairwise disjoint sets $A^{v}_1, A^{v}_2, \cdots,
A^{v}_g$, $v=1,2$. If
$$\left\{
\begin{array}{lll}
(a).\ \ \{t_1\}\cup A^{1}_{i_1}=\{t_2\}\cup
A^{2}_{j_1}, \ \ for \ t_1\in A^{2}_{j_1}, t_2\in A^{1}_{i_1};\ moreover,\
\ t_1\ and \ t_2\ define\ the\ set\ of\\ \ \ \ \ \ \ \ cardinality\ s,\ in\ which\ arbitrary\ two\
distinct\ treatments\ are\ first\ associates\ of\ the\\ \ \ \ \ \ \ \ pseudo$-$L_{g}(s)\ association\ scheme.\\
(b).\ \ For\ any\ i_2\in \{1,2,\cdots,g\}\setminus\{i_1\},\\
\ \ \ \ \ \ (\{t_1\}\cup A^{1}_{i_2})\cap(\{t_2\}\cup A^{2}_{j_2})=\emptyset, \ \ a\ single\ j_2\in \{1,2,\cdots,g\}\setminus\{j_1\}.\\
\ \ \ \ \ \ (\{t_1\}\cup A^{1}_{i_2})\cap(\{t_2\}\cup
A^{2}_{j'})=\{t_{1j'}\},\ \ each\ j'\in \{1,2,\cdots,g\}\setminus\{j_1, j_2\}.\\
(c). \ \ For\ any\ j_3\in \{1,2,\cdots,g\}\setminus\{j_1\},\\
\ \ \ \ \ \ (\{t_2\}\cup A^{2}_{j_3})\cap(\{t_1\}\cup A^{1}_{i_3})=\emptyset,\ \ a\ single\ i_3\in \{1,2,\cdots,g\}\setminus\{i_1\}.\\
\ \ \ \ \ \ (\{t_2\}\cup A^{2}_{j_3})\cap(\{t_1\}\cup A^{1}_{i'})=\{t_{2i'}\},
\ each\ i'\in \{1,2,\cdots,g\}\setminus\{i_1, i_3\}.\\
\ \ Here,\ t_{1j'}\ and\ t_{2i'}\ are\ commonly\ first\ associates \ of\ both \ t_1\ and\ t_2.
\end{array}
\right. $$
are satisfied.

$\mathbf{(vi)}$ For every $a_i\in A_i$, $i=1,2,\cdots,g$; there are
exactly $(g-2)$ treatments in $A_j$ such that they are first
associates of $a_i$, any $j \in \{1,2,\cdots,g\}$$\setminus$$\{i\}$.

Then the following four consequences could be obtained:

I. If $\mathbf{(i)}$ and $\mathbf{(ii)}$ hold then $s^2$ treatments of the pseudo-$L_{g}(s)$
association scheme can form an $L_{g}(s)$ association scheme.

II. $\mathbf{(i)}$, $\mathbf{(iii)}$ and
$\mathbf{(iv)}$ are equivalent to $\mathbf{(i)}$ and $\mathbf{(ii)}$
in the pseudo-$L_{g}(s)$ association scheme.

III. $\mathbf{(v)}$ is equivalent to $\mathbf{(i)}$ and
$\mathbf{(ii)}$ in the pseudo-$L_{g}(s)$ association scheme.

IV. If $\mathbf{(iii)}$ and
$\mathbf{(vi)}$ with $s>(g-1)^2$ are both satisfied in the pseudo-$L_{g}(s)$
association scheme, then we have \\

$\left\{
\begin{array}{lll}
(1) \ \emph{the\ set\ of\ cardinality}\ s\ containing\ t\ and\ a_i\ is\
uniquely\ \{t\}\cup A_i,\ here\ a_i\in A_i;\\
(2) \ \mathbf{(v)}\ holds.\\
\end{array}
\right.$\\

Where arbitrary two distinct treatments are first associates of the
pseudo-$L_{g}(s)$ association scheme within each
$\{t\}\cup A_{i}$, $\{t_1\}\cup A^{1}_{i}$, $\{t_2\}\cup A^{2}_{i},
i=1,2,\cdots,g$; each $B^i_m, B^j_n, m, n=1,2,\cdots,s$,
respectively.

{\bf Consequence I } If $\mathbf{(i)}$ and $\mathbf{(ii)}$ hold
then $s^2$ treatments of the pseudo-$L_{g}(s)$
association scheme can form an $L_{g}(s)$ association scheme.

{\bf Proof: } By $\mathbf{(i)}$, $s^2$ treatments of the pseudo-$L_{g}(s)$
association scheme can be exactly divided into $s$ pairwise disjoint
sets $B^i_1, B^i_2, \cdots, B^i_s$ inside the $i$-th paralleling classification,
$i=1,2,\cdots,g$.

Let $B^1_m$ be the $m$-th row of an $s\times s$ array $T^{'},
m=1,2,\cdots,s$. Now we shall write $B^2_n$, any
$n\in\{1,2,\cdots,s\}$. Since $B^1_1, B^1_2, \cdots, B^1_s$ are
pairwise disjoint, each of them and $B^2_n$ have exactly one common
treatment, thus these $s$ treatments of intersection are exactly the whole
treatments of $B^2_n$. It is assumed that $s$ treatments of each $B^1_m$
of $T^{'}$ remain constant, we could exchange their positions on
the $m$-th row in order that the whole treatments of $B^2_n$ are together
located in the $n$-th column of a new $s\times s$ array $T$,
$n=1,2,\cdots,s$. It follows that $s$ rows and $s$ columns of $T$
are exactly $B^1_1, B^1_2, \cdots, B^1_s$ and $B^2_1, B^2_2, \cdots,
B^2_s$, respectively; moreover, $T$ is precisely filled with $s^2$ treatments
of the pseudo-$L_{g}(s)$ association scheme.

The next step is to construct $g-2$ mutually orthogonal Latin
squares of order $s$.

For $i=1,2,\cdots,g$, let $s$ treatments of $B^i_m$ correspond to
$s$ $(m-1)^{,}$s, $m=1,2,\cdots,s$. Such as,\begin{center}
take $i=1, \left(
\begin{array}{c}
B^1_1\\
B^1_2\\
\vdots\\
B^1_s\\
\end{array}
\right)$ corresponds to $C=\left(
\begin{smallmatrix}
0&0&\cdots&0\\
1&1&\cdots&1\\
\vdots&\vdots&\vdots&\vdots\\
s-1&s-1&\cdots&s-1\\
\end{smallmatrix}
\right)$;\\ take $i=2, \left(
\begin{array}{cccc}
B^2_1&B^2_2&\cdots&B^2_s\\
\end{array}
\right)$ corresponds to $D=\left(
\begin{smallmatrix}
0 & 1 & \cdots&  s-1 \\
0 & 1 & \cdots&  s-1 \\
\vdots&\vdots &\vdots& \vdots \\
0 & 1 & \cdots &  s-1 \\
\end{smallmatrix}
\right)$.\end{center}

Now we observe $s^2$ treatments of the pseudo-$L_{g}(s)$ association
scheme have already been fixed on $T$. Since $B^k_1, B^k_2, \cdots, B^k_s$
are pairwise disjoint, $i=k=3,4,\cdots,g$; we obviously display $s$
treatments of $B^k_m$ correspond to $s$ $(m-1)^{,}$s. Thus we show
$s$ treatments of $B^k_1$, $s$ treatments of $B^k_2, \cdots$, $s$
treatments of $B^k_s$ just correspond to an $s\times s$ array $L_{k}$,
next we shall prove $L_{k}$ is a Latin square of order $s$.

It is known from $\mathbf{(ii)}$ that $B^k_m$ and
each of $B^1_1, B^1_2, \cdots, B^1_s$ have exactly one common
treatment, one-to-one correspondingly, every pair of $(m-1,0),(m-1,1),\ldots,(m-1,s-1)$
occurs one time when $L_{k}$ is superimposed on $C$.
There are $s$ pairwise disjoint sets $B^k_1, B^k_2, \cdots, B^k_s$, we
comprehend $L_{k}$ and $C$ are orthogonal due to Definition 2.1.
Similarly, $L_{k}$ and $D$ are also orthogonal when $L_{k}$ is
superimposed on $D$. Also by Theorem 3.1, we verify that $L_{k}$ is a Latin square.

With respect to any $k_{1}, k_{2}\in\{3,4,\cdots,s\}, k_{1}\neq k_{2}$, by
$\mathbf{(ii)}$, each of $B^{k_{1}}_1, B^{k_{1}}_2, \cdots,
B^{k_{1}}_s$ and each of $B^{k_{2}}_1, B^{k_{2}}_2, \cdots,
B^{k_{2}}_s$ have exactly one common treatment, one-to-one correspondingly,
every pair of all possible $s^2$ ones occurs one time when
$L_{k_{1}}$ is superimposed on $L_{k_{2}}$. It follows from Definition 2.3 that
$L_{k_{1}}$ and $L_{k_{2}}$ are orthogonal.

In sum, we have constructed $g-2$ mutually orthogonal Latin squares,
denoted as $L_{3}, L_{4}, \cdots, L_{g}$, and set $POL(s,g-2)=\{L_{3}, L_{4},
\cdots, L_{g}\}$.

According to Definition 2.11, hence $s^2$ treatments of $T$ and
the above $POL(s,g-2)$ are used to
construct an $L_{g}(s)$ association scheme. That is to say, if
$\mathbf{(i)}$ and $\mathbf{(ii)}$ are satisfied in the
pseudo-$L_{g}(s)$ association scheme, then its $s^2$
treatments can form an $L_{g}(s)$ association scheme. The proof is
complete. $ \blacksquare $

Although we are able to choose $sg$ different sets of cardinality $s$
in each of which arbitrary two distinct treatments are first associates of a
pseudo-$L_{g}(s)$ association scheme, so that $\mathbf{(i)}$ and $\mathbf{(ii)}$
of Consequence I are satisfied, but the pseudo-$L_{g}(s)$
association scheme is not always an $L_{g}(s)$ association scheme.

{\bf Consequence II } $\mathbf{(i)}$, $\mathbf{(iii)}$ and
$\mathbf{(iv)}$ are equivalent to $\mathbf{(i)}$ and $\mathbf{(ii)}$
in the pseudo-$L_{g}(s)$ association scheme.

{\bf Proof: } $\mathbf{(i)}, \mathbf{(iii)}$ and $\mathbf{(iv)}
\Longrightarrow \mathbf{(i)}$ and $\mathbf{(ii)}$

Take any $m\in\{1,2,\cdots,s\}$, let $B^i_m=\{t_1, t_2, \cdots,
t_s\}$ in the $i$-th paralleling classification. By $\mathbf{(iii)}$,
we learn $g(s-1)$ treatments, which are
first associates of $t_v$, can be equally divided into $g$ pairwise
disjoint sets $A^{v}_1, A^{v}_2, \cdots, A^{v}_g$, $v=1,2,\cdots,s$.
Due to $t_v\in B^i_m$, also by $\mathbf{(iv)}$, there exactly exists
only one of $\{t_v\}\cup A^{v}_{1}, \{t_v\}\cup A^{v}_{2}, \cdots,
\{t_v\}\cup A^{v}_{g}$ such that it is equal to $B^i_m$. Without
loss of generality, set $B^i_m=\{t_v\}\cup A^{v}_{i}$, thus we
write $B^i_m=\{t_1\}\cup A^{1}_{i}=\{t_2\}\cup
A^{2}_{i}=\cdots=\{t_s\}\cup A^{s}_{i}$.

As far as $t_v$ located in the $j$-th paralleling classification $B^j_{1},
B^j_{2}, \cdots, B^j_{s}$ is concerned, $v=1,2,\cdots,s$, $j\neq i$, it is seen from
$\mathbf{(i)}$ that there is exactly one unique $B^j_{n_v}$ such
that $t_v\in B^j_{n_v}$ holds. Using $\mathbf{(iv)}$ again, there
exactly exists only one of $\{t_v\}\cup A^{v}_{1}, \{t_v\}\cup
A^{v}_{2}, \cdots, \{t_v\}\cup A^{v}_{g}$ such that it is equal to
$B^j_{n_v}$. Without loss of generality, we put
$B^j_{n_v}=\{t_v\}\cup A^{v}_{j}$. It is known from $\mathbf
{(iii)}$ that $A^{v}_{i}\cap A^{v}_{j}=\emptyset$ holds, hence
$B^i_m$ and $B^j_{n_v}$ have exactly one common treatment $t_v$,
i.e., $B^i_m\cap B^j_{n_v}=\{t_v\}$.

With respect to arbitrary $t_{v_1}, t_{v_2}$, here $v_1, v_2\in\{1,2,\cdots,s\},
v_1\neq v_2$, let $B^j_{n_{v_1}}=\{t_{v_1}\}\cup A^{v_1}_{j}$ and
$B^j_{n_{v_2}}=\{t_{v_2}\}\cup A^{v_2}_{j}$. Since $B^j_{n_{v_1}}$
and $B^j_{n_{v_2}}$ are derived from $B^j_{1}, B^j_{2}, \cdots,
B^j_{s}$, using $\mathbf{(i)}$, we can clarify either $B^j_{n_{v_1}}=B^j_{n_{v_2}}$
or $B^j_{n_{v_1}}\cap B^j_{n_{v_2}}=\emptyset$ holds.

Since $B^i_m=\{t_{v_1}\}\cup A^{v_1}_{i}=\{t_{v_2}\}\cup
A^{v_2}_{i}$, we have $t_{v_1}\in A^{v_2}_{i}$. Due to $A^{v_2}_{i}\cap
A^{v_2}_{j}=\emptyset$, we obtain $t_{v_1}\notin A^{v_2}_{j}$, this
indicates $t_{v_1}\notin B^j_{n_{v_2}}$. Thus we must have
$B^j_{n_{v_2}}\neq B^j_{n_{v_1}}$ because of $t_{v_1}\in
B^j_{n_{v_1}}$. Finally we acquire $B^j_{n_{v_1}}\cap
B^j_{n_{v_2}}=\emptyset$. It follows that $B^j_{n_1}, B^j_{n_2},
\cdots, B^j_{n_s}$ is just one permutation of $B^j_{1}, B^j_{2},
\cdots, B^j_{s}$. Therefore, $B^i_m$ and each of $B^j_{1}, B^j_{2},
\cdots, B^j_{s}$ have exactly one common treatment, that is, $B^i_m$
and $B^j_n$ have exactly one common treatment, $n=1,2,\cdots,s$.

It is also verified that $B^j_n$ and $B^i_m$ have exactly one common
treatment in a similar manner, $m=1,2,\cdots,s$. In the end, we may
summarize $\mathbf{(ii)}$ holds.

$\mathbf{(i)}, \mathbf{(iii)}$ and $\mathbf{(iv)} \Longleftarrow
\mathbf{(i)}$ and $\mathbf{(ii)}$

When $\mathbf{(i)}$ and $\mathbf{(ii)}$ hold, by Consequence I, $s^2$
treatments of the pseudo-$L_{g}(s)$ association scheme can form an
$L_{g}(s)$ association scheme. Also utilizing Theorem 3.4, we may elucidate
$\mathbf{(iii)}$ and $\mathbf{(iv)}$ hold. The proof is now
complete. $ \blacksquare $

As far as $(a)$ of $\mathbf{(v)}$ is concerned,
we say $t_1$ and $t_2$ define the set of cardinality $s$, which means,
once $t_1$ and $t_2$ are fixed, the set of cardinality $s$ containing them
will be restricted. Even though there are several sets of cardinality $s$
containing $t_1$ and $t_2$ in the pseudo-$L_{g}(s)$ association scheme, we
only select one of several sets satisfying $\mathbf{(v)}$, and take no account of the remainders again.
Because of the randomicity choice of $t_1$ and $t_2$, we comprehend $\mathbf{(v)}$ is still satisfied
after they are replaced with other two first-associates treatments.

Further, also we have seen from $\mathbf{(v)}$ that any set of $\{t_1\}\cup A^{1}_{1}, \{t_1\}\cup
A^{1}_{2}, \cdots, \{t_1\}\cup A^{1}_{g}$ and any set of
$\{t_2\}\cup A^{2}_{1}, \{t_2\}\cup A^{2}_{2}, \cdots, \{t_2\}\cup
A^{2}_{g}$ have exactly $\mathbf{M}$ common treatments, here
$\mathbf{M}\in\{0,1,s\}$.  Of course undeniably,
since $\mathbf{(I)}$ and $\mathbf{(III)}$ of the
pseudo-$L_{g}(s)$ association scheme probably hold in Theorem 3.5,
not only $\mathbf{(v)}$ of Theorem 3.6, but also it is satisfied that
there maybe exist other different sets of cardinality $s$ (which are distinct from
$\{t_1\}\cup A^{1}_{1}, \{t_1\}\cup A^{1}_{2}, \cdots, \{t_1\}\cup A^{1}_{g}$
and $\{t_2\}\cup A^{2}_{1}, \{t_2\}\cup A^{2}_{2}, \cdots, \{t_2\}\cup A^{2}_{g}$)
containing $t_1$ and $t_2$ in the pseudo-$L_{g}(s)$ association scheme,
such that arbitrary two different treatments
within each of other sets are first associates.

Even if $t_1$ and $t'_1$ are second associates, it is impossible that
they belong to the same set of cardinality $s$, in which
arbitrary two distinct treatments are first associates of the
pseudo-$L_{g}(s)$ association scheme. According to $\mathbf{(v)}$,
primarily we must deduce that any set of
$\{t_1\}\cup A^{1}_{1}, \{t_1\}\cup A^{1}_{2},
\cdots, \{t_1\}\cup A^{1}_{g}$ and any set of $\{t'_1\}\cup
A^{'1}_{1}, \{t'_1\}\cup A^{'1}_{2}, \cdots, \{t'_1\}\cup
A^{'1}_{g}$ have one or no common treatments.

Using proof by contradiction, now we suppose that $\{t_1\}\cup A^{1}_{i}$ and $\{t'_1\}\cup A^{'1}_{j}$
have greater than one and less than $s$ common treatments.
Without loss of generality, let two of these common treatments be $t_3$ and $t_4$,
which are obviously first associates of the pseudo-$L_{g}(s)$ association scheme.
Because arbitrary two distinct
treatments of $\{t_1\}\cup A^{1}_{i}$ (or $\{t'_1\}\cup A^{'1}_{j}$)
are first associates of the pseudo-$L_{g}(s)$ association scheme, there are two different sets
$\{t_1\}\cup A^{1}_{i}$ and $\{t'_1\}\cup A^{'1}_{j}$ of cardinality $s$ containing $t_3$ and $t_4$.
This contradicts $t_3$ and $t_4$ define the set of cardinality $s$
within $\mathbf{(v)}$.

Of course incontestably, since $\mathbf{(I)}$ and $\mathbf{(III)}$ of the
pseudo-$L_{g}(s)$ association scheme possibly hold in Theorem 3.5,
it is not only seen that each of $\{t_1\}\cup A^{1}_{1}, \{t_1\}\cup A^{1}_{2},
\cdots, \{t_1\}\cup A^{1}_{g}$ and each of $\{t'_1\}\cup
A^{'1}_{1}, \{t'_1\}\cup A^{'1}_{2}, \cdots, \{t'_1\}\cup
A^{'1}_{g}$ have one or no common treatments, but also satisfied that
there maybe exist other two different sets of cardinality $s$ (which are different from
$\{t_1\}\cup A^{1}_{1}, \{t_1\}\cup A^{1}_{2}, \cdots, \{t_1\}\cup A^{1}_{g}$
and $\{t'_1\}\cup A^{'1}_{1}, \{t'_1\}\cup A^{'1}_{2}, \cdots, \{t'_1\}\cup
A^{'1}_{g}$), such that arbitrary two distinct treatments within either of
these two sets are first associates in the pseudo-$L_{g}(s)$ association scheme.
Moreover, one of them contains $t_1$ and the other involves $t'_1$,
in the meantime, they have greater than one and less than $s$ common treatments.

{\bf Consequence III } $\mathbf{(v)}$ is equivalent to $\mathbf{(i)}$ and
$\mathbf{(ii)}$ in the pseudo-$L_{g}(s)$ association scheme.

{\bf Proof: } $\mathbf{(v)} \Longrightarrow \mathbf{(i)}$ and
$\mathbf{(ii)}$

When $\mathbf{(v)}$ holds at all times in the pseudo-$L_{g}(s)$
association scheme, thus we obtain immediately:

$\mathbf{(iii^{'})}$ For $g(s-1)$ treatments, which are first
associates of any treatment, they can be equally divided into
$g$ pairwise disjoint sets, and

$\mathbf{(vi^{'})}$ With respect to arbitrary two distinct
sets of these $g$ ones in $\mathbf{(iii^{'})}$, take any treatment
from one set, there are exactly $(g-2)$ treatments in
the other such that they are all first associates of this treatment.

With respect to arbitrary two treatments of first associates
in the pseudo-$L_{g}(s)$ association scheme,
so long as $\mathbf{(v)}$ of Theorem 3.6 is always satisfied,
these two treatments restrict the set of cardinality $s$, in which
any two different treatments are first associates.
We know $s^2$ treatments are set forth
in an $s\times s$ array so that they can form an $L_{g}(s)$ association scheme,
it is very necessary that there are $sg$ sets of $s$ distinct treatments such
that any two of these $s$ treatments are first associates, moreover,
those $sg$ sets fall into $g$ parallel classes of $s$ ones each, distinct sets of the same
parallel class have no common treatments, two sets of different classes
have exactly one common treatment. However, a pseudo-$L_{g}(s)$ association
scheme probably has much more than above $sg$ sets of $s$ distinct treatments required, so
that any two of these $s$ treatments are first associates of the
pseudo-$L_{g}(s)$ association scheme. Now the aim is to carefully select
above $sg$ sets required out from the pseudo-$L_{g}(s)$ association scheme,
hence we must proclaim the following four views.

Firstly, we desire each of $\{t_m\}\cup A^{m}_{1},
\{t_m\}\cup A^{m}_{2}, \cdots, \{t_m\}\cup A^{m}_{g}$ and each of
$\{t_n\}\cup A^{n}_{1}, \{t_n\}\cup A^{n}_{2}, \cdots,
\{t_n\}\cup A^{n}_{g}$ have exactly $\mathbf{M}$ common treatments, here
$\mathbf{M}\in\{0,1,s\}$. Secondly, if $B^i_m=\{t_m\}\cup A^{m}_i$ and $B^j_n=\{t_n\}\cup A^{n}_j$
have greater than one and less than $s$ common treatments, then we reselect
$g$ pairwise disjoint sets which those $g(s-1)$ treatments
(that are first associates of $t_m$ or $t_n$) can be divided into.
Thirdly, if above $g$ pairwise disjoint sets
satisfying $\mathbf{(v)}$ cannot be found, which means there exist distinct
$\{t_m\}\cup A^{m}_{i'}$ and $\{t_n\}\cup A^{n}_{j'}$ such that they
have greater than one and less than $s$ common treatments,
this contradicts the set of cardinality $s$
have been definitely standardized by two of these common treatments inside $\mathbf{(v)}$.
Fourthly, even though it perhaps appears that other two different sets of cardinality $s$
are distinct from $\{t_m\}\cup A^{m}_{1}, \{t_m\}\cup A^{m}_{2}, \cdots, \{t_m\}\cup A^{m}_{g}$
and $\{t_n\}\cup A^{n}_{1}, \{t_n\}\cup A^{n}_{2}, \cdots, \{t_n\}\cup A^{n}_{g}$.
Moreover, one of them contains $t_m$ and the other involves $t_n$,
in the meantime, they have two common treatments or more, we leave out of account.

According to $\mathbf{(iii^{'})}$, taking any treatment $t_1$ from
the pseudo-$L_{g}(s)$ association scheme, its first-associates
$g(s-1)$ treatments can be equally divided into $g$ pairwise
disjoint sets $A^{1}_1, A^{1}_2, \cdots, A^{1}_g$, we
put $B^1_1=\{t_1\}\cup A^{1}_1, B^2_1=\{t_1\}\cup A^{1}_2, \cdots,
B^g_1=\{t_1\}\cup A^{1}_g$.

Next we choose $B^i_1=\{t_1\}\cup A^{1}_i, B^j_1=\{t_1\}\cup
A^{1}_j$, here any $i,j\in\{1, 2, \cdots, g\}, i\neq j$.
Let $A^{1}_j=\{t_2, t_3, \cdots, t_s\}$. Due to
$\mathbf{(iii^{'})}$, it is assumed that $g(s-1)$ treatments, which are first
associates of $t_m$, can be equally divided into $g$ pairwise
disjoint sets $A^{m}_1, A^{m}_2, \cdots, A^{m}_g$, $m=2,3,\cdots,s$.
With respect to $t_1$ and each $t_m$,
they define the confined set $\{t_1, t_2, \cdots, t_s\}$.
Without loss of generality, let $\{t_{1}\}\cup
A^{1}_{j}=\{t_{m}\}\cup A^{m}_{j}=\{t_1, t_2, \cdots, t_s\}$.

Because of $\{t_{1}\}\cup A^{1}_{j}=\{t_{m}\}\cup A^{m}_{j}$, it is known
from $(b)$ of $\mathbf{(v)}$ that there
exists a single set of $\{t_m\}\cup A^{m}_1, \cdots, \{t_m\}\cup
A^{m}_{j-1}, \{t_m\}\cup A^{m}_{j+1}, \cdots, \{t_m\}\cup A^{m}_g$
such that it and $\{t_1\}\cup A^{1}_{i}=B^i_1$ have no common treatments.
Without loss of generality, we denote the single set of cardinality $s$ as $B^i_m=\{t_m\}\cup
A^{m}_{i}$, it follows that $B^i_m\cap B^i_1=\emptyset$.

The following arguments are to verify $B^i_{m_1}\cap
B^i_{m_2}=\emptyset$, where any $m_1, m_2\in\{2,3,\cdots,s\},\\ m_1\neq
m_2$, the entire process of proof may be shown as follows.

\begin{center}$\begin{array}{cc}
B^i_1&\\
\|&\\
A^1_i&\\
\supset&\\
t_1&\cup
\end{array}\underbrace{\begin{array}{ccccccccc}
&B^i_2\;\;\;\;\;\;\;\;&&B^i_{m_1}&&B^i_{m_2}&&\;\;\;\;\;\;\;\;B^i_s&\\
&\|\;\;\;\;\;\;\;\;&&\|&&\|&&\;\;\;\;\;\;\;\;\|&\\
&A^{2}_i\;\;\;\;\;\;\;\;&&A^{m_1}_i&&A^{m_2}_i&&\;\;\;\;\;\;\;\;A^{s}_i&\\
&\supset\;\;\;\;\;\;\;\;&&\supset&&\supset&&\;\;\;\;\;\;\;\;\supset&\\
\{&t_2--&\cdots&--t_{m_1}--&\cdots&--t_{m_2}--&\cdots&--t_s&\}
\end{array}}$\\$\begin{array}{ccccc}
&&&&A^1_j\\
\end{array}$\end{center}

Since $t_{m_1}$ and $t_{m_2}$ have defined the restricted set $\{t_1, t_2, \cdots, t_s\}$,
based on $\mathbf{(v)}$, thus $\{t_{m_1}\}\cup A^{m_1}_{i}=B^i_{m_1}$ and
$\{t_{m_2}\}\cup A^{m_2}_{i}=B^i_{m_2}$ of cardinality $s$ have one or no common treatments.
Otherwise, we may readjust $g$ pairwise disjoint sets which those $g(s-1)$ treatments
(that are first associates of $t_{m_1}$ or $t_{m_2}$) can be divided into,
so that $B^i_{m_1}$ and $B^i_{m_2}$ are satisfied.

Utilize reduction to absurdity, now it is assumed that $B^i_{m_1}$ and $B^i_{m_2}$ have exactly
one common treatment, labeled by $t_{u_1}$. This will
eventually lead to a contradiction, it goes without saying
that $B^i_{m_1}\cap B^i_{m_2}=\emptyset$ holds.

Due to $\mathbf{(iii^{'})}$, its first-associates
$g(s-1)$ treatments can be equally divided into $g$ pairwise
disjoint sets $A^{u_1}_1, A^{u_1}_2, \cdots, A^{u_1}_g$, it is
very obvious that $B^i_{m_1}$ and
$B^i_{m_2}$ belong to the collection $\{\{t_{u_1}\}\cup A^{u_1}_1, \{t_{u_1}\}\cup
A^{u_1}_2, \cdots,\{t_{u_1}\}\cup A^{u_1}_g\}$.

(1) When $t_{u_1}$ and $t_{1}$ are first associates, this indicates,
there are two distinct sets $B^i_{m_1}$ and $B^i_{m_2}$ meeting in
$t_{u_1}$ such that they are both disjoint from $B^i_{1}$ through
$t_{1}$. That directly contradicts $t_2$ of
$(b)$ is replaced by $t_{u_1}$ inside $\mathbf{(v)}$.

(2) When $t_{u_1}$ and $t_{1}$ are second associates, we suppose
there exists one set of $\{t_{u_1}\}\cup A^{u_1}_1, \{t_{u_1}\}\cup
A^{u_1}_2, \cdots,\{t_{u_1}\}\cup A^{u_1}_g$ such that it and
$B^i_{1}$ have one common treatment, labeled by $t_{u_2}$.
It follows that there are two distinct sets $B^i_{m_1}$
and $B^i_{m_2}$ meeting in $t_{u_1}$ such that they are both
disjoint from $B^i_{1}$ through $t_{u_2}$. This also contradicts
$t_1$ and $t_2$ of $(b)$ are replaced with $t_{u_2}$ and $t_{u_1}$
inside $\mathbf{(v)}$, respectively.

Consequently, we presume arbitrary set of $\{t_{u_1}\}\cup A^{u_1}_1,
\{t_{u_1}\}\cup A^{u_1}_2, \cdots,\{t_{u_1}\}\cup A^{u_1}_g$ and
$B^i_{1}$ have no common treatments while $t_{u_1}$ and $t_{1}$ are
second associates. Since two first-associates treatments
define the set of cardinality $s$, thus each of $\{t_{u_1}\}\cup
A^{u_1}_1, \{t_{u_1}\}\cup A^{u_1}_2, \cdots,\{t_{u_1}\}\cup
A^{u_1}_g$ and each of $B^1_1, \cdots, B^{i-1}_1, B^{i+1}_1, \cdots,
B^{g}_1$ have at most one common treatment, which implies there are
at most $g$ intersection "point"s among $\{t_{u_1}\}\cup A^{u_1}_1,
\{t_{u_1}\}\cup A^{u_1}_2, \cdots,\{t_{u_1}\}\cup A^{u_1}_g$ and
$B^{j^{'}}_{1}$, here any $j^{'}\in \{1,2,\cdots,g\}$$\setminus$$\{i\}$. Due to
$p_{11}^{2}=g(g-1)$. Hence there are exactly $g$ treatments in
$B^{j}_1$ such that they are first associates of $t_{u_1}$, here
$j\neq i$.

On the other hand, since $t_{m_1}$ and $t_{u_1}$ together belong to
$B^i_{m_1}$, due to $\mathbf{(vi^{'})}$, there are exactly $(g-2)$
treatments in $A^{m_1}_{j}$ such that they are first associates of
$t_{u_1}$ for $t_{u_1}\in A^{m_1}_{i}$, here $j\neq i$. Adding $t_{m_1}$, it follows
that there are exactly $(g-1)$ treatments in $B^{j}_1$ such that
they are first associates of $t_{u_1}$. This results in a
contradiction.

On balance, there exists no intersection "point" $t_{u_1}$ of
$B^i_{m_1}$ and $B^i_{m_2}$. Hence we demonstrate there truly exist
$s$ pairwise disjoint sets $B^i_1, B^i_2, \cdots, B^i_s$.
That is to say, $s^2$ treatments of the
pseudo-$L_{g}(s)$ association scheme can be exactly divided into $s$
pairwise disjoint sets $B^i_1, B^i_2, \cdots, B^i_s$.

Notice every set of $B^i_1, B^i_2, \cdots, B^i_s$, it is very obvious that arbitrary two
treatments are different among it. Because of non-intersect of
$A^{1}_i$ and $A^{1}_j$, we may pick $i=1, j=2$ out, $s^2$ treatments
of the pseudo-$L_{g}(s)$ association scheme can be exactly divided
into $s$ pairwise disjoint sets $B^1_1, B^1_2, \cdots, B^1_s$ inside
the first paralleling classification. Take $i=2,3,\cdots,g$ and
$B^1_1=\{t_1\}\cup A^{1}_1$ ($j=1$) again, we obtain
its $s^2$ treatments can be exactly divided into $s$ pairwise
disjoint sets $B^i_1, B^i_2, \cdots, B^i_s$ inside the $i$-th
paralleling classification. It follows that $\mathbf{(i)}$ holds.

Now we fix the $j$-th paralleling classification $B^j_1, B^j_2, \cdots, B^j_s$,
the variable $i$-th paralleling classification is focused attention on, here any
$i\in \{1,2,\cdots,g\}$$\setminus$$\{j\}$.
Still tagging $B^j_1\cap B^i_{m}=\{t_{m}\}$ as above, $m=1,2,\cdots,s$, next we investigate
arbitrary $B^j_{n'}, n'=2,3,\cdots,s$.

Since $s^2$ treatments of
the pseudo-$L_{g}(s)$ association scheme can be divided
into $s$ pairwise disjoint sets $B^i_1, B^i_2, \cdots, B^i_s$,
also because of $B^j_{n'}\cap B^j_{1}=\emptyset$ and $B^i_{m}\cap B^j_1=\{t_{m}\}$,
it is very certain that $B^j_{n'}$ is impossibly coincident
with anyone of $B^i_1, B^i_2, \cdots, B^i_s$.

Since $s$ treatments of $B^j_{n'}$ are distributed onto $B^i_1, B^i_2, \cdots$, $B^i_s$.
If $B^j_{n'}$ and someone of $B^i_1, B^i_2, \cdots, B^i_s$ have no common treatments,
then there exists at least another of them such that
it and $B^j_{n'}$ have two common treatments. That contradicts
these two common treatments define the decided set of cardinality
$s$ within $\mathbf{(v)}$.

It is impossible that $B^j_{n'}$ and each of $B^i_1, B^i_2, \cdots, B^i_s$
have greater than one and less than $s$ common treatments.
Otherwise, this also contradicts two of these common treatments
define the determined set of cardinality $s$ within $\mathbf{(v)}$.

It is summarized that $B^j_n$ and each of $B^i_1, B^i_2, \cdots, B^i_s$ have exactly
one common treatment. It follows that $\mathbf{(ii)}$ holds.

$\mathbf{(v)} \Longleftarrow \mathbf{(i)}$ and $\mathbf{(ii)}$

When $\mathbf{(i)}$ and $\mathbf{(ii)}$ hold, by Consequence I, $s^2$
treatments of the pseudo-$L_{g}(s)$ association scheme can form an
$L_{g}(s)$ association scheme. Also applying Theorem 3.4, we may clarify
$\mathbf{(v)}$ holds. The consequence is therefore completely
proved. $ \blacksquare $

With respect to any two distinct treatments $t_1$ and $t_2$,
which are first associates of a pseudo-$L_{g}(s)$
association scheme, we consider the set of cardinality $s$ containing them, in which
arbitrary two different treatments are first associates of the
pseudo-$L_{g}(s)$ association scheme.

Though we are able to carefully select desired $g$
classifications from much more than these $g$ ones of the
pseudo-$L_{g}(s)$ association scheme, so that $\mathbf{(v)}$
of Theorem 3.6 is satisfied, it follows that $s^2$
treatments of the pseudo-$L_{g}(s)$ association scheme could form an
$L_{g}(s)$ association scheme. But the pseudo-$L_{g}(s)$
association scheme is not always an $L_{g}(s)$ association scheme,
mainly because there possibly exist two first-associates treatments $t_3$ and $t_4$
such that the set of cardinality $s$ containing them is non-unique.

According to Consequence III, we therefore have the following

{\bf Corollary 3.7 } If the set of cardinality $s$ containing $t_1$ and
$t_2$ is unique, moreover, $\mathbf{(v)}$
of Theorem 3.6 is always satisfied, then the pseudo-$L_{g}(s)$
association scheme must be an $L_{g}(s)$ association scheme.

{\bf Consequence IV } If $\mathbf{(iii)}$ and
$\mathbf{(vi)}$ with $s>(g-1)^2$ are both satisfied in the pseudo-$L_{g}(s)$
association scheme, then we have

(1) the set of cardinality $s$ containing $t$ and $a_i$ is
uniquely $\{t\}\cup A_i$, where $a_i\in A_i$;  and

(2) $\mathbf{(v)}$ holds.

{\bf Proof: } (1) With respect to any treatment $t$, we have
$\{t\}\cup A_i$ from $\mathbf{(iii)}$, each $i\in\{1, 2, \cdots,
g\}$. For one fixed $a_i\in A_i$, it has already been known that the cardinality of
$\{t\}\cup A_i$ is $s$.

Next we suppose there exists another collection $\{t\}\cup C$ satisfying
$a_i\in C$ and $C\neq A_i$, such that arbitrary two treatments of
$\{t\}\cup C$ are first associates of the pseudo-$L_{g}(s)$
association scheme. Observing the maximal value of cardinality of
$C$, we detect how many it is.

Since $C\neq A_i$, there exists at least $c\in C$, but $c\notin
A_i$. Also since $c$ and $t$ are first associates, we must have
$c\in A_{j_1}$, here some $j_1\in\{1, 2, \cdots,
g\}$$\setminus$$\{i\}$. It is very obvious that $a_i$ and $c$ are
first associates, and that the remaining treatments of $C$ are first
associates of both $a_i$ and $c$. Due to $a_i\in A_{i}\cap C$, by
$\mathbf{(vi)}$, $C$ contains at most $(g-2)$ treatments in
$A_{j_1}$. Obviously, $c$ belongs to these $(g-2)$ treatments. Due
to $c\in A_{j_1}\cap C$, also by $\mathbf{(vi)}$, $C$ contains at most
$(g-2)$ treatments in each $A_{j'}, j' \in
\{1,2,\cdots,g\}$$\setminus$$\{j_1\}$. In particular, $a_i$ belongs
to those $(g-2)$ treatments in $A_i$.

In total, $C$ possesses at most $g(g-2)$ treatments that are all
first associates of $t$. When any two of them are first associates,
the cardinality of $C$ is $g(g-2)$; when there exist two treatments
of them such that they are second associates, its cardinality is less than
$g(g-2)$. In a word, the maximal value of cardinality of $C$ is
$g(g-2)$.

According to $(g-1)^2<s$, i.e., $g(g-2)<s-1$, it follows that the
maximal value of cardinality of $C$ is less than $s-1$.

With respect to any treatment $t$ and one fixed treatment $a_i$,
here $a_i\in A_i$, they define one sole set $\{t\}\cup A_i$ of
cardinality $s$ such that arbitrary two distinct treatments of $\{t\}\cup
A_i$ are first associates of the pseudo-$L_{g}(s)$ association
scheme. In other words, with respect to $t$ and one fixed $a_i$,
since the cardinality of $C$ is less than $s-1$,
it is certain that the set of cardinality $s$ containing
them is uniquely $\{t\}\cup A_i$.

(2) $\mathbf{(iii)}$ and $\mathbf{(vi)}$ with $s>(g-1)^2$
$\Longrightarrow \mathbf{(v)}$

With respect to arbitrary two distinct treatments $t_1$ and $t_2$, which are
first associates. It is already known from $\mathbf{(iii)}$ that
$g(s-1)$ treatments, which are first associates of $t_v$, can be
equally divided into $g$ pairwise disjoint sets $A^{v}_1, A^{v}_2,
\cdots, A^{v}_g, v=1,2$.

For $s>(g-1)^2$, with respect to $t_1$ and a given treatment of
$A^{1}_{i}$, due to (1), we understand the set of cardinality $s$ containing
them is uniquely $\{t_1\}\cup A^{1}_{i}$; so does $\{t_2\}\cup A^{2}_{j}$. If
$\{t_1\}\cup A^{1}_{i}$ and $\{t_2\}\cup A^{2}_{j}$ have two common
treatments, then these two ones define one unique set of cardinality $s$ such that
arbitrary two distinct treatments of this set are first associates of the
pseudo-$L_{g}(s)$ association scheme. It follows that $\{t_1\}\cup
A^{1}_{i}=\{t_2\}\cup A^{2}_{j}$, which means two sets of cardinality $s$ are coincident.

Since $t_1$ and $t_2$ are first associates, there exist $t_1\in
A^{2}_{j_1}$ and $t_2\in A^{1}_{i_1}$, here $j_1, i_1\in
\{1,2,\cdots,g\}$. It is very obvious that
$\{t_1, t_2\}\subset(\{t_1\}\cup A^{1}_{i_1})\cap(\{t_2\}\cup
A^{2}_{j_1})$ holds, since $t_1$ and $t_2$ have already
confined the unique set of cardinality $s$, thus we have
$\{t_1\}\cup A^{1}_{i_1}=\{t_2\}\cup A^{2}_{j_1}$.
It follows that $(a)$ is satisfied.

Choose $\{t_1\}\cup A^{1}_{i_2}$, here any $i_2\in
\{1,2,\cdots,g\}$$\setminus$$\{i_1\}$. By $\mathbf{(iii)}$, we have
$A^{1}_{i_1}\neq A^{1}_{i_2}$, i.e., $\{t_1\}\cup A^{1}_{i_1}\neq
\{t_1\}\cup A^{1}_{i_2}$. On the one hand, if there is one set of $\{t_2\}\cup
A^{2}_{1}, \cdots, \{t_2\}\cup A^{2}_{j_1-1}, \{t_2\}\cup
A^{2}_{j_1+1}, \cdots, \{t_2\}\cup A^{2}_{g}$ such that it and
$\{t_1\}\cup A^{1}_{i_2}$ have two common treatments, which define
another unique set of cardinality $s$, then it is equal to $\{t_1\}\cup
A^{1}_{i_2}$. We thereby have $t_2\in A^{1}_{i_2}$, in this case $\{t_1, t_2\}$ belongs to
$(\{t_1\}\cup A^{1}_{i_1})\cap (\{t_1\}\cup A^{1}_{i_2})$. That contradicts
$\{t_1\}\cup A^{1}_{i_1}\neq \{t_1\}\cup A^{1}_{i_2}$. Therefore,
$\{t_1\}\cup A^{1}_{i_2}$ and each of $\{t_2\}\cup A^{2}_{1},
\cdots, \{t_2\}\cup A^{2}_{j_1-1}, \{t_2\}\cup A^{2}_{j_1+1},
\cdots, \{t_2\}\cup A^{2}_{g}$ have at most one common treatment.

On the other hand, it is assumed that $\{t_1\}\cup A^{1}_{i_2}$ and each of them are
all intersecting, thus there are $(g-1)$ treatments in $A^{1}_{i_2}$
that are first associates of $t_2$. Whereas, it is
seen from $\mathbf{(vi)}$ that there are exactly $(g-2)$ treatments in
$A^{1}_{i_2}$ that are first associates of $t_2$ (annotate: $a_i, A_i, A_j$ are
replaced by $t_2, A^{1}_{i_1}, A^{1}_{i_2}$ within $\mathbf{(vi)}$, respectively).
This leads to a contradiction.
It follows that there exists at least one set $\{t_2\}\cup A^{2}_{j_2}$ such that
$(\{t_2\}\cup A^{2}_{j_2})\cap(\{t_1\}\cup A^{1}_{i_2})=\emptyset$
holds, here $j_2\in \{1,2,\cdots,g\}$$\setminus$$\{j_1\}$.

In addition, if there is another set $\{t_2\}\cup A^{2}_{j_3}$
such that $(\{t_2\}\cup A^{2}_{j_3})\cap(\{t_1\}\cup
A^{1}_{i_2})=\emptyset$ holds, $j_3\neq j_2$, then there are at most
$(g-3)$ treatments in $A^{1}_{i_2}$ that are first associates of
$t_2$. This contradicts there are exactly $(g-2)$ treatments in
$A^{1}_{i_2}$.

In sum, we have $(\{t_1\}\cup A^{1}_{i_2})\cap(\{t_2\}\cup
A^{2}_{j_2})=\emptyset$ for one single $j_2\in
\{1,2,\cdots,g\}$$\setminus$$\{j_1\}$ and $(\{t_1\}\cup
A^{1}_{i_2})\cap(\{t_2\}\cup A^{2}_{j'})=\{t_{1j'}\}$ for each
$j'\in \{1,2,\cdots,g\}$$\setminus$$\{j_1, j_2\}$, here $t_{1j'}$ is
commonly first associates of both $t_1$ and $t_2$.
It follows that $(b)$ is satisfied.

Select $\{t_2\}\cup A^{2}_{j_3}$ again, here any $j_3\in
\{1,2,\cdots,g\}$$\setminus$$\{j_1\}$. Because of the symmetry of
$t_2$ and $t_1$, applying the same method, we acquire
$(c)$ is satisfied. The proof is complete. $ \blacksquare $

When the condition $s>(g-1)^2$ is omitted, provided we are able to
choose desired $g$ paralleling classifications
from much more than these $g$ ones of a pseudo-$L_{g}(s)$ association
scheme, so that $\mathbf{(v)}$ can be satisfied inside these $g$
classifications chosen. In this case, with respect to $t$ and
another fixed $a_i$, here $a_i\in A_i$; it is possible that
the set of cardinality $s$ containing $t$ and $a_i$ is non-unique in the
pseudo-$L_{g}(s)$ association scheme.

If $\mathbf{(iii)}$ and $\mathbf{(vi)}$ with $s>(g-1)^2$ of the
pseudo-$L_{g}(s)$ association scheme are satisfied in Theorem 3.6, by
Consequences IV, III and I, then $s^2$ treatments of the pseudo-$L_{g}(s)$
association scheme can form an $L_{g}(s)$ association scheme.
Furthermore, with respect to $t$ and $a_i$, here $a_i\in A_i$; due to
the uniqueness of $\{t\}\cup A_i$, the pseudo-$L_{g}(s)$ association
scheme must be an $L_{g}(s)$ association scheme. In particular, take
$g=3$ and $s>4$, the result still holds; meanwhile $\mathbf{(iii)}$ and
$\mathbf{(vi)}$ in the pseudo-$L_{3}(s)$ association scheme are
transformed into $\mathbf{(I)}$ and $\mathbf{(I^{'})}$ of Lemma 2.21.
With reference to $\mathbf{(I)}$ and $\mathbf{(I^{'})}$ with $s>4$ $\Longrightarrow
\mathbf{(i)}$ and $\mathbf{(ii)}$, we provide another new
proof in the following Theorem 3.8, which is different from
the demonstration of the paper \cite{lwr}, but analogous to the above verification
combining Consequences IV and III for $g=3$.

{\bf Theorem 3.8 } If $\mathbf{(I)}$ and
$\mathbf{(I^{'})}$ with $s>4$ of Lemma 2.21 hold in the pseudo-$L_{3}(s)$
association scheme then so do $\mathbf{(i)}$ and $\mathbf{(ii)}$.

Here $\mathbf{(i)}$ $s^2$ treatments of the pseudo-$L_{3}(s)$
association scheme can be exactly divided into $s$
pairwise disjoint sets $B^v_1, B^v_2, \cdots, B^v_s$ inside the
$v$-th paralleling classification, $v=1,2,3$.

$\mathbf{(ii)}$ With respect to two different classifications
$B^v_1, B^v_2, \cdots, B^v_s$ and $B^u_1, B^u_2, \cdots, B^u_s$, $v,
u\in \{1,2,3\}$; for any $i', j'\in\{1,2,\cdots,s\}$, $B^v_{i'}$ and
$B^u_{j'}$ have exactly one common treatment.

Wherein arbitrary two distinct treatments are first associates of the
pseudo-$L_{3}(s)$ association scheme in each $B^v_{i'}$ and $B^u_{j'}$,
$i', j'=1,2,\cdots,s$, respectively.

{\bf Proof: } Take any treatment $x$, it is known from
$\mathbf{(I)}$ that $3(s-1)$ treatments, which are first associates
of $x$, can be equally divided into three pairwise disjoint sets $Y,
Z, W$.

For $s>4$, with respect to $x$ and $w_i$, here $w_i\in W$, we now
prove that the set of cardinality $s$ containing $x$ and $w_i$ is uniquely $\{x\}\cup W$.

According to $\mathbf{(I^{'})}$, we could exchange $s-1$
treatments$^{,}$ positions separately within $Y$ and $Z$ in order to satisfy the following
conditions: $w_{i}$ and either of
$y_{j}, z_{j}$ are $\left\{\begin{array}{l}
\hbox{first associates, if $i=j$,}\\
\hbox{second associates, if $i\neq j$.}\\
\end{array}\right.$, where \begin{center}$\left\{\begin{array}{l}
\hbox{$Y=\{y_{1},y_{2},\cdots,y_{s-1}\}$,}\\
\hbox{$Z=\{z_{1},z_{2},\cdots,z_{s-1}\}$,}\\
\hbox{$W=\{w_{1},w_{2},\cdots,w_{s-1}\}$.}\\
\end{array}\right.$. \end{center}

When $w_{i}$ and either of $y_{i}, z_{i}$ are first associates,
$i=1,2,\cdots,s-1$. In this case, if $y_{i}$ and $z_{i}$ are second associates,
then any two of $x, w_{i}, y_{i}$ or $x, w_{i}, z_{i}$ are first
associates; if $y_{i}$ and $z_{i}$ are first associates, then any
two of $x, w_{i}, y_{i}, z_{i}$ are first associates.

Next we suppose there exists another collection $\{x\}\cup W'$ of
cardinality $s$ satisfying $w_i\in W'$ and $W'\neq W$, such that arbitrary
two treatments of $\{x\}\cup W'$ are first associates of the
pseudo-$L_{3}(s)$ association scheme. Because of $w_i\in W'$ and $W'\neq
W$, also since $w_{i}$ and either of $y_{j}, z_{j}$ are
second associates ($i\neq j$), hence $\{x\}\cup W'$ of cardinality $s$
with $s> 4$ includes at least
these five treatments $x, w_{i}, y_{i}, w_{j_1}, w_{j_2}$ or $x,
w_{i}, z_{i}, w_{j_1}, w_{j_2}$ or $x, w_{i}, y_{i}, z_{i}, w_{j_1}$,
here $j_1, j_2\in\{1,2,\cdots,s-1\}$$\setminus$$\{i\}, j_1\neq
j_2$.

On the other hand, since $w_{j}$ and either of $y_{i}, z_{i}$ are
second associates ($j\neq i$), thus it will not appear that any two of $x, w_{i},
y_{i}, w_{j_1}, w_{j_2}$ or $x, w_{i}, z_{i}, w_{j_1}, w_{j_2}$ or
$x, w_{i}, y_{i}, z_{i}, w_{j_1}$ are first associates. This leads to
a contradiction.

Hence $x$ and $w_i$ define a unique set $\{x\}\cup W$ of cardinality
$s$ with $s> 4$ such that arbitrary two distinct treatments of $\{x\}\cup W$ are first
associates of the pseudo-$L_{3}(s)$ association scheme. It is also
verified that two sets of cardinality $s$ with $s> 4$ separately comprising
$x, y_i$ and $x, z_i$ are uniquely $\{x\}\cup Y$ and $\{x\}\cup Z$ in a
similar manner, here $y_i\in Y, z_i\in Z$.

Let $B^1_1=\{x\}\cup Y, B^2_1=\{x\}\cup Z, B^3_1=\{x\}\cup W$.
Single out $3(s-1)$ treatments that are first associates of $w_i$,
here $w_i\in W$, according to $\mathbf{(I)}$, it may be assumed that they can be equally
divided into three pairwise disjoint sets $U_i, V_i, W_i, i=1,2,\cdots,s-1$.
With respect to $x$ and $w_i$, because of the uniqueness of $\{x\}\cup W$,
without loss of generality, we put $\{w_i\}\cup W_i=\{x\}\cup W$, where $x\in W_i$.

Using $\{x\}\cup W=\{w_i\}\cup W_i$, we search for common treatments
among $\{x\}\cup Y, \{x\}\cup Z$ and $\{w_i\}\cup U_i, \{w_i\}\cup V_i$. Due to
$p_{11}^{1}=s$, thus the number of treatments which are commonly first
associates of both $x$ and $w_i$ totals $s$. That is to say, there are $s-2$
treatments of $W\cap W_i$ such that they are commonly first associates
of both $x$ and $w_i$, the remaining two treatments are
$y_{i}$ and $z_{i}$. If $y_{i}$ and $z_{i}$ together belong to $U_i$ (or $V_i$),
for $x\in W_i$, then there are two treatments in
$U_i$ (or $V_i$) such that they are first associates of $x$. This contradicts
$\mathbf{(I^{'})}$. Therefore, one of $y_{i}$ and $z_{i}$ must belong to $U_i$, the
other belongs to $V_i$. Although we have fixed $y_i\in Y$ and $z_i\in Z$, but $U_i$ and $V_i$
are unconfined; without loss of generality, we put $z_i\in U_i$
and $y_i\in V_i$.

Since the set of cardinality $s$ with $s> 4$ containing two treatments
is unique, thus we have $(\{w_i\}\cup U_i)\cap (\{x\}\cup Z)=\{z_i\}$ and
$(\{w_i\}\cup V_i)\cap (\{x\}\cup Y)=\{y_i\}$.
Set $B^1_{i+1}=\{w_i\}\cup U_i, B^2_{i+1}=\{w_i\}\cup V_i$ again,
$i=1,2,\cdots,s-1$. Hence we obtain
\begin{center}$\left\{\begin{array}{l}
\hbox{$B^1_1\cap B^1_{i+1}=\emptyset, B^2_1\cap B^2_{i+1}=\emptyset$.}\\
\hbox{$B^1_1\cap B^2_{i+1}=\{y_{i}\}, B^2_1\cap B^1_{i+1}=\{z_{i}\}$.}\\
\end{array}\right.$.\end{center}

Next we choose any $w_{i_1}, w_{j_1}$ from $W$, here $i_1,
j_1\in\{1,2,\cdots,s-1\}, i_1\neq j_1$. According to $\{w_{i_1}\}\cup
W_{i_1}=\{w_{j_1}\}\cup W_{j_1}=\{x\}\cup W$, with the addition of $x\in W_{i_1}\cap
W_{j_1}$, we have $w_{i_1}\in W_{j_1}$ and
$w_{j_1}\in W_{i_1}$. It follows that we
must have $z_{i_1}\in U_{i_1}, y_{i_1}\in V_{i_1}$ and $z_{j_1}\in
U_{j_1}, y_{j_1}\in V_{j_1}$, as well as the corresponding
$B^1_{i_1+1}=\{w_{i_1}\}\cup U_{i_1}, B^2_{i_1+1}=\{w_{i_1}\}\cup
V_{i_1}$ and $B^1_{j_1+1}=\{w_{j_1}\}\cup U_{j_1},
B^2_{j_1+1}=\{w_{j_1}\}\cup V_{j_1}$.

Now we suppose $B^1_{i_1+1}\cap B^1_{j_1+1}\neq \emptyset$, this
indicates $B^1_{i_1+1}\cap B^1_{j_1+1}$ contains the single treatment,
labeled by $t_1$, here $t_1\notin\{w_{i_1}, w_{j_1}\}, w_{i_1}\notin U_{j_1},
w_{j_1}\notin U_{i_1}$. That will eventually lead to a contradiction.

Because of $B^1_{i_1+1}\cap B^1_{j_1+1}\neq \emptyset$,
we must have $t_1\in U_{i_1}\cap U_{j_1}$.
It is known from $\mathbf{(I)}$ that $3(s-1)$ treatments, which are first associates
of $t_1$, can be equally divided into three pairwise disjoint sets
$A_1, A_2, A_3$. Let $B^1_{i_1+1}=\{t_1\}\cup A_1$ and $B^1_{j_1+1}=\{t_1\}\cup
A_2$, noting $B^1_{i_1+1}, B^1_{j_1+1}$
and $B^1_1, B^2_1, B^3_1$, we found:
\begin{center}$\left\{\begin{array}{cc}
\hbox{$B^1_{i_1+1}\cap B^1_1=\emptyset,$}
&\hbox{$B^1_{j_1+1}\cap B^1_1=\emptyset,$}\\
\hbox{$B^1_{i_1+1}\cap B^2_1=\{z_{i_1}\},$}
&\hbox{$B^1_{j_1+1}\cap B^2_1=\{z_{j_1}\},$}\\
\hbox{$B^1_{i_1+1}\cap B^3_1=\{w_{i_1}\}.$}
&\hbox{$B^1_{j_1+1}\cap B^3_1=\{w_{j_1}\}.$}\\
\end{array}\right.$.\end{center}
We obviously see $z_{i_1}\neq z_{j_1}$ and
$t_1\notin\{z_{i_1}, z_{j_1}\}$, otherwise,
there are two treatments $w_{i_1}$ and $w_{j_1}$ in $W$ such that
they are first associates of
$z_{i_1}=z_{j_1}$ in $Z$. It is gradually realized that $t_1$ does not
belong to anyone of $Y, Z, W$, thus $t_1$ and $x$ are second associates.
Due to $p_{11}^{2}=6$, thus the number of treatments which are commonly
first associates of both $t_1$ and $x$ totals six.

It is assumed again that $\{t_1\}\cup A_3$ and $B^2_1$ have one common treatment,
denoted by $t_2$. We obviously have $t_2\notin\{z_{i_1}, z_{j_1}\}$.
Take $t_2$ as an example, since $t_1$ belongs to someone of three pairwise disjoint
sets which are first associates of $t_2$, there are two treatments $z_{i_1}$ and
$z_{j_1}$ in another of these three sets such that they are first associates of
$t_1$. This contradicts $\mathbf{(I^{'})}$. It follows that $\{t_1\}\cup A_3$
and $B^2_1$ have no common treatments, similarly, so do $\{t_1\}\cup
A_3$ and $B^3_1$.

Subtracting $z_{i_1}, z_{j_1}, w_{i_1}, w_{j_1}$ from six treatments,
we figure out $\{t_1\}\cup A_3$
and $B^1_1$ (=$\{x\}\cup Y$) have exactly two common treatments,
which define a unique set of cardinality $s$ ($s>4$) such that arbitrary two
distinct treatments of this set are first associates of the
pseudo-$L_{3}(s)$ association scheme. It goes without saying that $\{t_1\}\cup
A_3=\{x\}\cup Y$ holds, which indicates $t_1$ and $x$ are first associates.
This contradicts none of $Y, Z, W$ contains $t_1$.

The verification of $ B^2_{i_1+1}\cap B^2_{j_1+1}=\emptyset$ is similar to that of
$B^1_{i_1+1}\cap B^1_{j_1+1}=\emptyset$. It follows that $U_{i_1}\cap U_{j_1}=\emptyset$
and $V_{i_1}\cap V_{j_1}=\emptyset$ hold. Besides $s-2$ treatments of
$W_{i_1}\cap W_{j_1}$, which are commonly first
associates of both $w_{i_1}$ and $w_{j_1}$, by $\mathbf{(I^{'})}$,
there are the remaining two treatments such that one of them must
belong to $U_{i_1}$, the other belongs to $V_{i_1}$. Likewise,
these two treatments separately belong to $V_{j_1}$ and $U_{j_1}$.
Therefore, labelling these two treatments by $w_{i_1j_1}$ and $w_{j_1i_1}$,
we learn $U_{i_1}\cap V_{j_1}=\{w_{i_1j_1}\}$ and
$U_{j_1}\cap V_{i_1}=\{w_{j_1i_1}\}$, here $w_{i_1j_1}$
is obviously first associates of both $z_{i_1}$ and $y_{j_1}$; $w_{j_1i_1}$ is obviously
first associates of both $z_{j_1}$ and $y_{i_1}$. In a word, we have
\begin{center}$\left\{\begin{array}{l}
\hbox{$B^1_{i_1+1}\cap B^1_{j_1+1}=\emptyset, B^2_{i_1+1}\cap B^2_{j_1+1}=\emptyset$.}\\
\hbox{$B^1_{i_1+1}\cap B^2_{j_1+1}=\{w_{i_1j_1}\}, B^1_{j_1+1}\cap B^2_{i_1+1}=\{w_{j_1i_1}\}$.}\\
\end{array}\right.$.\end{center}

Because of $B^1_1\cap B^1_{i_1+1}=\emptyset, B^2_1\cap
B^2_{j_1+1}=\emptyset$ and $B^1_{i_1+1}\cap B^1_{j_1+1}=\emptyset,
B^2_{i_1+1}\cap B^2_{j_1+1}=\emptyset$, foresaid
$x, w_{i_1}, w_{j_1}, w_{i_1j_1}, z_{i_1}, y_{j_1}, w_{j_1i_1}, z_{j_1}, y_{i_1}$
are arranged onto an $s\times s$ matrix
\begin{center}$T=\left(
\begin{smallmatrix}
x&y_1&y_2&y_3&\cdots&y_{s-1}\\
z_1&w_1&w_{12}&w_{13}&\cdots&w_{1(s-1)}\\
z_2&w_{21}&w_2&w_{23}&\cdots&w_{2(s-1)}\\
z_3&w_{31}&w_{32}&w_3&\cdots&w_{3(s-1)}\\
\vdots&\vdots&\vdots&\vdots&\vdots&\vdots\\
z_{s-1}&w_{(s-1)1}&w_{(s-1)2}&w_{(s-1)3}&\cdots&w_{s-1}\\
\end{smallmatrix}
\right)$.
\end{center}

It is obviously seen that $s^2$ treatments of the pseudo-$L_{3}(s)$
association scheme are filled $T$ with,
and the $i'$-th row and the $j'$-th column of $T$
are exactly $B^1_{i'}$ and $B^2_{j'}$, respectively, moreover, $B^1_{i'}$ and
$B^2_{j'}$ have exactly one common treatment, $i', j'=1,2,\cdots,s$.

Labelling $z_0=x$ to state the following argument conveniently, finally we
select any $z_{i_1}, z_{j_1}\in B^2_1, z_{i_1}\neq z_{j_1}$, where
$B^2_1=\{z_{0},z_{1},\cdots,z_{s-1}\}$. Observe $T$, it is convinced from $\mathbf{(I)}$
that there exist $B^3_{i_1+1}$ and $B^3_{j_1+1}$ such that $B^2_1\cap
B^1_{i_1+1}\cap B^3_{i_1+1}=\{z_{i_1}\}$ and $B^2_1\cap
B^1_{j_1+1}\cap B^3_{j_1+1}=\{z_{j_1}\}$ hold, here $i_1, j_1=0,1,\cdots,s-1, i_1\neq j_1$.
The next statement is how to construct and explicate $B^3_{i'}, i'=i_1+1$.

Take each $z_{j_1}$ as an instance, since $z_{i_1}$ belongs to
$B^2_1$$\setminus$$\{z_{j_1}\}$ of three pairwise disjoint sets which
are first associates of $z_{j_1}$, also by $\mathbf{(I^{'})}$,
there exists exactly a sole treatment in $B^1_{j_1+1}$$\setminus$$\{z_{j_1}\}$ of these three
sets such that it is first associates of $z_{i_1}$. Denoting the sole treatment as
$b^{i_1+1}_{j_1+1}$, because of $b^{i_1+1}_{j_1+1}\notin B^2_1\cup B^1_{i_1+1}$,
thus we have $b^{i_1+1}_{j_1+1}\in
B^1_{j_1+1}\cap B^3_{i_1+1}$. Put $b^{i_1+1}_{i_1+1}=z_{i_1}$, thus we obtain
$B^3_{i_1+1}=\{b^{i_1+1}_{1}, \cdots, b^{i_1+1}_{i_1+1}, \cdots, b^{i_1+1}_{s}\}=B^3_{i'}$.
This implies that there is exactly one treatment
of $B^3_{i_1+1}$ on every row of $T$.

Noticing $T$, we have
$B^1_{i_1+1}=\{z_{i_1}, w_{i_11}, \cdots, w_{i_1}, \cdots, w_{i_1(s-1)}\}$.
It is very obvious that $z_{i_1}$ at the first column of $T$ belongs to $B^3_{i_1+1}$.
Due to $B^1_{i_1+1}\cap B^2_{j_1+1}=\{w_{i_1j_1}\}$. Fixing $B^1_{i_1+1}$
and varying $B^2_{j_1+1}, j_1=1,2,\cdots,s-1$, we set
$B^1_{i_1+1}\cap B^2_{j_1+1}=\{w_{i_1}\}$ if $i_1=j_1$,
and see $w_{i_1j_1}$ is located on the $(j_1+1)$-th column of $T$.
We choose $w_{i_1j_1}$ as a sample, for $z_{i_1}\in B^1_{i_1+1}$$\setminus$$\{w_{i_1j_1}\}$,
where $B^1_{i_1+1}$$\setminus$$\{w_{i_1j_1}\}$ is someone of three pairwise disjoint
sets which are first associates of $w_{i_1j_1}$,
there exists exactly a sole treatment in $B^2_{j_1+1}$$\setminus$$\{w_{i_1j_1}\}$
such that it is first associates of $z_{i_1}$. Because of $B^2_1\cap
B^2_{j_1+1}=\emptyset$, thus the sole treatment
is not contained in $B^2_1\cup B^1_{i_1+1}$, but belongs to
$B^2_{j_1+1}\cap B^3_{i_1+1}$. It follows that
there is exactly one treatment of $B^3_{i_1+1}$ on every column of $T$.

For $i', j'=1,2,\cdots,s$, it goes without saying that $B^3_{i'}$ and either of
$B^1_{j'}$ and $B^2_{j'}$ have exactly one common treatment, in addition,
$B^1_{i'}$ and $B^2_{j'}$ have exactly one common treatment. That is,
$\mathbf{(ii)}$ holds.

Because of $B^2_1\cap B^1_{i_1+1}\cap B^3_{i_1+1}=\{z_{i_1}\}$ and $B^2_1\cap
B^1_{j_1+1}\cap B^3_{j_1+1}=\{z_{j_1}\}$, here $z_{i_1}\neq z_{j_1}$, next we prove
$B^3_{i_1+1}\cap B^3_{j_1+1}=\emptyset$ by contradiction.

We assume $B^3_{i_1+1}\cap B^3_{j_1+1}\neq \emptyset$ again. Due to
$B^3_{i_1+1}\cap B^1_{i_1+1}=\{z_{i_1}\}, B^3_{j_1+1}\cap
B^1_{j_1+1}=\{z_{j_1}\}$ and $B^1_{i_1+1}\cap
B^1_{j_1+1}=\emptyset$, therefore, neither of $B^3_{i_1+1}\cap
B^1_{j_1+1}$ and $B^3_{j_1+1}\cap B^1_{i_1+1}$ is contained in
$B^3_{i_1+1}\cap B^3_{j_1+1}$. After $B^3_{i_1+1}\cap B^1_{j_1+1}$
and $B^3_{j_1+1}\cap B^1_{i_1+1}$ are added to $s-2$ treatments of
$B^2_1$$\setminus$$\{z_{i_1}, z_{j_1}\}$, $B^3_{i_1+1}\cap B^3_{j_1+1}$
is superadded. We aggregate the number of treatments that are
commonly first associates of both $z_{i_1}$
and $z_{j_1}$, it is greater than $s$. This contradicts
$p_{11}^{1}=s$.

In sum, we learn from $T$ that $s^2$ treatments of the
pseudo-$L_{3}(s)$ association scheme can be divided into $s$
pairwise disjoint sets $B^v_1, B^v_2, \cdots, B^v_s$, $v=1,2,3$.
That is, $\mathbf{(i)}$ holds. Now the proof is
finished. $ \blacksquare $

If $\mathbf{(I)}$ and $\mathbf{(I^{'})}$ with $s>4$ of the
pseudo-$L_{3}(s)$ association scheme are satisfied, using Theorem 3.8
and Consequence I of Theorem 3.6, then its $s^2$ treatments can form an
$L_{3}(s)$ association scheme. We consider three pairwise disjoint
sets $X, Y, Z$ of cardinality $s-1$ which are first associates of $x$,
because of the uniqueness of $\{x\}\cup Y, \{x\}\cup Z, \{x\}\cup
W$, thus the pseudo-$L_{3}(s)$ association scheme must be an
$L_{3}(s)$ association scheme.

It is assumed that there exists an $L_{g}(s)$ association scheme,
using Definition 2.16, a pseudo-$L^{*}_{[s+1-g]}(s)$ association scheme
is induced by the $L_{g}(s)$ association scheme. For instance,
if any two distinct treatments are first associates in the
pseudo-$L^{*}_{[s+1-g]}(s)$ association scheme,
then they are second associates in the $L_{g}(s)$
association scheme; and vice versa. Utilizing Lemma 2.20, we know the
pseudo-$L^{*}_{[s+1-g]}(s)$ association scheme is a pseudo-$L_{s+1-g}(s)$ association scheme.
For $g\geq 3, s\geq g+2$, after $g$ is replaced with $s+1-g$ in Theorem 3.6,
we therefore have the following

$\mathbf{(i)}$ $s^2$ treatments of the pseudo-$L^{*}_{[s+1-g]}(s)$
association scheme can be exactly divided into $s$
pairwise disjoint sets $B^{*i}_1, B^{*i}_2, \cdots, B^{*i}_s$ inside
the $i$-th paralleling classification, $i=1,2,\cdots,s+1-g$.

$\mathbf{(ii)}$ With respect to any two different classifications
$B^{*i}_1, \cdots, B^{*i}_s$ and $B^{*j}_1, \cdots, B^{*j}_s$,
$i, j\in \{1,2,\cdots,s+1-g\}$; for any $m,
n\in\{1,\cdots,s\}$, $B^{*i}_m$ and $B^{*j}_n$ have exactly one
common treatment.

$\mathbf{(iii)}$ For $(s+1-g)(s-1)$ treatments, which are first associates
of any treatment $t$ in the pseudo-$L^{*}_{[s+1-g]}(s)$
association scheme, they can be equally divided into $(s+1-g)$
pairwise disjoint sets $A^{*}_1, A^{*}_2, \cdots, A^{*}_{s+1-g}$.

$\mathbf{(iv)}$ With respect to $B^{*i}_{m_i}$ containing the treatment $t$,
where $B^{*i}_{m_i}\subset\{B^{*i}_1, B^{*i}_2, \cdots, B^{*i}_s\}$,
$i=1,2,\cdots,s+1-g$; there exactly exists a one to one mapping of
$B^{*i}_{m_i}$ onto $A^{*}_{k_i}$ satisfying $B^{*i}_{m_i}=\{t\}\cup A^{*}_{k_i},
k_i\in\{1,2,\cdots,s+1-g\}$.

$\mathbf{(v)}$ With respect to arbitrary two distinct treatments $t^{*}_1$ and $t^{*}_2$,
they are first associates of the pseudo-$L^{*}_{[s+1-g]}(s)$ association scheme, it is assumed that
$(s+1-g)(s-1)$ treatments, which are first associates of $t^{*}_v$,
can be equally divided into $(s+1-g)$ pairwise disjoint sets $A^{*v}_1,
A^{*v}_2, \cdots, A^{*v}_{s+1-g}$, $v=1,2$. If
$$\left\{
\begin{array}{lll}
(a).\ \ \{t_1^{*}\}\cup A^{*1}_{i_1}=\{t_2^{*}\}\cup
A^{*2}_{j_1}, \ \ for \ t_1^{*}\in A^{*2}_{j_1}, t_2^{*}\in A^{*1}_{i_1};\ moreover,\
\ t_1^{*}\ and \ t_2^{*}\ define\ the\ set\ of\\ \ \ \ \ \ \ \ cardinality\ s,\ in\ which\ arbitrary\ two\
distinct\ treatments\ are\ first\ associates\ of\ the\\ \ \ \ \ \ \ \ pseudo$-$L^{*}_{[s+1-g]}(s)\ association\ scheme.\\
(b).\ \ For\ any\ i_2\in \{1,2,\cdots,s+1-g\}\setminus\{i_1\},\\
\ \ \ \ \ \ (\{t_1^{*}\}\cup A^{*1}_{i_2})\cap(\{t_2^{*}\}\cup A^{*2}_{j_2})=\emptyset, \ \ a\ single\ j_2\in \{1,2,\cdots,s+1-g\}\setminus\{j_1\}.\\
\ \ \ \ \ \ (\{t_1^{*}\}\cup A^{*1}_{i_2})\cap(\{t_2^{*}\}\cup
A^{*2}_{j'})=\{t^{*}_{1j'}\},\ \ each\ j'\in \{1,2,\cdots,s+1-g\}\setminus\{j_1, j_2\}.\\
(c). \ \ For\ any\ j_3\in \{1,2,\cdots,s+1-g\}\setminus\{j_1\},\\
\ \ \ \ \ \ (\{t_2^{*}\}\cup A^{*2}_{j_3})\cap(\{t_1^{*}\}\cup A^{*1}_{i_3})=\emptyset,\ \ a\ single\ i_3\in \{1,2,\cdots,s+1-g\}\setminus\{i_1\}.\\
\ \ \ \ \ \ (\{t_2^{*}\}\cup A^{*2}_{j_3})\cap(\{t_1^{*}\}\cup A^{*1}_{i'})=\{t_{2i'}^{*}\},
\ each\ i'\in \{1,2,\cdots,s+1-g\}\setminus\{i_1, i_3\}.\\
\ \ Here,\ t_{1j'}^{*}\ and\ t_{2i'}^{*}\ are\ commonly\ first\ associates \ of\ both \ t_1^{*}\ and\ t_2^{*}.
\end{array}
\right. $$
are satisfied.

$\mathbf{(vi)}$ For every $a^{*}_i\in A^{*}_i$,
$i=1,2,\cdots,s+1-g$; there are exactly $[(s+1-g)-2]$ treatments in
$A^{*}_j$ such that they are first associates of $a^{*}_i$, any
$j\in \{1,2,\cdots,s+1-g\}$$\setminus$$\{i\}$.

Then the following four consequences could be obtained:

I. If $\mathbf{(i)}$ and $\mathbf{(ii)}$ hold then $s^2$ treatments of the pseudo-$L^{*}_{[s+1-g]}(s)$
association scheme can form an $L_{s+1-g}(s)$ association scheme.

II. $\mathbf{(i)}$, $\mathbf{(iii)}$ and
$\mathbf{(iv)}$ are equivalent to $\mathbf{(i)}$ and $\mathbf{(ii)}$
in the pseudo-$L^{*}_{[s+1-g]}(s)$ association scheme.

III. $\mathbf{(v)}$ is equivalent to $\mathbf{(i)}$ and
$\mathbf{(ii)}$ in the pseudo-$L^{*}_{[s+1-g]}(s)$ association scheme.

IV. If $\mathbf{(iii)}$ and
$\mathbf{(vi)}$ with $s>(s-g)^2$ are both satisfied in the pseudo-$L^{*}_{[s+1-g]}(s)$
association scheme, then we have

(1) the set of cardinality $s$ containing $t$ and $a^{*}_i$ is
uniquely $\{t\}\cup A^{*}_{i}$, where $a^{*}_i\in A^{*}_i$; and

(2) $\mathbf{(v)}$ holds.

Where arbitrary two distinct treatments are first associates of the
pseudo-$L^{*}_{[s+1-g]}(s)$ association scheme within each
$\{t\}\cup A^{*}_{i}$, $\{t^{*}_1\}\cup A^{*1}_{i}$,
$\{t^{*}_2\}\cup A^{*2}_{i}, i=1,2,\cdots,s+1-g$; each $B^{*i}_m,
B^{*j}_n, m, n=1,2,\cdots,s$, respectively.

It has been known that $s^2$
treatments of an $L_{g}(s)$ association scheme may be used to
arrange a net $N$ of order $s$, degree $g$. As far as one given
$L_{g}(s)$ association scheme is concerned, it can been assumed
from Definition 2.11 that its $s^2$ treatments are set forth in an $s\times s$ array and
its $POL(s,g-2)$ is $\{\mathbf{L_{1}, L_{2}, \cdots,
L_{g-2}}\}$. Next we review how to arrange a net $N$ with
$s^2$ treatments of the $L_{g}(s)$ association scheme.
Every treatment of it can be
regarded as each point of $N$. Firstly every row of the array
is called as each line of the first parallel class of $N$, every
column of the array is called as each line of the second
parallel class of $N$. Finally $s$ different treatments of
the array, which altogether correspond to the same symbol of
$\mathbf{L_k}$, lie on the same "line" of the ($\mathbf{k}+2$)-th
parallel class of $N$, $\mathbf{k}=1,2,\cdots,g-2$. It is
obviously verified that Statements (I)-(IV) of Definition 2.12 can be
satisfied, hence we utilize $s^2$ treatments of the $L_{g}(s)$ association
scheme to arrange the net $N$. Conversely, we
can obtain any net of order $s$, degree $g$ ($g\geq 1$) in the
manner indicated, usually in many ways.

We suppose that there exists the above net $N$ of order $s$, degree $g$, and define that
two distinct points of $N$ are first associates if and only if they are
joined in $N$, and second associates otherwise. According to the definition, it is obtained that
$s^2$ points of $N$ can form a pseudo-$L_{g}(s)$ association scheme.

Using Definition 2.17, a pseudo-net-$N^{*}$ of order $s$, degree $s+1-g$ is induced
by the above net $N$. For example, if any two distinct points are joined in pseudo-net-$N^{*}$,
then they are not joined in $N$; and vice versa. It is defined that two distinct points of
pseudo-net-$N^{*}$ are first associates if and only if they are
joined in pseudo-net-$N^{*}$, and second associates otherwise. Under this definition,
it is known from Lemma 2.20 that $s^2$ points of pseudo-net-$N^{*}$ can form a
pseudo-$L^{*}_{s+1-g}(s)$ association scheme, which is surely a
pseudo-$L_{s+1-g}(s)$ association scheme. It is obviously seen that the "line" of
pseudo-net-$N^{*}$ induced is exactly like the transversal of $N$,
after any treatment $t$ in the pseudo-$L^{*}_{[s+1-g]}(s)$
association scheme is replaced with any point $P$ in pseudo-net-$N^{*}$,
for $g\geq 3, s\geq g+2$, thus we obtain:

$\mathbf{(i)}$ $s^2$ points of the pseudo-net-$N^{*}$ be exactly
distributed into $s$ pairwise parallel lines $\beta^{*i}_1, \beta^{*i}_2,\\ \cdots,
\beta^{*i}_s$ inside the $i$-th parallel class, $i=1,2,\cdots,s+1-g$.

$\mathbf{(ii)}$ With respect to any two different parallel classes
$\beta^{*i}_1, \cdots, \beta^{*i}_s$ and
$\beta^{*j}_1, \cdots, \beta^{*j}_s$, $i, j\in
\{1,2,\cdots,s+1-g\}$; for any $m, n\in\{1,\cdots,s\}$,
$\beta^{*i}_m$ and $\beta^{*j}_n$ have exactly one common point.

$\mathbf{(iii)}$ For $(s+1-g)(s-1)$ points, which are all joined
to any point $P$ in pseudo-net-$N^{*}$, they can be equally
distributed into $(s+1-g)$ distinct lines $\alpha^{*}_1, \alpha^{*}_2, \cdots,
\alpha^{*}_{s+1-g}$, here $\alpha^{*}_{i}\cap \alpha^{*}_{j}=\{P\}$,
any $i, j\in\{1,2,\cdots,s+1-g\}, i\neq j$.

$\mathbf{(iv)}$ With respect to $\beta^{*i}_{m_i}$ containing the point $P$,
where $\beta^{*i}_{m_i}\subset\{\beta^{*i}_1, \beta^{*i}_2, \cdots,
\beta^{*i}_s\}$, $i=1,2,\cdots,s+1-g$; there exactly exists a one to one mapping of
$\beta^{*i}_{m_i}$ onto $\alpha^{*}_{k_i}$ satisfying $\beta^{*i}_{m_i}=\alpha^{*}_{k_i},
k_i\in\{1,2,\cdots,s+1-g\}$.

$\mathbf{(v)}$ With respect to arbitrary two distinct points $Q^{*}_1$ and $Q^{*}_2$,
they are joined in pseudo-net-$N^{*}$, it is assumed that $(s+1-g)(s-1)$
points, which are all joined to $Q^{*}_v$ in pseudo-net-$N^{*}$, can be equally
distributed into $(s+1-g)$ distinct lines $\alpha^{*v}_1,
\alpha^{*v}_2, \cdots, \alpha^{*v}_{s+1-g}$, $v=1,2$, here
$\alpha^{*v}_{i}\cap \alpha^{*v}_{j}=\{Q^{*}_v\}$, for any $i,
j\in\{1,2,\cdots,s+1-g\}, i\neq j$. If
$$\left\{
\begin{array}{lll}
(a).\ \  \alpha^{*1}_{i_1}=\alpha^{*2}_{j_1},\ for\ \ Q^{*}_1\ on \ \ \alpha^{*2}_{j_1},\ Q^{*}_2\ \ on\ \ \alpha^{*1}_{i_1};
\ moreover, Q^{*}_1\ and Q^{*}_2 \ define\ a\ ``line"\ having\\
\ \ \ \ \ \ s\ distinct\ points,
\ among \ which\ arbitrary\ two\  ones\ are\ joined\ in\ pseudo$-$net$-$N^{*}.\\
(b).\ \ For\ any\ i_2\in \{1,2,\cdots,s+1-g\}\setminus\{i_1\},\\
\ \ \ \ \ \ \alpha^{*1}_{i_2}\cap \alpha^{*2}_{j_2}=\emptyset, \ \ a\ single\ j_2\in \{1,2,\cdots,s+1-g\}\setminus\{j_1\}.\\
\ \ \ \ \ \ \alpha^{*1}_{i_2}\cap \alpha^{*2}_{j'}=\{Q^{*}_{1j'}\},\ \ each\ j'\in \{1,2,\cdots,s+1-g\}\setminus\{j_1, j_2\}.\\
(c). \ \ For\ any\ j_3\in \{1,2,\cdots,s+1-g\}\setminus\{j_1\},\\
\ \ \ \ \ \ \alpha^{*2}_{j_3}\cap \alpha^{*1}_{i_3}=\emptyset,\ \ a\ single\ i_3\in \{1,2,\cdots,s+1-g\}\setminus\{i_1\}.\\
\ \ \ \ \ \ \alpha^{*2}_{j_3}\cap \alpha^{*1}_{i'}=\{Q^{*}_{2i'}\},
\ each\ i'\in \{1,2,\cdots,s+1-g\}\setminus\{i_1, i_3\}.\\
\ \ Here \ Q^{*}_{1j'}\ and\ Q^{*}_{2i'}\ are\ commonly\ joined\ to\ both\ Q^{*}_1\ and\ Q^{*}_2\ in\ pseudo$-$net$-$N^{*}.
\end{array}
\right. $$
are satisfied.

$\mathbf{(vi)}$ For every point $P^{*}_i$, it lies on
$\alpha^{*}_i$ in addition to $P$, $i=1,2,\cdots,s+1-g$; there are
exactly $[(s+1-g)-2]$ points on $\alpha^{*}_j$ in addition to $P$
such that they are all joined to $P^{*}_i$ in pseudo-net-$N^{*}$, any $j \in
\{1,2,\cdots,s+1-g\}$$\setminus$$\{i\}$.

Then the following four consequences could be obtained:

I. If $\mathbf{(i)}$ and $\mathbf{(ii)}$ hold then $s^2$ points
of the pseudo-net-$N^{*}$ can be arranged a net of order $s$, degree $s+1-g$.

II. $\mathbf{(i)}$, $\mathbf{(iii)}$ and
$\mathbf{(iv)}$ are equivalent to $\mathbf{(i)}$ and $\mathbf{(ii)}$
in pseudo-net-$N^{*}$.

III. $\mathbf{(v)}$ is equivalent to $\mathbf{(i)}$ and
$\mathbf{(ii)}$ in pseudo-net-$N^{*}$.

IV. If $\mathbf{(iii)}$ and
$\mathbf{(vi)}$ with $s>(s-g)^2$ are both satisfied in pseudo-net-$N^{*}$, then we have

(1) possessing $s$ distinct points, the line which $P$ and $P^{*}_i$ lie on is
uniquely $\alpha^{*}_i$; and

(2) $\mathbf{(v)}$ holds.

Where each $\alpha^{*}_i, \alpha^{*1}_i, \alpha^{*2}_i$
and each $\beta^{*i}_m, \beta^{*j}_n$ are the transversals of $N$,
$i=1,2,\cdots,s+1-g; m, n=1,2,\cdots,s$.

Set $w=g-2$, we may utilize $s^2$ treatments of $\mathbf{T}$ and
$\mathbf{L_{1}, \cdots, L_{w}}$ of a $POL(s,w)$ in Definition 2.7 to construct an $L_{w+2}(s)$
association scheme here $w\geq 1, s\geq w+4$. Next $s^2$ treatments of the
$L_{w+2}(s)$ association scheme are used to arrange a net $N$ of
order $s$, degree $w+2$ without doubt, every treatment of it can be regarded as each
point of $N$. Due to $s+1-g=s-1-w$, applying Definitions 2.16 and 2.17,
we acquire one inducing pseudo-$L^{*}_{[s-1-w]}(s)$ association
scheme and one inducing pseudo-net-$N^{*}$ of order $s$, degree $s-1-w$.
It has already been known that the inducing pseudo-$L^{*}_{[s-1-w]}(s)$
association scheme is a pseudo-$L_{s-1-w}(s)$ association scheme. With respect to two
distinct points of the inducing pseudo-net-$N^{*}$, when we define they
are first associates if and only if they are joined in
pseudo-net-$N^{*}$, and second associates otherwise; thus $s^2$ points of
pseudo-net-$N^{*}$ also can form a pseudo-$L_{s-1-w}(s)$
association scheme.

It is always assumed there exist common transversals of
$\mathbf{L_{1}, \cdots, L_{w}}$ containing any treatment $t$,
we therefore have

{\bf Theorem 3.9 } Let $POL(s,w)=\{\mathbf{L_{1}, \cdots, L_{w}}\}$.
Taking any treatment $t$ and any point $P$, we have \\
(1) for every $a^{*}_i\in A^{*}_i$, there are exactly $[(s-1-w)-2]$ treatments in
$A^{*}_j$ such that they are first associates of $a^{*}_i$ in a
pseudo-$L^{*}_{[s-1-w]}(s)$ association scheme.\\
(2) for every point $P^{*}_i$, it lies on
$\alpha^{*}_i$ in addition to $P$; there are
exactly $[(s-1-w)-2]$ points on $\alpha^{*}_j$ in addition to $P$
such that they are all joined to $P^{*}_i$ in pseudo-net-$N^{*}$.

Here $w\geq 1, s\geq w+4$, $A^{*}_i\cap A^{*}_j=\emptyset$,
$\{t\}\cup A^{*}_i$ and $\{t\}\cup A^{*}_j$
are the common transversals of $\mathbf{L_{1}, \cdots, L_{w}}$;
$\alpha^{*}_i\cap \alpha^{*}_j=\{P\}$,
$\alpha^{*}_i$ and $\alpha^{*}_j$ are the transversals of $N$, $i\neq j$.

{\bf Proof: } We consider every $a^{*}_i\in A^{*}_i$. Since $a^{*}_i$ and $t$
first associates in a pseudo-$L^{*}_{[s-1-w]}(s)$ association
scheme, thus they appear neither in the same row nor in the same
column of $\mathbf{T}$, but correspond to two distinct symbols of
each $\mathbf{L_{k}}, k=1,\cdots,w$.

Since $\{t\}\cup A^{*}_j$ is exactly the common transversal of
$\mathbf{L_{1}, \cdots, L_{w}}$, we understand that there is exactly
one treatment of $\{t\}\cup A^{*}_j$ in each row and column of
$\mathbf{T}$, moreover, its $s$ treatments correspond to $s$
distinct symbols in each $\mathbf{L_{k}}, k=1,\cdots,w$.

Therefore, there are just two treatments $t_{w+1}, t_{w+2}$ in
$A^{*}_j$ such that $t_{w+1}, t_{w+2}$ and $a^{*}_i$ lie on
the same row and the same column of $\mathbf{T}$, respectively.
There is one unique treatment $t_{k}$ in $A^{*}_j$ such that $t_{k}$ and $a^{*}_i$
correspond to the same symbol of each $\mathbf{L_{k}},
k=1,\cdots,w$. It follows that these aggregating $w+2$ treatments of $A^{*}_j$
and $a^{*}_i$ are second associates in the pseudo-$L^{*}_{[s-1-w]}(s)$
association scheme.

The next step is to illustrate any two of $t_{1}, \cdots, t_{w},
t_{w+1}, t_{w+2}$ are different. By Definition 2.2, it is very obvious that $t_{k}\notin
\{t_{w+1},t_{w+2}\}$ holds, here $k\in\{1,\cdots,w\}$.
Because of pairwise orthogonality of $\mathbf{L_{1}, \cdots,
L_{w}}$, we must have $t_{k_1}\neq t_{k_2}$, here any $k_1,
k_2\in\{1,\cdots,w\}, k_1\neq k_2$.

Subtracting these $w+2$ treatments, there are the remaining
$[(s-1-w)-2]$ ones in $A^{*}_j$ such that every one of them and
$a^{*}_i$ correspond to distinct symbols of each $\mathbf{L_{k}}$,
which means they are first associates of $a^{*}_i$ in the
pseudo-$L^{*}_{[s-1-w]}(s)$ association scheme. Hence (1) holds.

After $s^2$ treatments of the $L_{w+2}(s)$ association scheme are
used to arrange a net $N$ of order $s$, degree $w+2$, any treatment $t$ can be
replaced with any point $P$, it follows that (2) holds. The proof is
finished. $ \blacksquare $

{\bf Remark 3.10. } With respect to $(s-1-w)(s-1)$ treatments, which are first associates
of some treatment $t_1$ in the pseudo-$L^{*}_{[s-1-w]}(s)$ association scheme,
either they cannot be equally divided into $(s-1-w)$ pairwise
disjoint sets, or there are more than one partitions of them such that
they can be equally divided into $(s-1-w)$ pairwise
disjoint sets, but there possibly exists $A^{*1}_i\cap A^{*1}_j=\emptyset$
such that $\{t_1\}\cup A^{*1}_i$ and $\{t_1\}\cup A^{*1}_j$
are the common transversals of $\mathbf{L_{1}, \cdots, L_{w}}$
whichever case is satisfied. Similarly, there probably exists
$\alpha^{*}_i\cap \alpha^{*}_j=\{P\}$
such that $\alpha^{*}_i$ and $\alpha^{*}_j$ are the transversals of $N$.

Even if the inducing pseudo-$L^{*}_{[s-1-w]}(s)$
association scheme satisfies conclusion (1) of Theorem 3.9. However,
after $g$ in Theorem 3.5 is replaced by $s-1-w$, five corresponding conditions
of the pseudo-$L^{*}_{[s-1-w]}(s)$ association scheme
may be obtained. So long as any one of these five conditions is satisfied,
the pseudo-$L^{*}_{[s-1-w]}(s)$ association scheme must not
be an $L_{s-1-w}(s)$ association scheme. Similarly, the inducing
pseudo-net-$N^{*}$ of order $s$, degree $s-1-w$ satisfies conclusion (2) of Theorem 3.9.
Although the inducing pseudo-net-$N^{*}$ is not a net of order $s$,
degree $s-1-w$, either, but its $s^2$ points are possibly arranged this net.

After $g$ in Theorem 3.6 is replaced with $s-1-w$, these matching conditions
and conclusions of the inducing pseudo-$L^{*}_{[s-1-w]}(s)$ association scheme
still hold. Researching the extension of a $POL(s,w)$, now we consider

{\bf Theorem 3.11 } If $s^2$ treatments of the inducing
pseudo-$L^{*}_{[s-1-w]}(s)$ association scheme can form an
$L_{s-1-w}(s)$ association scheme, then \\
(1) the net $N$ of order $s$, degree $w+2$ has at least one complementary net;\\
(2) the $POL(s,w)$ may be extended to a $POL(s,s-1)$.

{\bf Proof: } (1) Since $s^2$ treatments of the
$L_{w+2}(s)$ association scheme are arranged a net $N$ of
order $s$, degree $w+2$. Also since $s^2$ treatments of the inducing
pseudo-$L^{*}_{[s-1-w]}(s)$ association scheme can form an
$L_{s-1-w}(s)$ association scheme, each treatment of the
pseudo-$L^{*}_{[s-1-w]}(s)$ association scheme can be regarded as each point of
the inducing pseudo-net-$N^{*}$. Hence $s^2$ points of pseudo-net-$N^{*}$
could be arranged at least a net of order $s$, degree $s-1-w$,
which is exactly one complementary net of the above net $N$.
It follows that Result (1) holds.

(2) Now we put $h=i+w$ for $i=1,2,\cdots,s-1-w$. Using Theorem 3.4, we
obtain:

$\mathbf{(i)}$ $s^2$ treatments of the inducing
pseudo-$L^{*}_{[s-1-w]}(s)$ association scheme can be exactly divided into $s$
pairwise parallel common transversals $B^{*h}_1, B^{*h}_2, \cdots,
B^{*h}_s$ of $\mathbf{L_1, \cdots, L_{w}}$ inside the $h$-th paralleling
classification, $h=w+1,w+2,\cdots,s-1$.

$\mathbf{(ii)}$ With respect to $B^{*h_1}_1, B^{*h_1}_2, \cdots,
B^{*h_1}_s$ and $B^{*h_2}_1, B^{*h_2}_2, \cdots, B^{*h_2}_s$, $h_1,
h_2\in \{w+1,w+2,\cdots,s-1\}, h_1\neq h_2$; for any $m,
n\in\{1,2,\cdots,s\}$, $B^{*h_1}_m$ and $B^{*h_2}_n$ have exactly
one common treatment. Where $B^{*h_1}_m$ and $B^{*h_2}_n$ are both
common transversals of $\mathbf{L_1, \cdots, L_{w}}; m,
n=1,2,\cdots,s$.

For $h=w+1,w+2,\cdots,s-1$, it is assumed that $s$ treatments of
$B^{*h}_m$ correspond to $s$ $(m-1)^{,}$s, $m=1,2,\cdots,s$.
Since $B^{*h}_1, B^{*h}_2, \cdots, B^{*h}_s$ are pairwise parallel,
observing $s^2$ treatments of $\mathbf{T}$, thus we exhibit $s$
treatments of $B^{*h}_1$, $s$ treatments of $B^{*h}_2, \cdots$, $s$
treatments of $B^{*h}_s$ just correspond to an $s\times s$
matrix $\mathbf{L_{h}}$ with entries from $\{0,1,\cdots,s-1\}$.

Since $B^{*h}_m$ is common transversal of $\mathbf{L_1, \cdots, L_{w}}$,
$m=1,2,\cdots,s$, there is exactly one treatment of $B^{*h}_m$ in each row and each
column of $\mathbf{T}$; correspondingly, there is exactly one symbol
$m-1$ in every row and every column of $\mathbf{L_{h}}$. Plug the value of $m$ in order,
it can be discovered that each symbol of $\{0,1,\cdots,s-1\}$ appears once
in every row and every column of $\mathbf{L_{h}}$. It follows from Definition 2.2 that
$\mathbf{L_{h}}$ is a Latin square of order $s$, moreover, it and each of
$\mathbf{L_1, \cdots, L_{w}}$ are orthogonal.

With respect to any $h_1, h_2\in\{w+1,w+2,\cdots,s-1\}, h_1\neq h_2$, it is seen
from $\mathbf{(ii)}$ that each of $B^{*h_1}_1, B^{*h_1}_2, \cdots,
B^{*h_1}_s$ and each of $B^{*h_2}_1, B^{*h_2}_2, \cdots, B^{*h_2}_s$
have exactly one common treatment, which corresponds to the ordered pair composed of
one symbol of $\mathbf{L_{h_1}}$ and the other of $\mathbf{L_{h_2}}$.
According to Consequence I of Theorem 3.6, any treatment of the
pseudo-$L^{*}_{[s-1-w]}(s)$ association scheme
can uniquely be expressed as $B^{*h_1}_{m_1}\cap B^{*h_2}_{m_2}$,
here $m_1, m_2\in\{1,2,\cdots,s\}$. Therefore, $s^2$ treatments of $\mathbf{T}$
one-to-one correspond all possible $s^2$ pairs, which could be constituted
with $s$ distinct symbols of $\mathbf{L_{h_1}}$ and $s$ distinct
symbols of $\mathbf{L_{h_2}}$ when $\mathbf{L_{h_1}}$ is superimposed
on $\mathbf{L_{h_2}}$. It follows from Definition 2.3 that
$\mathbf{L_{h_1}}$ and $\mathbf{L_{h_2}}$ are orthogonal.

Hence $\mathbf{L_{w+1}, L_{w+2}, \cdots, L_{s-1}}$ can be added to
the $POL(s,w)$ to obtain a $POL(s,s-1)$. The proof is complete. $ \blacksquare $

{\bf Theorem 3.12 } Let $POL(s,w)=\{\mathbf{L_1, \cdots, L_{w}}\}$,
if either $\mathbf{(iii)}$ with $s>(s-1-w-1)^2$ or $\mathbf{(v)}$
is satisfied, then the $POL(s,w)$ can be extended to a $POL(s,s-1)$.
Where $w\geq 1, s\geq w+4$,\\
$\mathbf{(iii)}$ For $(s-1-w)(s-1)$ treatments, which are first associates
of any treatment $t$ in the pseudo-$L^{*}_{[s-1-w]}(s)$
association scheme, they can be equally divided into $(s-1-w)$
pairwise disjoint sets $A^{*}_1, A^{*}_2, \cdots, A^{*}_{s-1-w}$.\\
$\mathbf{(v)}$ With respect to arbitrary two distinct treatments $t^{*}_1$ and $t^{*}_2$,
they are first associates of the pseudo-$L^{*}_{[s-1-w]}(s)$ association scheme, it is assumed that
$(s-1-w)(s-1)$ treatments, which are first associates of $t^{*}_v$,
can be equally divided into $(s-1-w)$ pairwise disjoint sets $A^{*v}_1,
A^{*v}_2, \cdots, A^{*v}_{s-1-w}$, $v=1,2$. If
$$\left\{
\begin{array}{lll}
(a).\ \ \{t_1^{*}\}\cup A^{*1}_{i_1}=\{t_2^{*}\}\cup
A^{*2}_{j_1}, \ \ for \ t_1^{*}\in A^{*2}_{j_1}, t_2^{*}\in A^{*1}_{i_1};\ moreover,\
\ t_1^{*}\ and \ t_2^{*}\ define\ the\ set\ of\\ \ \ \ \ \ \ \ cardinality\ s,\ in\ which\ arbitrary\ two\
distinct\ treatments\ are\ first\ associates\ of\ the\\ \ \ \ \ \ \ \ pseudo$-$L^{*}_{[s+1-w]}(s)\ association\ scheme.\\
(b).\ \ For\ any\ i_2\in \{1,2,\cdots,s+1-w\}\setminus\{i_1\},\\
\ \ \ \ \ \ (\{t_1^{*}\}\cup A^{*1}_{i_2})\cap(\{t_2^{*}\}\cup A^{*2}_{j_2})=\emptyset, \ \ a\ single\ j_2\in \{1,2,\cdots,s+1-w\}\setminus\{j_1\}.\\
\ \ \ \ \ \ (\{t_1^{*}\}\cup A^{*1}_{i_2})\cap(\{t_2^{*}\}\cup
A^{*2}_{j'})=\{t^{*}_{1j'}\},\ \ each\ j'\in \{1,2,\cdots,s+1-g\}\setminus\{j_1, j_2\}.\\
(c). \ \ For\ any\ j_3\in \{1,2,\cdots,s+1-w\}\setminus\{j_1\},\\
\ \ \ \ \ \ (\{t_2^{*}\}\cup A^{*2}_{j_3})\cap(\{t_1^{*}\}\cup A^{*1}_{i_3})=\emptyset,\ \ a\ single\ i_3\in \{1,2,\cdots,s+1-w\}\setminus\{i_1\}.\\
\ \ \ \ \ \ (\{t_2^{*}\}\cup A^{*2}_{j_3})\cap(\{t_1^{*}\}\cup A^{*1}_{i'})=\{t_{2i'}^{*}\},
\ each\ i'\in \{1,2,\cdots,s+1-w\}\setminus\{i_1, i_3\}.\\
\ \ Here,\ t_{1j'}^{*}\ and\ t_{2i'}^{*}\ are\ commonly\ first\ associates \ of\ both \ t_1^{*}\ and\ t_2^{*}.
\end{array}
\right. $$
are satisfied.

{\bf Proof: } When the condition $s>(s-1-w-1)^2$ is satisfied, so long as
$\mathbf{(iii)}$ is also satisfied in the inducing pseudo-$L^{*}_{[s-1-w]}(s)$ association
scheme. Superadding Theorem 3.9, by Theorem 3.6, we firstly learn the
set of cardinality $s$ containing $t$ and $a^{*}_i$ is
uniquely $\{t\}\cup A^{*}_{i}$, here $a^{*}_i\in A^{*}_i$;
we are ultimately able to infer the pseudo-$L^{*}_{[s-1-w]}(s)$ association
scheme must be an $L_{s-1-w}(s)$ association scheme.

If $\mathbf{(v)}$ is satisfied, by Theorem 3.6, then $s^2$ treatments of
the inducing pseudo-$L^{*}_{[s-1-w]}(s)$ association scheme can form an
$L_{s-1-w}(s)$ association scheme. Though it is possible that the set of cardinality $s$
containing $t$ and $a^{*}_i$ is non-unique in the pseudo-$L^{*}_{[s-1-w]}(s)$
association scheme, however, this does not affect that desired $s-1-w$
paralleling classifications are selected from much more than these $s-1-w$ ones of the
pseudo-$L^{*}_{[s-1-w]}(s)$ association scheme, so that
$\mathbf{(v)}$ holds.
By Theorem 3.11, hence this theorem holds. $ \blacksquare $

\section{ Several specific examples }

Even though the pseudo-$L_{g}(s)$ association scheme is not an $L_{g}(s)$
association scheme, but its $s^2$ treatments probably form an $L_{g}(s)$ association scheme.
Now we suppose that $s^2$ treatments of the pseudo-$L_{g}(s)$ association scheme
may form an $L_{g}(s)$ association scheme. If every treatment of the pseudo-$L_{g}(s)$
association scheme is regarded as each point, then we will utilize its $s^2$ points to
arrange a net $N$ of order $s$, degree $g$. It is defined that two distinct treatments of
the pseudo-$L_{g}(s)$ association scheme are first associates if and only if
they are joined in $N$; and second associates otherwise, but the pseudo-$L_{g}(s)$
association scheme may not be the net $N$, either. Hence we will not explore it in the
following Examples, we only research when $s^2$ treatments of the pseudo-$L_{g}(s)$
association scheme can form an $L_{g}(s)$ association scheme.

{\bf Example 4.1 } Let \begin{center}$\mathbf{T}=\left(
\begin{array}{ccc}
1&2&3\\
4&5&6\\
7&8&9\\
\end{array}
\right),\ \ \mathbf{L_{1}}=\left(
\begin{array}{ccc}
a&b&c\\
b&c&a\\
c&a&b\\
\end{array}
\right)$.\end{center}

According to Definition 2.11, we firstly utilize nine treatments of $\mathbf{T}$ and the above
$\mathbf{L_{1}}$ to construct one $\mathbf{L_{3}(3)}$ association
scheme, it is very obvious that the $\mathbf{L_{3}(3)}$ association
scheme constructed is a pseudo-$L_{3}(3)$ association scheme. Next
we research how $3^2$ treatments of this
pseudo-$L_{3}(3)$ association scheme are arranged onto a $3\times 3$ matrix.

Take the treatment $1$ as an example, there are six treatments $2,3,4,6,7,8$,
which are first associates of it in this pseudo-$L_{3}(3)$
association scheme, such that they can be equally divided into three
pairwise disjoint sets, moreover, two treatments are first
associates within each of these three sets. There are six distinct types of
three pairwise disjoint sets, which are exhibited in
\begin{center}Type $1 \left\{\begin{array}{l}
\hbox{$\nearrow\{2$--3$\}$}\\
\hbox{1$\rightarrow\{$4--7$\}$}\\
\hbox{$\searrow\{$6--8$\}$}\\
\end{array}\right.$, Type $2 \left\{\begin{array}{l}
\hbox{$\nearrow\{2$--3$\}$}\\
\hbox{1$\rightarrow\{$4--6$\}$}\\
\hbox{$\searrow\{$7--8$\}$}\\
\end{array}\right.$, Type $3 \left\{\begin{array}{l}
\hbox{$\nearrow\{2$--4$\}$}\\
\hbox{1$\rightarrow\{$3--6$\}$}\\
\hbox{$\searrow\{$7--8$\}$}\\
\end{array}\right.$,

Type $4 \left\{\begin{array}{l}
\hbox{$\nearrow\{2$--4$\}$}\\
\hbox{1$\rightarrow\{$3--7$\}$}\\
\hbox{$\searrow\{$6--8$\}$}\\
\end{array}\right.$, Type $5 \left\{\begin{array}{l}
\hbox{$\nearrow\{2$--8$\}$}\\
\hbox{1$\rightarrow\{$4--6$\}$}\\
\hbox{$\searrow\{$7--3$\}$}\\
\end{array}\right.$, Type $6 \left\{\begin{array}{l}
\hbox{$\nearrow\{2$--8$\}$}\\
\hbox{1$\rightarrow\{$4--7$\}$}\\
\hbox{$\searrow\{$6--3$\}$}\\
\end{array}\right.$.
\end{center}

Select the remaining treatments except $1$ from this pseudo-$L_{3}(3)$
association scheme, those corresponding results are also acquired in a similar manner.

Now we found: there exist two distinct treatments, which are first
associates, such that the set of cardinality three containing them is non-unique,
arbitrary two of the set are first associates of this pseudo-$L_{3}(3)$ association scheme.
For example, taking two treatments $1$ and $2$, we have $\{$1--2--3$\}$, $\{$1--2--4$\}$,
$\{$1--2--8$\}$.

Displaying six distinct types of three pairwise disjoint sets concerning
every treatment, we shall examine whether this
pseudo-$L_{3}(3)$ association scheme is inversely an $L_{3}(3)$ association
scheme.

(a) When the following restrictive conditions are satisfied, its nine treatments
can form twelve distinct $L_{3}(3)$ association schemes.

If two different treatments are first associates, then we arrange
them to either lie on the same row or column of $\mathbf{T_u,
u}=1,2,\cdots,12$, or correspond to the same symbol of
$\mathbf{L_1}$; and elsewhere otherwise. Therefore, nine treatments
of this pseudo-$L_{3}(3)$ association scheme can exactly be equally divided
into three pairwise disjoint sets $B^i_1, B^i_2, B^i_3$ inside the
$i$-th parallel class, where $i=3u-2,3u-1,3u; u=1,2,\cdots,12$,
which are shown as follows: \begin{center}
$\begin{array}{ccccccccccccc}
\mathbf{T_1}=\mathbf{T},&&&&&&&&\mathbf{T_2}=\left(
\begin{array}{ccc}
1&2&3\\
4&9&6\\
7&8&5\\
\end{array}
\right),&&\mathbf{T_3}=\left(
\begin{array}{ccc}
1&2&3\\
4&5&7\\
6&8&9\\
\end{array}
\right),\\
\end{array}$

Class $1 \left\{\begin{array}{l}
\hbox{$B^{1}_1:\{$1--2--3.$\}$}\\
\hbox{$B^{1}_2:\{$4--5--6.$\}$}\\
\hbox{$B^{1}_3:\{$7--8--9.$\}$}\\
\end{array}\right.$,\;\; Class $4 \left\{\begin{array}{l}
\hbox{$B^{4}_1:\{$1--2--3.$\}$}\\
\hbox{$B^{4}_2:\{$4--9--6.$\}$}\\
\hbox{$B^{4}_3:\{$7--8--5.$\}$}\\
\end{array}\right.$,\;\; Class $7 \left\{\begin{array}{l}
\hbox{$B^{7}_1:\{$1--2--3.$\}$}\\
\hbox{$B^{7}_2:\{$4--5--7.$\}$}\\
\hbox{$B^{7}_3:\{$6--8--9.$\}$}\\
\end{array}\right.$,

Class $2 \left\{\begin{array}{l}
\hbox{$B^{2}_1:\{$1--4--7.$\}$}\\
\hbox{$B^{2}_2:\{$2--5--8.$\}$}\\
\hbox{$B^{2}_3:\{$3--6--9.$\}$}\\
\end{array}\right.$,\;\; Class $5 \left\{\begin{array}{l}
\hbox{$B^{5}_1:\{$1--4--7.$\}$}\\
\hbox{$B^{5}_2:\{$2--9--8.$\}$}\\
\hbox{$B^{5}_3:\{$3--6--5.$\}$}\\
\end{array}\right.$,\;\; Class $8 \left\{\begin{array}{l}
\hbox{$B^{8}_1:\{$1--4--6.$\}$}\\
\hbox{$B^{8}_2:\{$2--5--8.$\}$}\\
\hbox{$B^{8}_3:\{$3--7--9.$\}$}\\
\end{array}\right.$,

Class $3 \left\{\begin{array}{l}
\hbox{$B^{3}_1:\{$1--6--8.$\}$}\\
\hbox{$B^{3}_2:\{$2--4--9.$\}$}\\
\hbox{$B^{3}_3:\{$3--5--7.$\}$}\\
\end{array}\right.$.\;\; Class $6 \left\{\begin{array}{l}
\hbox{$B^{6}_1:\{$1--6--8.$\}$}\\
\hbox{$B^{6}_2:\{$2--4--5.$\}$}\\
\hbox{$B^{6}_3:\{$3--9--7.$\}$}\\
\end{array}\right.$.\;\; Class $9 \left\{\begin{array}{l}
\hbox{$B^{9}_1:\{$1--7--8.$\}$}\\
\hbox{$B^{9}_2:\{$2--4--9.$\}$}\\
\hbox{$B^{9}_3:\{$3--5--6.$\}$}\\
\end{array}\right.$.

$\begin{array}{ccccccccccccc} \mathbf{T_4}=\left(
\begin{array}{ccc}
1&2&3\\
4&9&7\\
6&8&5\\
\end{array}
\right),&&&&&\mathbf{T_5}=\left(
\begin{array}{ccc}
1&2&4\\
3&5&7\\
6&8&9\\
\end{array}
\right),&&&&&\mathbf{T_6}=\left(
\begin{array}{ccc}
1&2&4\\
3&9&7\\
6&8&5\\
\end{array}
\right),\\
\end{array}$

Class $10 \left\{\begin{array}{l}
\hbox{$B^{10}_1:\{$1--2--3.$\}$}\\
\hbox{$B^{10}_2:\{$4--9--7.$\}$}\\
\hbox{$B^{10}_3:\{$6--8--5.$\}$}\\
\end{array}\right.$,\;\; Class $13 \left\{\begin{array}{l}
\hbox{$B^{13}_1:\{$1--2--4.$\}$}\\
\hbox{$B^{13}_2:\{$3--5--7.$\}$}\\
\hbox{$B^{13}_3:\{$6--8--9.$\}$}\\
\end{array}\right.$,\;\; Class $16 \left\{\begin{array}{l}
\hbox{$B^{16}_1:\{$1--2--4.$\}$}\\
\hbox{$B^{16}_2:\{$3--9--7.$\}$}\\
\hbox{$B^{16}_3:\{$6--8--5.$\}$}\\
\end{array}\right.$,

Class $11 \left\{\begin{array}{l}
\hbox{$B^{11}_1:\{$1--4--6.$\}$}\\
\hbox{$B^{11}_2:\{$2--9--8.$\}$}\\
\hbox{$B^{11}_3:\{$3--7--5.$\}$}\\
\end{array}\right.$,\;\; Class $14 \left\{\begin{array}{l}
\hbox{$B^{14}_1:\{$1--3--6.$\}$}\\
\hbox{$B^{14}_2:\{$2--5--8.$\}$}\\
\hbox{$B^{14}_3:\{$4--7--9.$\}$}\\
\end{array}\right.$,\;\; Class $17 \left\{\begin{array}{l}
\hbox{$B^{17}_1:\{$1--3--6.$\}$}\\
\hbox{$B^{17}_2:\{$2--9--8.$\}$}\\
\hbox{$B^{17}_3:\{$4--7--5.$\}$}\\
\end{array}\right.$,

Class $12 \left\{\begin{array}{l}
\hbox{$B^{12}_1:\{$1--7--8.$\}$}\\
\hbox{$B^{12}_2:\{$2--4--5.$\}$}\\
\hbox{$B^{12}_3:\{$3--9--6.$\}$}\\
\end{array}\right.$.\;\; Class $15 \left\{\begin{array}{l}
\hbox{$B^{15}_1:\{$1--7--8.$\}$}\\
\hbox{$B^{15}_2:\{$2--3--9.$\}$}\\
\hbox{$B^{15}_3:\{$4--5--6.$\}$}\\
\end{array}\right.$.\;\; Class $18 \left\{\begin{array}{l}
\hbox{$B^{18}_1:\{$1--7--8.$\}$}\\
\hbox{$B^{18}_2:\{$2--3--5.$\}$}\\
\hbox{$B^{18}_3:\{$4--9--6.$\}$}\\
\end{array}\right.$.

$\begin{array}{ccccccccccccc} \mathbf{T_7}=\left(
\begin{array}{ccc}
1&2&4\\
3&5&6\\
7&8&9\\
\end{array}
\right),&&&&&\mathbf{T_8}=\left(
\begin{array}{ccc}
1&2&4\\
3&9&6\\
7&8&5\\
\end{array}
\right),&&&&&\mathbf{T_9}=\left(
\begin{array}{ccc}
1&2&8\\
4&5&7\\
6&3&9\\
\end{array}
\right),\\
\end{array}$

Class $19 \left\{\begin{array}{l}
\hbox{$B^{19}_1:\{$1--2--4.$\}$}\\
\hbox{$B^{19}_2:\{$3--5--6.$\}$}\\
\hbox{$B^{19}_3:\{$7--8--9.$\}$}\\
\end{array}\right.$,\;\; Class $22 \left\{\begin{array}{l}
\hbox{$B^{22}_1:\{$1--2--4.$\}$}\\
\hbox{$B^{22}_2:\{$3--9--6.$\}$}\\
\hbox{$B^{22}_3:\{$7--8--5.$\}$}\\
\end{array}\right.$,\;\; Class $25 \left\{\begin{array}{l}
\hbox{$B^{25}_1:\{$1--2--8.$\}$}\\
\hbox{$B^{25}_2:\{$4--5--7.$\}$}\\
\hbox{$B^{25}_3:\{$6--3--9.$\}$}\\
\end{array}\right.$,

Class $20 \left\{\begin{array}{l}
\hbox{$B^{20}_1:\{$1--3--7.$\}$}\\
\hbox{$B^{20}_2:\{$2--5--8.$\}$}\\
\hbox{$B^{20}_3:\{$4--6--9.$\}$}\\
\end{array}\right.$,\;\; Class $23 \left\{\begin{array}{l}
\hbox{$B^{23}_1:\{$1--3--7.$\}$}\\
\hbox{$B^{23}_2:\{$2--9--8.$\}$}\\
\hbox{$B^{23}_3:\{$4--6--5.$\}$}\\
\end{array}\right.$,\;\; Class $26 \left\{\begin{array}{l}
\hbox{$B^{26}_1:\{$1--4--6.$\}$}\\
\hbox{$B^{26}_2:\{$2--5--3.$\}$}\\
\hbox{$B^{26}_3:\{$8--7--9.$\}$}\\
\end{array}\right.$,

Class $21 \left\{\begin{array}{l}
\hbox{$B^{21}_1:\{$1--6--8.$\}$}\\
\hbox{$B^{21}_2:\{$2--3--9.$\}$}\\
\hbox{$B^{21}_3:\{$4--5--7.$\}$}\\
\end{array}\right.$.\;\; Class $24 \left\{\begin{array}{l}
\hbox{$B^{24}_1:\{$1--6--8.$\}$}\\
\hbox{$B^{24}_2:\{$2--3--5.$\}$}\\
\hbox{$B^{24}_3:\{$4--9--7.$\}$}\\
\end{array}\right.$.\;\; Class $27 \left\{\begin{array}{l}
\hbox{$B^{27}_1:\{$1--7--3.$\}$}\\
\hbox{$B^{27}_2:\{$2--4--9.$\}$}\\
\hbox{$B^{27}_3:\{$8--5--6.$\}$}\\
\end{array}\right.$.

$\begin{array}{ccccccccccccc} \mathbf{T_{10}}=\left(
\begin{array}{ccc}
1&2&8\\
4&9&7\\
6&3&5\\
\end{array}
\right),&&&&&\mathbf{T_{11}}=\left(
\begin{array}{ccc}
1&2&8\\
4&5&6\\
7&3&9\\
\end{array}
\right),&&&&&\mathbf{T_{12}}=\left(
\begin{array}{ccc}
1&2&8\\
4&9&6\\
7&3&5\\
\end{array}
\right),\\
\end{array}$

Class $28 \left\{\begin{array}{l}
\hbox{$B^{28}_1:\{$1--2--8.$\}$}\\
\hbox{$B^{28}_2:\{$4--9--7.$\}$}\\
\hbox{$B^{28}_3:\{$6--3--5.$\}$}\\
\end{array}\right.$, Class $31 \left\{\begin{array}{l}
\hbox{$B^{31}_1:\{$1--2--8.$\}$}\\
\hbox{$B^{31}_2:\{$4--5--6.$\}$}\\
\hbox{$B^{31}_3:\{$7--3--9.$\}$}\\
\end{array}\right.$, Class $34 \left\{\begin{array}{l}
\hbox{$B^{34}_1:\{$1--2--8.$\}$}\\
\hbox{$B^{34}_2:\{$4--9--6.$\}$}\\
\hbox{$B^{34}_3:\{$7--3--5.$\}$}\\
\end{array}\right.$,

Class $29 \left\{\begin{array}{l}
\hbox{$B^{29}_1:\{$1--4--6.$\}$}\\
\hbox{$B^{29}_2:\{$2--9--3.$\}$}\\
\hbox{$B^{29}_3:\{$8--7--5.$\}$}\\
\end{array}\right.$, Class $32 \left\{\begin{array}{l}
\hbox{$B^{32}_1:\{$1--4--7.$\}$}\\
\hbox{$B^{32}_2:\{$2--5--3.$\}$}\\
\hbox{$B^{32}_3:\{$8--6--9.$\}$}\\
\end{array}\right.$, Class $35 \left\{\begin{array}{l}
\hbox{$B^{35}_1:\{$1--4--7.$\}$}\\
\hbox{$B^{35}_2:\{$2--9--3.$\}$}\\
\hbox{$B^{35}_3:\{$8--6--5.$\}$}\\
\end{array}\right.$,

Class $30 \left\{\begin{array}{l}
\hbox{$B^{30}_1:\{$1--7--3.$\}$}\\
\hbox{$B^{30}_2:\{$2--4--5.$\}$}\\
\hbox{$B^{30}_3:\{$8--9--6.$\}$}\\
\end{array}\right.$. Class $33 \left\{\begin{array}{l}
\hbox{$B^{33}_1:\{$1--6--3.$\}$}\\
\hbox{$B^{33}_2:\{$2--4--9.$\}$}\\
\hbox{$B^{33}_3:\{$8--5--7.$\}$}\\
\end{array}\right.$. Class $36 \left\{\begin{array}{l}
\hbox{$B^{36}_1:\{$1--6--3.$\}$}\\
\hbox{$B^{36}_2:\{$2--4--5.$\}$}\\
\hbox{$B^{36}_3:\{$8--9--7.$\}$}\\
\end{array}\right.$.\end{center}

Test $B^i_1, B^i_2, B^i_3$ and $B^j_1, B^j_2, B^j_3$, here $i,
j\in \{3u-2,3u-1,3u\}, i\neq j$; for any $m, n\in\{1,2,3\}$, it is
obviously seen that $B^i_m$ and $B^j_n$ have exactly one common
treatment. By Consequence I of Theorem 3.6, it is easily verified that $3^2$ treatments of this
pseudo-$L_{3}(3)$ association scheme form the $u$-th $L_{3}(3)$
association scheme, $u=1,2,\cdots,12$.

(b) When anyone of five conditions of Theorem 3.5 holds,
this pseudo-$L_{3}(3)$ association scheme is not an $L_{3}(3)$ association scheme.

Now we consider the first $L_{3}(3)$ association scheme
and the second $L_{3}(3)$ association scheme,
Classes 1 to 6 of this pseudo-$L_{3}(3)$ association scheme are selected.
It is immediately discovered that $B^1_2$ and $B^4_2$ have two common treatments.
Thus we illustrate Conditions $\mathbf{(I)}$ and $\mathbf{(II)}$ are satisfied.

In order to detect that $\{t_1\}\cup A^{1}_{i_1}$ and $\{t_2\}\cup A^{2}_{j_1}$ have two common
treatments, we take $t_1=4$ and $t_2=5$ as an example, they are first associates of this
pseudo-$L_{3}(3)$ association scheme. Six treatments that are first associates of $t_v$, $v=1,2$,
can be equally divided into three pairwise disjoint sets, shown in \begin{center}
$\left\{\begin{array}{l}
\hbox{$\{4\}\cup\{$1--7$\}$}\\
\hbox{$\{4\}\cup\{$5--6$\}$}\\
\hbox{$\{4\}\cup\{$2--9$\}$}\\
\end{array}\right.$; $\left.\begin{array}{l}
\hbox{$\{4\}\cup\{$1--7$\}$}\\
\hbox{$\{4\}\cup\{$9--6$\}$}\\
\hbox{$\{4\}\cup\{$2--5$\}$}\\
\end{array}\right.$ and $\left\{\begin{array}{l}
\hbox{$\{5\}\cup\{$7--8$\}$}\\
\hbox{$\{5\}\cup\{$3--6$\}$}\\
\hbox{$\{5\}\cup\{$2--4$\}$}\\
\end{array}\right.$.
\end{center}
None of $\{4\}\cup\{$1--7$\}, \{4\}\cup\{$5--6$\}, \{4\}\cup\{$2--9$\}$ is equal
to $\{$4--9--6.$\}$, it is obviously seen that Conditions $\mathbf{(III)}$
and $\mathbf{(IV)}$ hold in Theorem 3.5. We notice $(\{4\}\cup\{$2--9$\})\cap (\{5\}\cup\{$7--8$\})=\emptyset$
and $(\{4\}\cup\{$2--9$\})\cap (\{5\}\cup\{$3--6$\})=\emptyset$ hold; and observe $\{4\}\cup\{$2--9$\}$
and $\{5\}\cup\{$2--4$\}$ have two common treatments, that is to say,
the set of cardinality $3$ containing two treatments $4$ and $2$ is non-unique.
It follows that $\mathbf{(V)}$ is satisfied. Clearly,
let $s=3$ in Lemma 2.21, we could verify that three pairwise disjoint sets obtained
satisfy both $\mathbf {(I)}$ and $\mathbf {(I^{'})}$ of Lemma 2.21
when taking $t_1=4$ and $t_2=5$ of
this pseudo-$L_{3}(3)$ association scheme. However, by Theorem 3.5, this
pseudo-$L_{3}(3)$ association scheme is not an $L_{3}(3)$ association scheme.

Although two treatments $3$ and $4$ are second associates of this
pseudo-$L_{3}(3)$ association scheme, we apparently have
$(\{3\}\cup\{$1--2$\})\cap (\{4\}\cup\{$1--2\}$)
=\{1, 2\}$.

Every treatment of the pseudo-$L_{3}(3)$ association scheme is considered as each point,
two distinct points are joined if and only if they are first associates of
the pseudo-$L_{3}(3)$ association scheme. Even if $3^2$ points of the pseudo-$L_{3}(3)$
association scheme are arranged a net $N$ of order $3$, degree $3$,
but it is not the net,
mainly because the line of cardinality $3$ containing two points is non-unique.

{\bf Example 4.2 } Let \begin{center}$\mathbf{T}=\left(
\begin{array}{cccc}
1 &2 &3 &4 \\
5 &6 &7 &8 \\
9 &10&11&12\\
13&14&15&16\\
\end{array}
\right),\ \ \mathbf{L_{1}}=\left(
\begin{array}{cccc}
a&c&d&b\\
d&b&a&c\\
b&d&c&a\\
c&a&b&d\\
\end{array}
\right)$.\end{center}

By Definition 2.11, now sixteen treatments of $\mathbf{T}$ and the above
$\mathbf{L_{1}}$ are used to construct one $\mathbf{L_{3}(4)}$
association scheme, which is obviously a
pseudo-$L_{3}(4)$ association scheme. Next we study how $4^2$ treatments of this
pseudo-$L_{3}(4)$ association scheme are arranged onto a $4\times 4$ matrix.

Choose the treatment $1$, there are nine treatments that are first
associates of it in this pseudo-$L_{3}(4)$ association scheme, such
that they can be equally divided into three pairwise disjoint sets,
further, arbitrary two different treatments are first associates within each
of these three sets. There are two distinct types of three pairwise
disjoint sets, which are shown in
\begin{center}Type $1 \left\{\begin{array}{l}
\hbox{$\nearrow\{2$--3--4$\}$}\\
\hbox{1$\rightarrow\{$13--5--9$\}$}\\
\hbox{$\searrow\{$14--7--12$\}$}\\
\end{array}\right.$ and Type $2 \left\{\begin{array}{l}
\hbox{$\nearrow\{2$--13--14$\}$}\\
\hbox{1$\rightarrow\{$3--5--7$\}$}\\
\hbox{$\searrow\{$4--9--12$\}$}\\
\end{array}\right.$.\end{center}

Using the definition of the $\mathbf{L_{3}(4)}$ association scheme,
we directly write Type $1$; applying the definition of this
pseudo-$L_{3}(4)$ association scheme and the property of the Latin
square $\mathbf{L_{1}}$, we display Type $2$.

Select the remainders except $1$ from this pseudo-$L_{3}(4)$ association
scheme, those corresponding conclusions are also obtained in a similar manner.

We know: there exist two distinct treatments, which are first
associates, such that the set of cardinality four containing them is non-unique,
any two of the set are first associates of this pseudo-$L_{3}(4)$ association scheme.
For instance, taking $1$ and $2$, we have $\{$1--2--3--4$\}$ and $\{$1--2--13--14$\}$.

Manifesting two distinct types of three pairwise disjoint sets about
every treatment, we shall determine if this pseudo-$L_{3}(4)$ association
scheme is conversely an $L_{3}(4)$ association scheme.

(a) When the following restrictive conditions are satisfied, its sixteen treatments
can form two distinct $L_{3}(4)$ association schemes.

Set $\mathbf{T_1}=\mathbf{T}, \mathbf{T_2}=\left(
\begin{array}{cccc}
1 &2 &13&14\\
5 &6 &9 &10\\
7 &8 &11&12\\
3 &4 &15&16\\
\end{array}
\right)$. If two different treatments are first associates, then we put
them to either lie on the same row or column of $\mathbf{T_u,
u}=1,2$, or correspond to the same symbol of
$\mathbf{L_1}$; and elsewhere otherwise. Hence sixteen treatments
of this pseudo-$L_{3}(4)$ association scheme can exactly be equally divided
into four pairwise disjoint sets $B^i_1, B^i_2, B^i_3, B^i_4$ inside the
$i$-th parallel class, here $i=3u-2,3u-1,3u; u=1,2$,
which are exhibited as follows: \begin{center}
Class $1 \left\{\begin{array}{l}
\hbox{$B^{1}_1:\{$1--2--3--4.$\}$}\\
\hbox{$B^{1}_2:\{$5--6--7--8.$\}$}\\
\hbox{$B^{1}_3:\{$9--10--11--12.$\}$}\\
\hbox{$B^{1}_4:\{$13--14--15--16.$\}$}\\
\end{array}\right.$, Class $2 \left\{\begin{array}{l}
\hbox{$B^{2}_1:\{$1--5--9--13.$\}$}\\
\hbox{$B^{2}_2:\{$2--6--10--14.$\}$}\\
\hbox{$B^{2}_3:\{$3--7--11--15.$\}$}\\
\hbox{$B^{2}_4:\{$4--8--12--16.$\}$}\\
\end{array}\right.$,\\ Class $3 \left\{\begin{array}{l}
\hbox{$B^{3}_1:\{$1--7--12--14.$\}$}\\
\hbox{$B^{3}_2:\{$2--8--11--13.$\}$}\\
\hbox{$B^{3}_3:\{$3--5--10--16.$\}$}\\
\hbox{$B^{3}_4:\{$4--6--9--15.$\}$}\\
\end{array}\right.$. Class $4 \left\{\begin{array}{l}
\hbox{$B^{4}_1:\{$1--2--13--14.$\}$}\\
\hbox{$B^{4}_2:\{$5--6--9--10.$\}$}\\
\hbox{$B^{4}_3:\{$7--8--11--12.$\}$}\\
\hbox{$B^{4}_4:\{$3--4--15--16.$\}$}\\
\end{array}\right.$,\\ Class $5 \left\{\begin{array}{l}
\hbox{$B^{5}_1:\{$1--5--7--3.$\}$}\\
\hbox{$B^{5}_2:\{$2--6--8--4.$\}$}\\
\hbox{$B^{5}_3:\{$13--9--11--15.$\}$}\\
\hbox{$B^{5}_4:\{$14--10--12--16.$\}$}\\
\end{array}\right.$, Class $6 \left\{\begin{array}{l}
\hbox{$B^{6}_1:\{$1--9--12--4.$\}$}\\
\hbox{$B^{6}_2:\{$2--10--11--3.$\}$}\\
\hbox{$B^{6}_3:\{$13--5--8--16.$\}$}\\
\hbox{$B^{6}_4:\{$14--6--7--15.$\}$}\\
\end{array}\right.$.\end{center}

Check $B^i_1, B^i_2, B^i_3, B^i_4$ and $B^j_1, B^j_2, B^j_3, B^j_4$, here $i,
j\in \{3u-2,3u-1,3u\}, i\neq j$; for any $m, n\in\{1,2,3,4\}$, it is
obviously shown that $B^i_m$ and $B^j_n$ have exactly one common
treatment. By Consequence I of Theorem 3.6, it is easily illustrated that $4^2$ treatments of this
pseudo-$L_{3}(4)$ association scheme form the $u$-th $L_{3}(4)$
association scheme, $u=1,2$.

(b) When anyone of five conditions of Theorem 3.5 holds,
this pseudo-$L_{3}(4)$ association scheme is not an $L_{3}(4)$ association scheme.

Analyzing two different $L_{3}(4)$ association schemes,
we observe Classes 1 to 6 of this pseudo-$L_{3}(4)$ association scheme.
It is immediately found that $B^1_1$ and $B^4_1$ have two common treatments.
Thus we illustrate Conditions $\mathbf{(I)}$ and $\mathbf{(II)}$ are satisfied.

In order to explore that $\{t_1\}\cup A^{1}_{i_1}$ and $\{t_2\}\cup A^{2}_{j_1}$ have two common
treatments, we take $t_1=1$ and $t_2=2$ as an instance, they are first associates of this
pseudo-$L_{3}(4)$ association scheme. Nine treatments that are first associates of $t_v$, $v=1,2$,
can be equally divided into three pairwise disjoint sets, exhibited in \begin{center}
$\left\{\begin{array}{l}
\hbox{$\{1\}\cup\{$2--3--4$\}$}\\
\hbox{$\{1\}\cup\{$13--5--9$\}$}\\
\hbox{$\{1\}\cup\{$14--7--12$\}$}\\
\end{array}\right.$; $\left.\begin{array}{l}
\hbox{$\{1\}\cup\{$2--13--14$\}$}\\
\hbox{$\{1\}\cup\{$3--5--7$\}$}\\
\hbox{$\{1\}\cup\{$4--9--12$\}$}\\
\end{array}\right.$ and $\left\{\begin{array}{l}
\hbox{$\{2\}\cup\{$1--13--14$\}$}\\
\hbox{$\{2\}\cup\{$4--8--6$\}$}\\
\hbox{$\{2\}\cup\{$3--11--10$\}$}\\
\end{array}\right.$.\end{center}
None of $\{1\}\cup\{$2--3--4$\}, \{1\}\cup\{$13--5--9$\}, \{1\}\cup\{$14--7--12$\}$ is equal
to $\{$1--3--5--7.$\}$, it is obviously learned that Conditions $\mathbf{(III)}$
and $\mathbf{(IV)}$ hold in Theorem 3.5. We note $(\{1\}\cup\{$13--5--9$\})\cap
(\{2\}\cup\{$4--8--6$\})=\emptyset$
and $(\{1\}\cup\{$13--5--9$\})\cap (\{2\}\cup\{$3--11--10\}$)=\emptyset$
hold; and detect $\{1\}\cup\{$14--7--12$\}$
and $\{2\}\cup\{$1--13--14$\}$ have two common treatments, that is,
the set of cardinality $4$ containing two treatments $1$ and $14$ is non-unique.
It follows that $\mathbf{(V)}$ is satisfied. Visibly,
let $s=4$ in Lemma 2.21, we may verify that three pairwise disjoint sets acquired
satisfy both $\mathbf {(I)}$ and $\mathbf {(I^{'})}$ of Lemma 2.21
while taking $t_1=1$ and $t_2=2$ of
this pseudo-$L_{3}(4)$ association scheme. Due to Theorem 3.5, but this
pseudo-$L_{3}(4)$ association scheme is not an $L_{3}(4)$ association scheme.

Though two treatments $3$ and $13$ are second associates of this
pseudo-$L_{3}(4)$ association scheme, we obviously have
$(\{3\}\cup\{$1--2--4$\})\cap (\{13\}\cup\{$1--2--14\}$)
=\{1, 2\}$.

Similarly, even though $4^2$ treatments of the pseudo-$L_{3}(4)$
association scheme are arranged a net $N$ of order $4$, degree $3$,
but it is not the net,
primarily because the line of cardinality $4$ containing two points is non-unique.

{\bf Example 4.3 } Take $w=1, s=w+4=5$ in Theorem 3.12, let \begin{center}
$\mathbf{T}=\left(\begin{array}{ccccc}
1&2&3&4&5\\
6&7&8&9&10\\
11&12&13&14&15\\
16&17&18&19&20\\
21&22&23&24&25\\
\end{array}
\right),\ \ \mathbf{L_{1}}=\left(
\begin{array}{ccccc}
a&d&b&e&c\\
d&b&e&c&a\\
b&e&c&a&d\\
e&c&a&d&b\\
c&a&d&b&e\\
\end{array}
\right)$.
\end{center}

We utilize twenty-five treatments of $\mathbf{T}$ and the above
$\mathbf{L_{1}}$ to construct one $\mathbf{L_{3}(5)}$ association
scheme. By Definition 2.16, one
pseudo-$L^{*}_{[5-1-1]}(5)$ association scheme can be induced by the
$\mathbf{L_{3}(5)}$ association scheme, also by Lemma 2.20, it is a
pseudo-$L_{3}(5)$ association scheme. In other words, its first
associates are regarded as second associates of the $\mathbf{L_{3}(5)}$
association scheme, and vice versa. Next we consider how
$5^2$ treatments of the pseudo-$L^{*}_{[5-1-1]}(5)$
association scheme, which are also those of the $\mathbf{L_{3}(5)}$
association scheme, can be distributed onto a $5\times 5$ matrix.

Taking any treatment $t$, by Definition 2.7, we may calculate $\mathbf{C^{T;
L_{1}}_{t}}$. Now we choose $t=1$ as an example, there are twelve treatments
which are first associates of $1$ in the pseudo-$L^{*}_{[5-1-1]}(5)$ association
scheme, that is, they are exactly second associates of $1$
in the $\mathbf{L_{3}(5)}$ association scheme, such that they and
$1$ are together located in $\mathbf{\Delta^{T;
L_{1}}_{1}}$. Subtract these thirteen treatments of $\mathbf{\Delta^{T;
L_{1}}_{1}}$ from $\mathbf{T}$, those remaining treatments deleted are second associates
of $1$ in the pseudo-$L^{*}_{[5-1-1]}(5)$ association scheme. According to
\begin{center}$\mathbf{\Delta^{T;
L_{1}}_{1}}=\left(\begin{array}{ccccc}
1& & & & \\
 &7 &8 &9 &  \\
 &12&13&  &15\\
 &17&  &19&20\\
 &  &23&24&25\\
\end{array}
\right),\ \ \mathbf{\Delta^{T;
L_{1}}_{1--7}}=\left(\begin{array}{ccccc}
1& & & & \\
 &7 &  &  &  \\
 &  &13&  &15\\
 &  &  &19&  \\
 &  &23&  &25\\
\end{array}
\right)$, $\mathbf{\Delta^{T;
L_{1}}_{1--8}}=\left(\begin{array}{ccccc}
1& & & & \\
 &  &8 &  &  \\
 &  &  &  &15\\
 &17&  &19&20\\
 &  &  &24&  \\
\end{array}
\right)$,\ \ $\mathbf{\Delta^{T;
L_{1}}_{1--9}}=\left(\begin{array}{ccccc}
1& & & & \\
 &  &  &9 &  \\
 &12&  &  &15\\
 &  &  &  &20\\
 &  &23&  &25\\
\end{array}
\right)$.\end{center} Hence we will write the collection of all the transversals
of $\mathbf{L_{1}}$ which contain the treatment $1$, i.e.,
\begin{center}$$\mathbf{C^{T; L_{1}}_{1}}=\mathbf{C^{T;
L_{1}}_{1--7}}\cup \mathbf{C^{T; L_{1}}_{1--8}}\cup \mathbf{C^{T;
L_{1}}_{1--9}}$$ $=\{1-7-13-19-25\}\cup \{1-8-15-17-24\}\cup
\{1-9-12-20-23\}$\\$=\{(1)\;1-7-13-19-25.\;\;(2)\;1-8-15-17-24.\;\;
(3)\;1-9-12-20-23.\}$.

We set $\begin{array}{ccc}
&\nearrow&Y=\{7,13,19,25\}\\
x=1&\rightarrow&Z=\{9,12,20,23\} \\
&\searrow&W=\{8,15,17,24\}
\end{array}$ in Theorem 3.8.\end{center}

Via computing, there are totally fifteen distinct transversals of
$\mathbf{L_{1}}$, such that they fall into three
paralleling classifications of five transversals each.
In this case, let each transversal $B^{*1}_m$ and each transversal
$B^{*2}_n$ be every row and every column of $\mathbf{T^{*}}$, respectively,
$m, n=1,2,3,4,5$, these obtain Classifications 1 and 2; let each
transversal $B^{*3}_m$ be the set of five treatments of $\mathbf{T^{*}}$
corresponding to the same symbol of $\mathbf{L^{*}_{1}}$, this
obtains Classification 3. Where \begin{center}
$\mathbf{T^{*}}=\left(\begin{array}{ccccc}
1&7&13&19&25\\
9&15&16&22&3\\
12&18&24&5&6\\
20&21&2&8&14\\
23&4&10&11&17\\
\end{array}
\right), \ \ \mathbf{L^{*}_{1}}=\left(
\begin{array}{ccccc}
a^{*}&b^{*}&c^{*}&d^{*}&e^{*}\\
e^{*}&a^{*}&b^{*}&c^{*}&d^{*}\\
d^{*}&e^{*}&a^{*}&b^{*}&c^{*}\\
c^{*}&d^{*}&e^{*}&a^{*}&b^{*}\\
b^{*}&c^{*}&d^{*}&e^{*}&a^{*}\\
\end{array}
\right)$.
\end{center}

On balance, $5^2$ treatments of the pseudo-$L^{*}_{[5-1-1]}(5)$ association
scheme can exactly be divided into five pairwise parallel
transversals $B^{*i}_1, B^{*i}_2, B^{*i}_3, B^{*i}_4, B^{*i}_5$ inside
the $i$-th classification, $i=1,2,3$. Furthermore, two transversals of different classifications
have exactly one common treatment. By Consequence I of Theorem 3.6,
it is seen from $\mathbf{T^{*}}$ and $\mathbf{L^{*}_{1}}$ that
its $5^2$ treatments can form one $L_{3}(5)$ association scheme.

It is necessary to explain that five treatments $1, 8, 15, 17, 24$
correspond to five symbols $a, e, d, c, b$ in
$\mathbf{L_{1}}$ of the $\mathbf{L_{3}(5)}$ association scheme
above constructed, respectively. However, they altogether correspond to
the same symbol $a^{*}$ in $\mathbf{L^{*}_1}$ of that $L_{3}(5)$
association scheme formed.

Due to $5>(5-1-1-1)^2$, calculating $\mathbf{C^{T; L_{1}}_{t}}$, we see that
twelve treatments, which are first associates of $t$ in the
pseudo-$L^{*}_{[5-1-1]}(5)$ association scheme, can be equally
divided into three pairwise disjoint sets, and pay attention to Theorem 3.9 again.
With respect to arbitrary two distinct treatments, which are first associates
of the pseudo-$L^{*}_{[5-1-1]}(5)$ association scheme, by Consequence IV of Theorem 3.6
or Theorem 3.8, it is known that the transversal of $\mathbf{L_{1}}$
containing them is sole one, that is, the set of cardinality $5$ containing
them is unique. Moreover, we may verify that $\mathbf{(v)}$ in Theorem 3.12 can be satisfied.
Apply Corollary 3.7, hence it is deduced that the inducing pseudo-$L^{*}_{[5-1-1]}(5)$ association
scheme must be an $L_{3}(5)$ association scheme.

Next we use $5^2$ treatments of the $\mathbf{L_{3}(5)}$
association scheme to arrange a net $N$ of order $5$, degree $3$.
By Definition 2.17, the inducing pseudo-net-$N^{*}$ of order $5$,
degree $5-1-1$ can be obtained. Since this inducing
pseudo-$L^{*}_{[5-1-1]}(5)$ association scheme is an $L_{3}(5)$
association scheme, the pseudo-net-$N^{*}$ is also a net of order $5$,
degree $3$, that is, it is a complementary net of $N$. It
follows that we may imbed them in an affine plane of order $5$.

By Theorem 3.12, $\mathbf{L_2, L_3, L_4}$ can be added to the
$POL(5,1)=\{\mathbf{L_1}\}$ to obtain a $POL(5,4)$. Here

$\mathbf{L_2}=\left(
\begin{array}{ccccc}
a&e&d&c&b\\
b&a&e&d&c\\
c&b&a&e&d\\
d&c&b&a&e\\
e&d&c&b&a\\
\end{array}
\right), \ \ \mathbf{L_3}=\left(
\begin{array}{ccccc}
a&e&d&c&b\\
c&b&a&e&d\\
e&d&c&b&a\\
b&a&e&d&c\\
d&c&b&a&e\\
\end{array}
\right),\ \ \mathbf{L_4}=\left(
\begin{array}{ccccc}
a&c&e&b&d\\
e&b&d&a&c\\
d&a&c&e&b\\
c&e&b&d&a\\
b&d&a&c&e\\
\end{array}
\right)$.

Every treatment is regarded as each point, thus $s^2$ treatments of the
$\mathbf{L_{w+2}(s)}$ association scheme are used to arrange a net $N$ of
order $s$, degree $w+2$, here $w\geq 1$. By Definitions 2.16 and 2.17,
one inducing pseudo-$L^{*}_{[s-1-w]}(s)$ association scheme and
one inducing pseudo-net-$N^{*}$ of order $s$, degree $s-1-w$ are obtained.
It is defined that two distinct treatments of
the pseudo-$L^{*}_{[s-1-w]}(s)$ association scheme are first associates if and only if
they are joined in pseudo-net-$N^{*}$; and second associates otherwise.
Due to Lemma 2.20, the pseudo-$L^{*}_{[s-1-w]}(s)$ association scheme
is a pseudo-$L_{s-1-w}(s)$ association scheme. Even if $s^2$ treatments of
the inducing pseudo-$L^{*}_{[s-1-w]}(s)$ association scheme can form an
$L_{s-1-w}(s)$ association scheme and $s^2$ points of the inducing pseudo-net-$N^{*}$ are
arranged a net of order $s$, degree $s-1-w$, however, the pseudo-$L^{*}_{[s-1-w]}(s)$
association scheme may not be an $L_{s-1-w}(s)$ association scheme
and the inducing pseudo-net-$N^{*}$ may not be the net. These are illustrated
in the following Examples. If $s^2$ treatments of
the inducing pseudo-$L^{*}_{[s-1-w]}(s)$ association scheme cannot form an
$L_{s-1-w}(s)$ association scheme in any case, then $s^2$ points of the inducing
pseudo-net-$N^{*}$ must not be arranged a net of order $s$, degree
$s-1-w$. Next we only study when $s^2$ treatments of the inducing
pseudo-$L^{*}_{[s-1-w]}(s)$ association scheme can form an $L_{s-1-w}(s)$
association scheme.

{\bf Example 4.4 } Let \begin{center} $\mathbf{T}=\left(
\begin{array}{ccccccc}
1 &2 &3 &4 &5 &6 &7 \\
8 &9 &10&11&12&13&14\\
15&16&17&18&19&20&21\\
22&23&24&25&26&27&28\\
29&30&31&32&33&34&35\\
36&37&38&39&40&41&42\\
43&44&45&46&47&48&49\\
\end{array}
\right)$.\ $\mathbf{L_{1}}=\left(
\begin{array}{ccccccc}
a&c&e&g&b&d&f\\
g&b&d&f&a&c&e\\
f&a&c&e&g&b&d\\
e&g&b&d&f&a&c\\
d&f&a&c&e&g&b\\
c&e&g&b&d&f&a\\
b&d&f&a&c&e&g\\
\end{array}
\right)$,

$\mathbf{L_{2}}=\left(
\begin{array}{ccccccc}
a&d&g&c&f&b&e\\
f&b&e&a&d&g&c\\
d&g&c&f&b&e&a\\
b&e&a&d&g&c&f\\
g&c&f&b&e&a&d\\
e&a&d&g&c&f&b\\
c&f&b&e&a&d&g\\
\end{array}
\right), \ \ \ \ \mathbf{L_{3}}=\left(
\begin{array}{ccccccc}
a&g&f&e&d&c&b\\
c&b&a&g&f&e&d\\
e&d&c&b&a&g&f\\
g&f&e&d&c&b&a\\
b&a&g&f&e&d&c\\
d&c&b&a&g&f&e\\
f&e&d&c&b&a&g\\
\end{array}
\right)$.
\end{center}
Here $POL(7,3)=\{\mathbf{L_{1}, L_{2}, L_{3}}\}$.

Take $s=7, w=1,2,3$ in Theorem 3.12, forty-nine treatments of $\mathbf{T}$ and the
above $POL(7,3)$ are used to construct a few $\mathbf{L_{w+2}(7)}$
association schemes. By Definition 2.16, we know a few
pseudo-$L^{*}_{[7-1-w]}(7)$ association schemes can be
induced by a few $\mathbf{L_{w+2}(7)}$ association schemes, moreover,
they are a few pseudo-$L_{7-1-w}(7)$ association schemes. Now we
study how $7^2$ treatments of a few pseudo-$L^{*}_{[7-1-w]}(7)$
association schemes, which are also those of a few $\mathbf{L_{w+2}(7)}$
association schemes, can be distributed onto a $7\times 7$ matrix.

{\bf Case 1: w=1. } Using $\mathbf{L_{1}}$, we may construct one
$\mathbf{L_{3}(7)}$ association scheme, thus one inducing
pseudo-$L^{*}_{[7-1-1]}(7)$ association scheme is obtained.

Choosing two treatments $1$ and $13$ as an example, which are first associates
of the pseudo-$L^{*}_{[7-1-1]}(7)$ association scheme, we may break $\mathbf{C^{T;
L_{1}}_{1}}$ and $\mathbf{C^{T; L_{1}}_{13}}$ down. Firstly, by
Definition 2.7, we show
$$\mathbf{\Delta^{T; L_{1}}_{1}}=\left(
\begin{array}{ccccccc}
1&  &  &  &  &  &  \\
 &9 &10&11&  &13&14\\
 &  &17&18&19&20&21\\
 &23&24&25&26&  &28\\
 &30&  &32&33&34&35\\
 &37&38&39&40&41&  \\
 &44&45&  &47&48&49\\
\end{array}
\right), \mathbf{\Delta^{T; L_{1}}_{13}}=\left(
\begin{array}{ccccccc}
1 &  &3 &4 &5 &  &7 \\
  &  &  &  &  &13&  \\
15&16&  &18&19&  &21\\
22&23&24&25&26&  &  \\
29&30&31&  &33&  &35\\
  &37&38&39&40&  &42\\
43&44&45&46&  &  &49\\
\end{array}
\right).$$

Secondly, by Definitions 2.8 and 2.9, we will exhibit the collection of
all the transversals of $\mathbf{L_{1}}$ which contain the treatment
$1$, i.e., $\mathbf{C^{T; L_{1}}_{1}}=\mathbf{C^{T;
L_{1}}_{1--9}}\cup \mathbf{C^{T; L_{1}}_{1--10}}\cup \mathbf{C^{T;
L_{1}}_{1--11}}\cup \mathbf{C^{T; L_{1}}_{1--13}}\cup \mathbf{C^{T;
L_{1}}_{1--14}}$. Similarly, we write again $\mathbf{C^{T;
L_{1}}_{13}}=\mathbf{C^{T; L_{1}}_{1--13}}\cup \mathbf{C^{T;
L_{1}}_{13--3}}\cup \mathbf{C^{T; L_{1}}_{13--4}}\cup \mathbf{C^{T;
L_{1}}_{13--5}}\cup \mathbf{C^{T; L_{1}}_{13--7}}$.
Here\begin{center} $\begin{array}{ccccccccccccccccccccc}
\mathbf{C^{T; L_{1}}_{1--9}}&&&&&&&&&\mathbf{C^{T; L_{1}}_{1--10}}
&&&&&&&&&\mathbf{C^{T; L_{1}}_{1--11}}\\
\end{array}$

=$\{$(1) 1--9--17--25--33--41--49.\;\;=$\{$(1)
1--10--19--28--30--39--48.\;\;=$\{$(1) 1--11--21--24--34--37--47.

(2) 1--9--18--28--34--40--45.\;\;\;\;(2)
1--10--18--23--35--41--47.\;\;\;\;(2) 1--11--20--28--33--38--44.

(3) 1--9--21--26--32--38--48.$\}$,\;\;(3)
1--10--20--26--32--37--49.$\}$,\;\;\;(3)
1--11--17--23--35--40--48.$\}$,

$\begin{array}{ccccccccccccccccccccc}
\mathbf{C^{T; L_{1}}_{13--3}}&&&&&&&&&\mathbf{C^{T; L_{1}}_{13--4}}
&&&&&&&&&\mathbf{C^{T; L_{1}}_{13--5}}\\
\end{array}$

=$\{$(1) 13--3--16--26--29--39--49.\;\;=$\{$(1)
13--4--15--24--33--42--44.\;=$\{$(1) 13--5--21--22--30--38--46.

(2) 13--3--15--23--35--40--46.\;\;\;\;(2)
13--4--16--22--35--40--45.\;\;\;\;(2) 13--5--15--25--31--37--49.

(3) 13--3--19--25--30--42--43.$\}$,\;\;(3)
13--4--21--26--31--37--43.$\}$,\;\;\;(3)
13--5--18--23--29--42--45.$\}$,

$\begin{array}{ccccccccccccccccccccc}
\mathbf{C^{T; L_{1}}_{1--14}}&&&&&&&&&\mathbf{C^{T; L_{1}}_{1--13}}
&&&&&&&&&\mathbf{C^{T; L_{1}}_{13--7}}\\
\end{array}$

=$\{$(1) 1--14--20--26--32--38--44.\;\;=$\{$(1)
1--13--18--23--35--40--45.\;=$\{$(1) 13--7--19--25--31--37--43.

(2) 1--14--20--23--32--40--45.\;\;\;\;(2)
1--13--18--26--35--38--44.\;\;\;\;(2) 13--7--19--24--29--37--46.

(3) 1--14--17--26--34--39--44.\;\;\;(3)
1--13--18--24--30--40--49.\;\;\;\;(3) 13--7--19--22--31--39--44.

(4) 1--14--19--24--32--41--44.\;\;\;\;(4)
1--13--19--25--35--37--45.\;\;\;\;(4) 13--7--18--23--31--40--43.

(5) 1--14--20--25--30--38--47.$\}$,\;\;(5)
1--13--21--23--33--39--45.$\}$,\;\;(5)
13--7--16--25--33--38--43.$\}$.
\end{center}

Thirdly, we analyze the combination of thirty treatments,
which are first associates of $1$ or $13$ in the pseudo-$L^{*}_{[7-1-1]}(7)$
association scheme. As for $1$, since these thirty treatments can be equally divided
into five pairwise disjoint sets, so that arbitrary two different
treatments inside each of these five sets are first associates of
the pseudo-$L^{*}_{[7-1-1]}(7)$ association scheme.
Furthermore, there are two distinct types of these five sets, shown
in \begin{center} Type $1 \left\{\begin{array}{l}
\hbox{(1) of $\mathbf{C^{T; L_{1}}_{1--9}}$}\\
\hbox{(1) of $\mathbf{C^{T; L_{1}}_{1--10}}$}\\
\hbox{(1) of $\mathbf{C^{T; L_{1}}_{1--11}}$}\\
\hbox{(1) of $\mathbf{C^{T; L_{1}}_{1--13}}$}\\
\hbox{(1) of $\mathbf{C^{T; L_{1}}_{1--14}}$}\\
\end{array}\right.$ and Type $2 \left\{\begin{array}{l}
\hbox{(1) of $\mathbf{C^{T; L_{1}}_{1--9}}$}\\
\hbox{(1) of $\mathbf{C^{T; L_{1}}_{1--10}}$}\\
\hbox{(1) of $\mathbf{C^{T; L_{1}}_{1--11}}$}\\
\hbox{(2) of $\mathbf{C^{T; L_{1}}_{1--13}}$}\\
\hbox{(2) of $\mathbf{C^{T; L_{1}}_{1--14}}$}\\
\end{array}\right.$. \end{center}

Similarly, there are two distinct types of five pairwise disjoint
sets concerning the treatment $13$, exhibited in \begin{center} Type
$I \left\{\begin{array}{l}
\hbox{(1) of $\mathbf{C^{T; L_{1}}_{13--3}}$}\\
\hbox{(1) of $\mathbf{C^{T; L_{1}}_{13--4}}$}\\
\hbox{(1) of $\mathbf{C^{T; L_{1}}_{13--5}}$}\\
\hbox{(1) of $\mathbf{C^{T; L_{1}}_{1--13}}$}\\
\hbox{(1) of $\mathbf{C^{T; L_{1}}_{13--7}}$}\\
\end{array}\right.$ and Type $II \left\{\begin{array}{l}
\hbox{(1) of $\mathbf{C^{T; L_{1}}_{13--3}}$}\\
\hbox{(1) of $\mathbf{C^{T; L_{1}}_{13--4}}$}\\
\hbox{(1) of $\mathbf{C^{T; L_{1}}_{13--5}}$}\\
\hbox{(4) of $\mathbf{C^{T; L_{1}}_{1--13}}$}\\
\hbox{(4) of $\mathbf{C^{T; L_{1}}_{13--7}}$}\\
\end{array}\right.$. \end{center}

Select the remaining treatments except $1$ and $13$ from the pseudo-$L^{*}_{[7-1-1]}(7)$
association scheme, those corresponding conclusions are obviously acquired in a similar manner.

We see: there exist two distinct treatments, which are first
associates of the pseudo-$L^{*}_{[7-1-1]}(7)$ association scheme,
such that the transversal of $\mathbf{L_{1}}$ containing them is non-unique.
For instance, these two treatments are $1$ and $13$.

Finally, observing two distinct types of five pairwise disjoint sets
about every treatment, we shall discuss whether the inducing pseudo-$L^{*}_{[7-1-1]}(7)$
association scheme is an $L_{5}(7)$ association scheme.

(a) When the following confining conditions are satisfied, its forty-nine
treatments can form one $L_{5}(7)$ association scheme.

After calculating $\mathbf{C^{T; L_{1}}_{t}}$, we are able to pick
thirty-five distinct transversals of $\mathbf{L_{1}}$ out from more
than these thirty-five ones of the pseudo-$L^{*}_{[7-1-1]}(7)$ association scheme,
so that they fall into five paralleling classifications of seven transversals each.
In this case, each transversal $B^{*1}_m$ and each transversal $B^{*2}_n$ are every row and every
column of $\mathbf{T^{*}}$, respectively, $m, n=1,\cdots,7$, these obtain
Classifications 1 and 2; each transversal $B^{*(k+2)}_m$ is the set of seven
treatments of $\mathbf{T^{*}}$ altogether corresponding to the same symbol
of $\mathbf{L^{*}_{k}}$, this obtains Classification ($k+2$), $k=1,2,3$. Where

$\mathbf{T^{*}}=\left(\begin{array}{ccccccc}
1 &40&23&13&45&35&18\\
26&9 &48&31&21&4 &36\\
44&34&17&7 &39&22&12\\
20&3 &42&25&8 &47&30\\
38&28&11&43&33&16&6 \\
14&46&29&19&2 &41&24\\
32&15&5 &37&27&10&49\\
\end{array}
\right), \ \ \mathbf{L^{*}_1}=\left(
\begin{array}{ccccccc}
a^{*}&g^{*}&f^{*}&e^{*}&d^{*}&c^{*}&b^{*}\\
b^{*}&a^{*}&g^{*}&f^{*}&e^{*}&d^{*}&c^{*}\\
c^{*}&b^{*}&a^{*}&g^{*}&f^{*}&e^{*}&d^{*}\\
d^{*}&c^{*}&b^{*}&a^{*}&g^{*}&f^{*}&e^{*}\\
e^{*}&d^{*}&c^{*}&b^{*}&a^{*}&g^{*}&f^{*}\\
f^{*}&e^{*}&d^{*}&c^{*}&b^{*}&a^{*}&g^{*}\\
g^{*}&f^{*}&e^{*}&d^{*}&c^{*}&b^{*}&a^{*}\\
\end{array}
\right)$,

$\mathbf{L^{*}_2}=\left(
\begin{array}{ccccccc}
a^{*}&b^{*}&c^{*}&d^{*}&e^{*}&f^{*}&g^{*}\\
f^{*}&g^{*}&a^{*}&b^{*}&c^{*}&d^{*}&e^{*}\\
d^{*}&e^{*}&f^{*}&g^{*}&a^{*}&b^{*}&c^{*}\\
b^{*}&c^{*}&d^{*}&e^{*}&f^{*}&g^{*}&a^{*}\\
g^{*}&a^{*}&b^{*}&c^{*}&d^{*}&e^{*}&f^{*}\\
e^{*}&f^{*}&g^{*}&a^{*}&b^{*}&c^{*}&d^{*}\\
c^{*}&d^{*}&e^{*}&f^{*}&g^{*}&a^{*}&b^{*}\\
\end{array}
\right), \ \ \mathbf{L^{*}_3}=\left(
\begin{array}{ccccccc}
a^{*}&d^{*}&g^{*}&c^{*}&f^{*}&b^{*}&e^{*}\\
c^{*}&f^{*}&b^{*}&e^{*}&a^{*}&d^{*}&g^{*}\\
e^{*}&a^{*}&d^{*}&g^{*}&c^{*}&f^{*}&b^{*}\\
g^{*}&c^{*}&f^{*}&b^{*}&e^{*}&a^{*}&d^{*}\\
b^{*}&e^{*}&a^{*}&d^{*}&g^{*}&c^{*}&f^{*}\\
d^{*}&g^{*}&c^{*}&f^{*}&b^{*}&e^{*}&a^{*}\\
f^{*}&b^{*}&e^{*}&a^{*}&d^{*}&g^{*}&c^{*}\\
\end{array}
\right)$.

It is very obvious that $\mathbf{L^{*}_1}, \mathbf{L^{*}_2}$ and
$\mathbf{L^{*}_3}$ are mutually orthogonal.

As a whole, $7^2$ treatments of the pseudo-$L^{*}_{[7-1-1]}(7)$ association
scheme can be exactly divided into seven pairwise parallel transversals
$B^{*i}_1, \cdots, B^{*i}_7$ inside the $i$-th classification, $i=1,2,3,4,5$.
Moreover, two transversals of different classifications
have exactly one common treatment. By Consequence I of Theorem 3.6,
we may verify that its $7^2$ treatments can form one $L_{5}(7)$ association scheme.

(b) When anyone of five conditions of Theorem 3.5 is established,
the pseudo-$L^{*}_{[7-1-1]}(7)$ association scheme is not an $L_{5}(7)$ association scheme.

Besides five paralleling classifications of that $L_{5}(7)$ association
scheme formed, we may provide extra four distinct classes of the
pseudo-$L^{*}_{[7-1-1]}(7)$ association scheme again:

Class $6 \left\{\begin{array}{l}
\hbox{$B^{*6}_1:\{$1--13--18--26--35--38--44.$\}$}\\
\hbox{$B^{*6}_2:\{$2--14--19--27--29--39--45.$\}$}\\
\hbox{$B^{*6}_3:\{$3--8--20--28--30--40--46.$\}$}\\
\hbox{$B^{*6}_4:\{$4--9--21--22--31--41--47.$\}$}\\
\hbox{$B^{*6}_5:\{$5--10--15--23--32--42--48.$\}$}\\
\hbox{$B^{*6}_6:\{$6--11--16--24--33--36--49.$\}$}\\
\hbox{$B^{*6}_7:\{$7--12--17--25--34--37--43.$\}$}\\
\end{array}\right.$,\ \  Class $7 \left\{\begin{array}{l}
\hbox{$B^{*7}_1:\{$1--14--20--23--32--40--45.$\}$}\\
\hbox{$B^{*7}_2:\{$2--8--21--24--33--41--46.$\}$}\\
\hbox{$B^{*7}_3:\{$3--9--15--25--34--42--47.$\}$}\\
\hbox{$B^{*7}_4:\{$4--10--16--26--35--36--48.$\}$}\\
\hbox{$B^{*7}_5:\{$5--11--17--27--29--37--49.$\}$}\\
\hbox{$B^{*7}_6:\{$6--12--18--28--30--38--43.$\}$}\\
\hbox{$B^{*7}_7:\{$7--13--19--22--31--39--44.$\}$}\\
\end{array}\right.$,

Class $8  \left\{\begin{array}{l}
\hbox{$B^{*8}_1:\{$1--13--18--24--30--40--49.$\}$}\\
\hbox{$B^{*8}_2:\{$2--14--19--25--31--41--43.$\}$}\\
\hbox{$B^{*8}_3:\{$3--8--20--26--32--42--44.$\}$}\\
\hbox{$B^{*8}_4:\{$4--9--21--27--33--36--45.$\}$}\\
\hbox{$B^{*8}_5:\{$5--10--15--28--34--37--46.$\}$}\\
\hbox{$B^{*8}_6:\{$6--11--16--22--35--38--47.$\}$}\\
\hbox{$B^{*8}_7:\{$7--12--17--23--29--39--48.$\}$}\\
\end{array}\right.$,\ \  Class $9 \left\{\begin{array}{l}
\hbox{$B^{*9}_1:\{$1--14--20--25--30--38--47.$\}$}\\
\hbox{$B^{*9}_2:\{$2--8--21--26--31--39--48.$\}$}\\
\hbox{$B^{*9}_3:\{$3--9--15--27--32--40--49.$\}$}\\
\hbox{$B^{*9}_4:\{$4--10--16--28--33--41--43.$\}$}\\
\hbox{$B^{*9}_5:\{$5--11--17--22--34--42--44.$\}$}\\
\hbox{$B^{*9}_6:\{$6--12--18--23--35--36--45.$\}$}\\
\hbox{$B^{*9}_7:\{$7--13--19--24--29--37--46.$\}$}\\
\end{array}\right.$.

It is obviously found that $\mathbf{(I)}$ and $\mathbf{(III)}$ have been
established in Theorem 3.5. We notice $B^{*2}_1\cap B^{*6}_1=\{$1--14--20--26--32--38--44.$\}\cap \{$1--13--18--26--35--38--44.$\}=\{1, 26, 38, 44\}$.
As for $t^{*}_1=1$, none of (1) in $\mathbf{C^{T; L_{1}}_{1--9}}$, (1) in
$\mathbf{C^{T; L_{1}}_{1--10}}$, (1) in
$\mathbf{C^{T; L_{1}}_{1--11}}$, (1) in
$\mathbf{C^{T; L_{1}}_{1--13}}$, (1) in
$\mathbf{C^{T; L_{1}}_{1--14}}$ is equal
to $B^{*6}_1$. Thus these illustrate Conditions $\mathbf{(II)}$ and $\mathbf{(IV)}$
of Theorem 3.5 are established. Set $t^{*}_1=1$ and $t^{*}_2=13$, it is seen that $(\{1\}\cup\{$9--17--25--33--41--49$\})\cap (\{13\}\cup\{$3--15--23--35--40--46$\})=\emptyset$
and $(\{1\}\cup\{$9--17--25--33--41--49$\})\cap (\{13\}\cup\{$4--16--22--35--40--45\}$)=\emptyset$
hold; and it is detected from $\mathbf{C^{T; L_{1}}_{1--13}}$ that
the transversal of $\mathbf{L_{1}}$ containing $1$ and $13$ is not sole one, that is,
the set of cardinality $7$ containing them is non-unique.
It follows that $\mathbf{(V)}$ holds. By Theorem 3.5, but the
inducing pseudo-$L^{*}_{[7-1-1]}(7)$ association scheme is not an $L_{5}(7)$ association scheme.

Even if two treatments $13$ and $14$ are second associates of the
pseudo-$L^{*}_{[7-1-1]}(7)$ association scheme, we obviously have
$(\{13\}\cup\{$1--18--26--35--38--44$\})\cap (\{14\}\cup\{$1--17--26--34--39--44\}$)
=\{1, 26, 44\}$.

{\bf Case 2: w=2. } Using $\mathbf{L_{1}, L_{2}}$, we may construct one
$\mathbf{L_{4}(7)}$ association scheme, thus one inducing
pseudo-$L^{*}_{[7-1-2]}(7)$ association scheme is acquired.

Taking any treatment $t$, we may calculate $\mathbf{C^{T; L_{1},
L_{2}}_{t}}$. We now choose the treatment $1$ as an example. There are
twenty-four treatments which are first associates of $1$ in the
pseudo-$L^{*}_{[7-1-2]}(7)$ association scheme, that is, they are
exactly second associates of $1$
in the $\mathbf{L_{4}(7)}$ association scheme, such that they and
$1$ are together located in
$$\mathbf{\Delta^{T; L_{1}, L_{2}}_{1}}=\left(
\begin{array}{ccccccc}
1&  &  &  &  &  &  \\
 &9 &10&  &  &13&14\\
 &  &17&18&19&20&  \\
 &23&  &25&26&  &28\\
 &30&  &32&33&  &35\\
 &  &38&39&40&41&  \\
 &44&45&  &  &48&49\\
\end{array}
\right).$$

Next we will show the collection of all the common transversals of
$\mathbf{L_{1}, L_{2}}$ containing $1$, i.e.,
$\mathbf{C^{T; L_{1}, L_{2}}_{1}}=\mathbf{C^{T; L_{1},
L_{2}}_{1--9}}\cup \mathbf{C^{T; L_{1}, L_{2}}_{1--10}}\cup
\mathbf{C^{T; L_{1}, L_{2}}_{1--13}}\cup \mathbf{C^{T; L_{1},
L_{2}}_{1--14}}$. Here
\begin{center} $\mathbf{C^{T; L_{1}, L_{2}}_{1--9}}=\{$ (1)
1--9--17--25--33--41--49.$\}$, $\mathbf{C^{T; L_{1},
L_{2}}_{1--10}}=\{$ (1) 1--10--19--28--30--39--48.$\}$,

$\mathbf{C^{T; L_{1}, L_{2}}_{1--13}}=\{$ (1)
1--13--18--23--35--40--45.$\}$, $\mathbf{C^{T; L_{1},
L_{2}}_{1--14}}=\{$ (1) 1--14--20--26--32--38--44.$\}$.

We put
$\begin{array}{ccc}
&\nearrow&A^{1}_1=\{9,17,25,33,41,49\}\\
t_1=1&\rightarrow&A^{1}_2=\{10,19,28,30,39,48\}\\
&\searrow&A^{1}_3=\{13,18,23,35,40,45\}\\
&\searrow&A^{1}_4=\{14,20,26,32,38,44\}
\end{array}$ in Consequence III of Theorem 3.6.\end{center}

It is emphasized that these twenty-four treatments can be
equally divided into four pairwise disjoint sets, so that arbitrary two different
treatments inside each of these four sets are first associates
of the pseudo-$L^{*}_{[7-1-2]}(7)$ association scheme.

Select the remainders except $1$ from the pseudo-$L^{*}_{[7-1-2]}(7)$
association scheme, those corresponding results are also obtained in a similar manner.

Calculating $\mathbf{C^{T; L_{1}, L_{2}}_{t}}$, there are entirely
twenty-eight distinct common transversals of $\mathbf{L_{1}, L_{2}}$,
such that they fall into four paralleling classifications of seven transversals each.
In this case, seven pairwise parallel transversals $B^{*i}_1, \cdots, B^{*i}_7$
are identical to those of (a) within {\bf Case 1}, $i=1,2,3,4$.

In sum, $7^2$ treatments of the pseudo-$L^{*}_{[7-1-2]}(7)$
association scheme can be exactly divided into seven pairwise
parallel transversals $B^{*i}_1, \cdots, B^{*i}_7$ inside the $i$-th
classification, $i=1,2,3,4$. Moreover, two transversals of different classifications
have exactly one common treatment. By Consequence I of Theorem 3.6,
it is verified that its $7^2$ treatments can form one $L_{4}(7)$ association scheme.

With respect to arbitrary two treatments, which are first
associates of the pseudo-$L^{*}_{[7-1-2]}(7)$ association scheme, via computing, it
is shown that the common transversal of $\mathbf{L_{1}, L_{2}}$
containing them is sole one, that is, the set of cardinality $7$ containing
them is unique. Further, it may be illustrated that
$\mathbf{(v)}$ in Consequence III of Theorem 3.6 can be satisfied.
By Corollary 3.7, it goes without saying that the inducing
pseudo-$L^{*}_{[7-1-2]}(7)$ association scheme must be
an $L_{4}(7)$ association scheme.

{\bf Case 3: w=3. } Using $\mathbf{L_{1}, L_{2}, L_{3}}$, we may
construct one $\mathbf{L_{5}(7)}$ association scheme, thus one inducing
pseudo-$L^{*}_{[7-1-3]}(7)$ association scheme is acquired.

Selecting $x=1$ in Theorem 3.8, we now count $\mathbf{C^{T; L_{1},
L_{2}, L_{3}}_{1}}$. We will display the
collection of all the common transversals of $\mathbf{L_{1}, L_{2}, L_{3}}$
containing the treatment $1$, i.e., $\mathbf{C^{T; L_{1}, L_{2},
L_{3}}_{1}}=\mathbf{C^{T; L_{1}, L_{2}, L_{3}}_{1--9}}\cup
\mathbf{C^{T; L_{1}, L_{2}, L_{3}}_{1--13}}\cup \mathbf{C^{T; L_{1},
L_{2}, L_{3}}_{1--14}}$. Here
$\mathbf{C^{T; L_{1}, L_{2}, L_{3}}_{1--9}}=\{$(1)
1--9--17--25--33--41--49.$\}$,\\ $\mathbf{C^{T; L_{1}, L_{2},
L_{3}}_{1--13}}=\{$(1) 1--13--18--23--35--40--45.$\}$,
$\mathbf{C^{T; L_{1}, L_{2}, L_{3}}_{1--14}}=\{$(1)
1--14--20--26--32--38--44.$\}$.
\begin{center}
Let $\begin{array}{ccc}
&\nearrow&Y=\{13,18,23,35,40,45\}\\
x=1&\rightarrow&Z=\{14,20,26,32,38,44\} \\
&\searrow&W=\{9,17,25,33,41,49\}
\end{array}$ in Theorem 3.8.\end{center}

It is underlined that there are eighteen treatments which are first associates of $1$ in the
pseudo-$L^{*}_{[7-1-3]}(7)$ association scheme, such that they can be
equally divided into three pairwise disjoint sets.
Take the remainders except $1$ from the pseudo-$L^{*}_{[7-1-3]}(7)$
association scheme, it may be known that those similar conclusions are also acquired.

Through calculations, there are totally twenty-one distinct common
transversals of $\mathbf{L_{1}, L_{2}, L_{3}}$, such that they
fall into three paralleling classifications of seven transversals each.
In the circumstances, Classifications 1 to 3 are identical to those of {\bf Case 2}.

In short, $7^2$ treatments of the pseudo-$L^{*}_{[7-1-3]}(7)$ association
scheme can be exactly divided into seven pairwise parallel transversals
$B^{*i}_1, \cdots, B^{*i}_7$ inside the $i$-th classification, $i=1,2,3$.
Moreover, two transversals of different classifications
have exactly one common treatment. By Consequence I of Theorem 3.6,
it is asserted that its $7^2$ treatments can form one $L_{3}(7)$ association scheme.

Because of $7>(7-1-3-1)^2$, calculating $\mathbf{C^{T; L_{1}, L_{2}, L_{3}}_{t}}$,
we know that eighteen treatments, which are first associates of $t$ in the
pseudo-$L^{*}_{[7-1-3]}(7)$ association scheme, can be equally divided
into three pairwise disjoint sets, and pay attention to Theorem 3.9 again.
With respect to arbitrary two distinct treatments, which are first associates of the
pseudo-$L^{*}_{[7-1-3]}(7)$ association scheme, according to Consequence IV of Theorem 3.6 or Theorem 3.8,
it is seen that the common transversal of $\mathbf{L_{1}, L_{2}, L_{3}}$
containing them is sole one, that is, the set of cardinality $7$ containing
them is unique. Further, we check $\mathbf{(v)}$ in Theorem 3.12 can be satisfied.
By Corollary 3.7, thus it is inferred that the inducing pseudo-$L^{*}_{[7-1-3]}(7)$ association
scheme must be an $L_{3}(7)$ association scheme.

By Theorem 3.12, $\mathbf{L_4, L_5, L_6}$ can be added to the
$POL(7,3)=\{\mathbf{L_1, L_2, L_3}\}$ to obtain a $POL(7,6)$. Here
\begin{center}
$\mathbf{L_4}=\left(
\begin{array}{ccccccc}
a&g&f&e&d&c&b\\
b&a&g&f&e&d&c\\
c&b&a&g&f&e&d\\
d&c&b&a&g&f&e\\
e&d&c&b&a&g&f\\
f&e&d&c&b&a&g\\
g&f&e&d&c&b&a\\
\end{array}
\right),\ \ \mathbf{L_5}=\left(
\begin{array}{ccccccc}
a&f&d&b&g&e&c\\
d&b&g&e&c&a&f\\
g&e&c&a&f&d&b\\
c&a&f&d&b&g&e\\
f&d&b&g&e&c&a\\
b&g&e&c&a&f&d\\
e&c&a&f&d&b&g\\
\end{array}
\right)$,

$\mathbf{L_6}=\left(
\begin{array}{ccccccc}
a&e&b&f&c&g&d\\
e&b&f&c&g&d&a\\
b&f&c&g&d&a&e\\
f&c&g&d&a&e&b\\
c&g&d&a&e&b&f\\
g&d&a&e&b&f&c\\
d&a&e&b&f&c&g\\
\end{array}
\right)$.\end{center}

{\bf Example 4.5 } Let \begin{center}$\mathbf{T}=\left(
\begin{array}{cccccccc}
1 &2 &3 &4 &5 &6 &7 &8 \\
9 &10&11&12&13&14&15&16\\
17&18&19&20&21&22&23&24\\
25&26&27&28&29&30&31&32\\
33&34&35&36&37&38&39&40\\
41&42&43&44&45&46&47&48\\
49&50&51&52&53&54&55&56\\
57&58&59&60&61&62&63&64\\
\end{array}
\right)$.

$\mathbf{L_{1}}=\left(
\begin{array}{cccccccc}
a&f&h&c&d&g&e&b\\
e&b&d&g&h&c&a&f\\
f&a&c&h&g&d&b&e\\
b&e&g&d&c&h&f&a\\
h&c&a&f&e&b&d&g\\
d&g&e&b&a&f&h&c\\
c&h&f&a&b&e&g&d\\
g&d&b&e&f&a&c&h\\
\end{array}
\right), \ \ \mathbf{L_2}=\left(
\begin{array}{cccccccc}
a&e&f&b&h&d&c&g\\
f&b&a&e&c&g&h&d\\
h&d&c&g&a&e&f&b\\
c&g&h&d&f&b&a&e\\
d&h&g&c&e&a&b&f\\
g&c&d&h&b&f&e&a\\
e&a&b&f&d&h&g&c\\
b&f&e&a&g&c&d&h\\
\end{array}
\right)$,

$\mathbf{L_3}=\left(
\begin{array}{cccccccc}
a&g&b&h&c&e&d&f\\
h&b&g&a&f&d&e&c\\
d&f&c&e&b&h&a&g\\
e&c&f&d&g&a&h&b\\
g&a&h&b&e&c&f&d\\
b&h&a&g&d&f&c&e\\
f&d&e&c&h&b&g&a\\
c&e&d&f&a&g&b&h\\
\end{array}
\right)$,\ \ $\mathbf{L_4}=\left(
\begin{array}{cccccccc}
a&h&d&e&g&b&f&c\\
g&b&f&c&a&h&d&e\\
b&g&c&f&h&a&e&d\\
h&a&e&d&b&g&c&f\\
c&f&b&g&e&d&h&a\\
e&d&h&a&c&f&b&g\\
d&e&a&h&f&c&g&b\\
f&c&g&b&d&e&a&h\\
\end{array}
\right)$.
\end{center}
Here $POL(8,4)=\{\mathbf{L_{1}, L_{2}, L_{3}, L_{4}}\}$.

Taking $s=8, w=1,2,3,4$ in Theorem 3.12, we apply sixty-four treatments of $\mathbf{T}$ and the
above $POL(8,4)$ to construct several $\mathbf{L_{w+2}(8)}$
association schemes. By Definition 2.16, we see that several
pseudo-$L^{*}_{[8-1-w]}(8)$ association schemes can be
induced by several $\mathbf{L_{w+2}(8)}$ association schemes, moreover,
they are several pseudo-$L_{8-1-w}(8)$ association schemes. Now we
research how $8^2$ treatments of several pseudo-$L^{*}_{[8-1-w]}(8)$
association schemes, which are also those of several $\mathbf{L_{w+2}(8)}$
association schemes, can be distributed onto an $8\times 8$ matrix.

{\bf Case 1: w=1. } Utilize $\mathbf{L_{1}}$, one $\mathbf{L_{3}(8)}$
association scheme can be constructed, thus we acquire one inducing
pseudo-$L^{*}_{[8-1-1]}(8)$ association scheme.

At first, we pitch on the treatment $1$ as an example. By
Definition 2.7, there are forty-two treatments which are first
associates of $1$ in the pseudo-$L^{*}_{[8-1-1]}(8)$ association
scheme, that is, they are exactly second associates of $1$
in the $\mathbf{L_{3}(8)}$ association scheme, such that they and $1$ are together located in
$$\mathbf{\Delta^{T; L_{1}}_{1}}=\left(
\begin{array}{cccccccc}
1&  &  &  &  &  &  &  \\
 &10&11&12&13&14&  &16\\
 &  &19&20&21&22&23&24\\
 &26&27&28&29&30&31&  \\
 &34&  &36&37&38&39&40\\
 &42&43&44&  &46&47&48\\
 &50&51&  &53&54&55&56\\
 &58&59&60&61&  &63&64\\
\end{array}
\right).$$

By Definitions 2.8 and 2.9, we will write the collection of all the transversals of
$\mathbf{L_{1}}$ which contain the treatment $1$, i.e.,
$\mathbf{C^{T; L_{1}}_{1}}=\mathbf{C^{T; L_{1}}_{1--10}}\cup
\mathbf{C^{T; L_{1}}_{1--11}}\cup\mathbf{C^{T; L_{1}}_{1--12}}\cup
\mathbf{C^{T; L_{1}}_{1--13}}\cup \mathbf{C^{T; L_{1}}_{1--14}}\cup
\mathbf{C^{T; L_{1}}_{1--16}}$. Where \begin{center}
$\begin{array}{cccccccccccccccccccccc}
\mathbf{C^{T; L_{1}}_{1--10}}&&&&&&&&&&&&&&&&&&&\mathbf{C^{T; L_{1}}_{1--11}}\\
\end{array}$

=$\{$(1) 1--10--19--28--37--46--55--64.
\;\;\;\;\;\;\;\;\;\;\;\;\;\;\;\;\;\;\;\;\;\;\;\; =$\{$(1)
1--11--21--31--38--48--50--60.

(2) 1--10--19--28--40--47--54--61.
\;\;\;\;\;\;\;\;\;\;\;\;\;\;\;\;\;\;\;\;\;\;\;\; (2)
1--11--24--30--36--42--53--63.

(3) 1--10--20--27--37--46--56--63.
\;\;\;\;\;\;\;\;\;\;\;\;\;\;\;\;\;\;\;\;\;\;\;\; (3)
1--11--21--31--34--44--54--64.

(4) 1--10--20--27--39--48--54--61.
\;\;\;\;\;\;\;\;\;\;\;\;\;\;\;\;\;\;\;\;\;\;\;\; (4)
1--11--23--29--40--46--50--60.

(5) 1--10--21--30--36--43--56--63.
\;\;\;\;\;\;\;\;\;\;\;\;\;\;\;\;\;\;\;\;\;\;\;\; (5)
1--11--20--26--38--48--55--61.

(6) 1--10--21--30--39--48--51--60.
\;\;\;\;\;\;\;\;\;\;\;\;\;\;\;\;\;\;\;\;\;\;\;\; (6)
1--11--23--29--36--42--54--64.

(7) 1--10--22--29--36--43--55--64.
\;\;\;\;\;\;\;\;\;\;\;\;\;\;\;\;\;\;\;\;\;\;\;\; (7)
1--11--20--26--40--46--53--63.

(8) 1--10--22--29--40--47--51--60.$\}$,
\;\;\;\;\;\;\;\;\;\;\;\;\;\;\;\;\;\;\;\;\;\;\; (8)
1--11--24--30--34--44--55--61.$\}$,

$\begin{array}{cccccccccccccccccccccc}
\mathbf{C^{T; L_{1}}_{1--12}}&&&&&&&&&&&&&&&&&&&\mathbf{C^{T; L_{1}}_{1--13}}\\
\end{array}$

=$\{$(1) 1--12--23--30--34--43--56--61.
\;\;\;\;\;\;\;\;\;\;\;\;\;\;\;\;\;\;\;\;\;\;\;\; =$\{$(1)
1--13--22--26--40--44--51--63.

(2) 1--12--22--31--37--48--50--59.
\;\;\;\;\;\;\;\;\;\;\;\;\;\;\;\;\;\;\;\;\;\;\;\; (2)
1--13--23--27--34--46--56--60.

(3) 1--12--24--29--38--47--51--58.
\;\;\;\;\;\;\;\;\;\;\;\;\;\;\;\;\;\;\;\;\;\;\;\; (3)
1--13--22--26--36--48--55--59.

(4) 1--12--22--31--34--43--53--64.
\;\;\;\;\;\;\;\;\;\;\;\;\;\;\;\;\;\;\;\;\;\;\;\; (4)
1--13--24--28--38--42--51--63.

(5) 1--12--24--29--39--46--50--59.
\;\;\;\;\;\;\;\;\;\;\;\;\;\;\;\;\;\;\;\;\;\;\;\; (5)
1--13--19--31--40--44--54--58.

(6) 1--12--19--26--38--47--56--61.
\;\;\;\;\;\;\;\;\;\;\;\;\;\;\;\;\;\;\;\;\;\;\;\; (6)
1--13--24--28--34--46--55--59.

(7) 1--12--23--30--37--48--51--58.
\;\;\;\;\;\;\;\;\;\;\;\;\;\;\;\;\;\;\;\;\;\;\;\; (7)
1--13--19--31--38--42--56--60.

(8) 1--12--19--26--39--46--53--64.$\}$,
\;\;\;\;\;\;\;\;\;\;\;\;\;\;\;\;\;\;\;\;\;\;\; (8)
1--13--23--27--36--48--54--58.$\}$,

$\begin{array}{cccccccccccccccccccccc}
\mathbf{C^{T; L_{1}}_{1--14}}&&&&&&&&&&&&&&&&&&&\mathbf{C^{T; L_{1}}_{1--16}}\\
\end{array}$

=$\{$(1) 1--14--24--27--36--47--53--58.
\;\;\;\;\;\;\;\;\;\;\;\;\;\;\;\;\;\;\;\;\;\;\;\; =$\{$(1)
1--16--20--29--39--42--54--59.

(2) 1--14--21--26--39--44--51--64.
\;\;\;\;\;\;\;\;\;\;\;\;\;\;\;\;\;\;\;\;\;\;\;\; (2)
1--16--20--29--38--43--55--58.

(3) 1--14--23--28--40--43--50--61.
\;\;\;\;\;\;\;\;\;\;\;\;\;\;\;\;\;\;\;\;\;\;\;\; (3)
1--16--19--30--39--42--53--60.

(4) 1--14--21--26--36--47--56--59.
\;\;\;\;\;\;\;\;\;\;\;\;\;\;\;\;\;\;\;\;\;\;\;\; (4)
1--16--19--30--37--44--55--58.

(5) 1--14--23--28--37--42--51--64.
\;\;\;\;\;\;\;\;\;\;\;\;\;\;\;\;\;\;\;\;\;\;\;\; (5)
1--16--22--27--34--47--53--60.

(6) 1--14--20--31--40--43--53--58.
\;\;\;\;\;\;\;\;\;\;\;\;\;\;\;\;\;\;\;\;\;\;\;\; (6)
1--16--22--27--37--44--50--63.

(7) 1--14--24--27--39--44--50--61.
\;\;\;\;\;\;\;\;\;\;\;\;\;\;\;\;\;\;\;\;\;\;\;\; (7)
1--16--21--28--34--47--54--59.

(8) 1--14--20--31--37--42--56--59.$\}$,
\;\;\;\;\;\;\;\;\;\;\;\;\;\;\;\;\;\;\;\;\;\;\; (8)
1--16--21--28--38--43--50--63.$\}$.
\end{center}

Subsequently, we deliberate the combination of those forty-two treatments
mentioned. They can be equally divided
into six pairwise disjoint sets, so that arbitrary two different
treatments inside each of these six sets are first associates of
the pseudo-$L^{*}_{[8-1-1]}(8)$ association scheme.
Furthermore, there are eight distinct types of these six sets, exhibited
in \begin{center}Type $u \left\{\begin{array}{l}
\hbox{(u) of $\mathbf{C^{T; L_{1}}_{1--10}}$}\\
\hbox{(u) of $\mathbf{C^{T; L_{1}}_{1--11}}$}\\
\hbox{(u) of $\mathbf{C^{T; L_{1}}_{1--12}}$}\\
\hbox{(u) of $\mathbf{C^{T; L_{1}}_{1--13}}$}\\
\hbox{(u) of $\mathbf{C^{T; L_{1}}_{1--14}}$}\\
\hbox{(u) of $\mathbf{C^{T; L_{1}}_{1--16}}$}\\
\end{array}\right.$, here $u=1,2,3,4,5,6,7,8$. \end{center}

Select the remaining treatments except $1$ from the pseudo-$L^{*}_{[8-1-1]}(8)$
association scheme, those corresponding results are obviously obtained in a similar manner.

It is seen from $\mathbf{C^{T; L_{1}}_{1--10}}$: there exist two distinct treatments, which are first
associates of the pseudo-$L^{*}_{[8-1-1]}(8)$ association scheme,
such that the transversal of $\mathbf{L_{1}}$ containing them is non-unique.

In order to illustrate if the inducing pseudo-$L^{*}_{[8-1-1]}(8)$ association scheme
is an $L_{6}(8)$ association scheme, we shall choose two distinct types matched from eight ones
concerning every treatment. For instance, as for $1$, we shall pick Types $1$ and $2$ up.

(a) When the following confining conditions are satisfied, its sixty-four
treatments can form two distinct $L_{6}(8)$ association schemes.

After calculating $\mathbf{C^{T; L_{1}}_{t}}$, we are able to screen
ninety-six distinct transversals of $\mathbf{L_{1}}$ out from much more than
these ones of the pseudo-$L^{*}_{[8-1-1]}(8)$ association scheme,
so that they fall into twelve paralleling classifications of eight transversals each.

In the circumstances, each transversal $B^{*1}_m$ and each transversal
$B^{*2}_n$ are every row and every
column of $\mathbf{T_{1}^{*}}$, respectively, $m, n=1,\cdots,8$, these obtain
Classifications 1 and 2; each transversal $B^{*(k+2)}_m$ is the set of eight
treatments of $\mathbf{T_{1}^{*}}$ altogether corresponding to the same symbol
of $\mathbf{L^{*1}_k}$, this obtains Classification ($k+2$), $k=1,2,3,4$.
Moreover, two transversals of different classifications
have exactly one common treatment. Afterwards,
each transversal $B^{*7}_m$ and each transversal $B^{*8}_n$ are every row and every
column of $\mathbf{T_{2}^{*}}$, respectively, $m, n=1,\cdots,8$, these obtain
Classifications 7 and 8; each transversal $B^{*(k+8)}_m$ is the set of eight
treatments of $\mathbf{T_{2}^{*}}$ altogether corresponding to the same symbol
of $\mathbf{L^{*2}_k}$, this obtains Classification ($k+8$), $k=1,2,3,4$.
Meanwhile, two transversals of distinct classifications
have exactly one common treatment. Where\begin{center}
$\mathbf{T_{1}^{*}}=\left(
\begin{array}{cccccccc}
1 &10&19&28&37&46&55&64\\
20&27&2 &9 &56&63&38&45\\
39&48&53&62&3 &12&17&26\\
54&61&40&47&18&25&4 &11\\
42&33&60&51&14&5 &32&23\\
59&52&41&34&31&24&13&6 \\
16&7 &30&21&44&35&58&49\\
29&22&15&8 &57&50&43&36\\
\end{array}
\right)$.

$\mathbf{L^{*1}_1}$=$\left(
\begin{array}{cccccccc}
a^{*}&e^{*}&f^{*}&b^{*}&h^{*}&d^{*}&c^{*}&g^{*}\\
e^{*}&a^{*}&b^{*}&f^{*}&d^{*}&h^{*}&g^{*}&c^{*}\\
f^{*}&b^{*}&a^{*}&e^{*}&c^{*}&g^{*}&h^{*}&d^{*}\\
b^{*}&f^{*}&e^{*}&a^{*}&g^{*}&c^{*}&d^{*}&h^{*}\\
h^{*}&d^{*}&c^{*}&g^{*}&a^{*}&e^{*}&f^{*}&b^{*}\\
d^{*}&h^{*}&g^{*}&c^{*}&e^{*}&a^{*}&b^{*}&f^{*}\\
c^{*}&g^{*}&h^{*}&d^{*}&f^{*}&b^{*}&a^{*}&e^{*}\\
g^{*}&c^{*}&d^{*}&h^{*}&b^{*}&f^{*}&e^{*}&a^{*}\\
\end{array}
\right),\ \ \mathbf{L^{*1}_2}$=$\left(
\begin{array}{cccccccc}
a^{*}&f^{*}&h^{*}&c^{*}&d^{*}&g^{*}&e^{*}&b^{*}\\
g^{*}&d^{*}&b^{*}&e^{*}&f^{*}&a^{*}&c^{*}&h^{*}\\
b^{*}&e^{*}&g^{*}&d^{*}&c^{*}&h^{*}&f^{*}&a^{*}\\
h^{*}&c^{*}&a^{*}&f^{*}&e^{*}&b^{*}&d^{*}&g^{*}\\
c^{*}&h^{*}&f^{*}&a^{*}&b^{*}&e^{*}&g^{*}&d^{*}\\
e^{*}&b^{*}&d^{*}&g^{*}&h^{*}&c^{*}&a^{*}&f^{*}\\
d^{*}&g^{*}&e^{*}&b^{*}&a^{*}&f^{*}&h^{*}&c^{*}\\
f^{*}&a^{*}&c^{*}&h^{*}&g^{*}&d^{*}&b^{*}&e^{*}\\
\end{array}
\right)$,

$\mathbf{L^{*1}_3}$=$\left(
\begin{array}{cccccccc}
a^{*}&c^{*}&e^{*}&g^{*}&f^{*}&h^{*}&b^{*}&d^{*}\\
f^{*}&h^{*}&b^{*}&d^{*}&a^{*}&c^{*}&e^{*}&g^{*}\\
h^{*}&f^{*}&d^{*}&b^{*}&c^{*}&a^{*}&g^{*}&e^{*}\\
c^{*}&a^{*}&g^{*}&e^{*}&h^{*}&f^{*}&d^{*}&b^{*}\\
d^{*}&b^{*}&h^{*}&f^{*}&g^{*}&e^{*}&c^{*}&a^{*}\\
g^{*}&e^{*}&c^{*}&a^{*}&d^{*}&b^{*}&h^{*}&f^{*}\\
e^{*}&g^{*}&a^{*}&c^{*}&b^{*}&d^{*}&f^{*}&h^{*}\\
b^{*}&d^{*}&f^{*}&h^{*}&e^{*}&g^{*}&a^{*}&c^{*}\\
\end{array}
\right),\ \ \mathbf{L^{*1}_4}$=$\left(
\begin{array}{cccccccc}
a^{*}&d^{*}&g^{*}&f^{*}&b^{*}&c^{*}&h^{*}&e^{*}\\
h^{*}&e^{*}&b^{*}&c^{*}&g^{*}&f^{*}&a^{*}&d^{*}\\
d^{*}&a^{*}&f^{*}&g^{*}&c^{*}&b^{*}&e^{*}&h^{*}\\
e^{*}&h^{*}&c^{*}&b^{*}&f^{*}&g^{*}&d^{*}&a^{*}\\
g^{*}&f^{*}&a^{*}&d^{*}&h^{*}&e^{*}&b^{*}&c^{*}\\
b^{*}&c^{*}&h^{*}&e^{*}&a^{*}&d^{*}&g^{*}&f^{*}\\
f^{*}&g^{*}&d^{*}&a^{*}&e^{*}&h^{*}&c^{*}&b^{*}\\
c^{*}&b^{*}&e^{*}&h^{*}&d^{*}&a^{*}&f^{*}&g^{*}\\
\end{array}
\right)$.
\end{center}
It is very obvious that $\mathbf{L^{*1}_1}, \mathbf{L^{*1}_2},
\mathbf{L^{*1}_3}, \mathbf{L^{*1}_4}$ are mutually orthogonal.
\begin{center}
$\mathbf{T_{2}^{*}}=\left(
\begin{array}{cccccccc}
1 &10&19&28&54&61&40&47\\
20&27&2 &9 &39&48&53&62\\
58&49&44&35&13&6 &31&24\\
43&36&57&50&32&23&14&5 \\
55&64&37&46&4 &11&18&25\\
38&45&56&63&17&26&3 &12\\
16&7 &30&21&59&52&41&34\\
29&22&15&8 &42&33&60&51\\
\end{array}
\right)$.

$\mathbf{L^{*2}_1}$=$\left(
\begin{array}{cccccccc}
a^{*}&e^{*}&g^{*}&c^{*}&h^{*}&d^{*}&b^{*}&f^{*}\\
h^{*}&d^{*}&b^{*}&f^{*}&a^{*}&e^{*}&g^{*}&c^{*}\\
g^{*}&c^{*}&a^{*}&e^{*}&b^{*}&f^{*}&h^{*}&d^{*}\\
b^{*}&f^{*}&h^{*}&d^{*}&g^{*}&c^{*}&a^{*}&e^{*}\\
e^{*}&a^{*}&c^{*}&g^{*}&d^{*}&h^{*}&f^{*}&b^{*}\\
d^{*}&h^{*}&f^{*}&b^{*}&e^{*}&a^{*}&c^{*}&g^{*}\\
c^{*}&g^{*}&e^{*}&a^{*}&f^{*}&b^{*}&d^{*}&h^{*}\\
f^{*}&b^{*}&d^{*}&h^{*}&c^{*}&g^{*}&e^{*}&a^{*}\\
\end{array}
\right),\ \ \mathbf{L^{*2}_2}$=$\left(
\begin{array}{cccccccc}
a^{*}&f^{*}&e^{*}&b^{*}&c^{*}&h^{*}&g^{*}&d^{*}\\
f^{*}&a^{*}&b^{*}&e^{*}&h^{*}&c^{*}&d^{*}&g^{*}\\
c^{*}&h^{*}&g^{*}&d^{*}&a^{*}&f^{*}&e^{*}&b^{*}\\
h^{*}&c^{*}&d^{*}&g^{*}&f^{*}&a^{*}&b^{*}&e^{*}\\
b^{*}&e^{*}&f^{*}&a^{*}&d^{*}&g^{*}&h^{*}&c^{*}\\
e^{*}&b^{*}&a^{*}&f^{*}&g^{*}&d^{*}&c^{*}&h^{*}\\
d^{*}&g^{*}&h^{*}&c^{*}&b^{*}&e^{*}&f^{*}&a^{*}\\
g^{*}&d^{*}&c^{*}&h^{*}&e^{*}&b^{*}&a^{*}&f^{*}\\
\end{array}
\right)$,

$\mathbf{L^{*2}_3}$=$\left(
\begin{array}{cccccccc}
a^{*}&c^{*}&h^{*}&f^{*}&e^{*}&g^{*}&d^{*}&b^{*}\\
g^{*}&e^{*}&b^{*}&d^{*}&c^{*}&a^{*}&f^{*}&h^{*}\\
d^{*}&b^{*}&e^{*}&g^{*}&h^{*}&f^{*}&a^{*}&c^{*}\\
f^{*}&h^{*}&c^{*}&a^{*}&b^{*}&d^{*}&g^{*}&e^{*}\\
h^{*}&f^{*}&a^{*}&c^{*}&d^{*}&b^{*}&e^{*}&g^{*}\\
b^{*}&d^{*}&g^{*}&e^{*}&f^{*}&h^{*}&c^{*}&a^{*}\\
e^{*}&g^{*}&d^{*}&b^{*}&a^{*}&c^{*}&h^{*}&f^{*}\\
c^{*}&a^{*}&f^{*}&h^{*}&g^{*}&e^{*}&b^{*}&d^{*}\\
\end{array}
\right),\ \ \mathbf{L^{*2}_4}$=$\left(
\begin{array}{cccccccc}
a^{*}&d^{*}&f^{*}&g^{*}&b^{*}&c^{*}&e^{*}&h^{*}\\
e^{*}&h^{*}&b^{*}&c^{*}&f^{*}&g^{*}&a^{*}&d^{*}\\
h^{*}&e^{*}&c^{*}&b^{*}&g^{*}&f^{*}&d^{*}&a^{*}\\
d^{*}&a^{*}&g^{*}&f^{*}&c^{*}&b^{*}&h^{*}&e^{*}\\
c^{*}&b^{*}&h^{*}&e^{*}&d^{*}&a^{*}&g^{*}&f^{*}\\
g^{*}&f^{*}&d^{*}&a^{*}&h^{*}&e^{*}&c^{*}&b^{*}\\
f^{*}&g^{*}&a^{*}&d^{*}&e^{*}&h^{*}&b^{*}&c^{*}\\
b^{*}&c^{*}&e^{*}&h^{*}&a^{*}&d^{*}&f^{*}&g^{*}\\
\end{array}
\right)$.
\end{center}
It is very apparent that $\mathbf{L^{*2}_1}, \mathbf{L^{*2}_2},
\mathbf{L^{*2}_3}, \mathbf{L^{*2}_4}$ are mutually orthogonal.

In sum, $8^2$ treatments of the pseudo-$L^{*}_{[8-1-1]}(8)$ association
scheme can be exactly divided into eight pairwise parallel transversals
$B^{*i}_1, \cdots, B^{*i}_8$ inside the $i$-th classification, here
$i=6v-5,6v-4,6v-3,6v-2,6v-1,6v; v=1,2$.
By Consequence I of Theorem 3.6, it may be verified that its $8^2$ treatments
can form two distinct $L_{6}(8)$ association schemes.

(b) When anyone of five conditions of Theorem 3.5 is established,
the pseudo-$L^{*}_{[8-1-1]}(8)$ association scheme is not an $L_{6}(8)$ association scheme.

It is apparently understood that $\mathbf{(I)}$ and $\mathbf{(III)}$ have already been
established in Theorem 3.5. We observe $B^{*1}_1\cap B^{*7}_1=\{$1--10--19--28--37--46--55--64.$\}\cap \{$1--10--19--28--40--47--54--61.$\}=\{1, 10, 19, 28\}$.
Let $t^{*}_2=1$, none of (1) in $\mathbf{C^{T; L_{1}}_{1--10}}$, (1) in
$\mathbf{C^{T; L_{1}}_{1--11}}$, (1) in
$\mathbf{C^{T; L_{1}}_{1--12}}$, (1) in
$\mathbf{C^{T; L_{1}}_{1--13}}$, (1) in
$\mathbf{C^{T; L_{1}}_{1--14}}$, (1) in
$\mathbf{C^{T; L_{1}}_{1--16}}$ is equal
to $B^{*7}_1$. Hence these clarify Conditions $\mathbf{(II)}$ and $\mathbf{(IV)}$
of Theorem 3.5 are established. Put $t^{*}_1=10$ and $t^{*}_2=1$, it is known that $(\{10\}\cup\{$7--22--27--33--48--52--61$\})\cap (\{1\}\cup\{$12--19--26--39--46--53--64$\})=\emptyset$
and $(\{10\}\cup\{$7--22--27--33--48--52--61$\})\cap (\{1\}\cup\{$13--24--28--34--46--55--59\}$)=\emptyset$
hold; and it is discovered from $\mathbf{C^{T; L_{1}, L_{2}}_{1--10}}$ that
the transversal of $\mathbf{L_{1}}$ containing $10$ and $1$ is not sole one, that is,
the set of cardinality $8$ containing them is non-unique.
It follows that $\mathbf{(V)}$ holds. According to Theorem 3.5, but the
inducing pseudo-$L^{*}_{[8-1-1]}(8)$ association scheme is not an $L_{6}(8)$ association scheme.

Even though two treatments $37$ and $40$ are second associates of the
pseudo-$L^{*}_{[8-1-1]}(8)$ association scheme, we apparently have
$(\{37\}\cup\{$1--10--20--27--46--56--63$\})\cap (\{40\}\cup\{$1--11--20--26--46--53--63\}$)
=\{1, 20, 46, 63\}$.

{\bf Case 2: w=2. } Utilize $\mathbf{L_{1}, L_{2}}$, one $\mathbf{L_{4}(8)}$
association scheme can be constructed, thus we obtain one inducing
pseudo-$L^{*}_{[8-1-2]}(8)$ association scheme.

Taking any treatment $t$, we may compute $\mathbf{C^{T; L_{1},
L_{2}}_{t}}$. Now we select the treatment $1$ as an example. There are
thirty-five treatments which are first associates of $1$ in the
pseudo-$L^{*}_{[8-1-2]}(8)$ association scheme, such that they and
$1$ are together located in $$\mathbf{\Delta^{T; L_{1}, L_{2}}_{1}}=\left(
\begin{array}{cccccccc}
1&  &  &  &  &  &  &  \\
 &10&  &12&13&14&  &16\\
 &  &19&20&  &22&23&24\\
 &26&27&28&29&30&  &  \\
 &34&  &36&37&  &39&40\\
 &42&43&44&  &46&47&  \\
 &  &51&  &53&54&55&56\\
 &58&59&  &61&  &63&64\\
\end{array}
\right).$$

Next we display the collection of all the common transversals of
$\mathbf{L_{1}, L_{2}}$ containing $1$, i.e.,
$\mathbf{C^{T; L_{1}, L_{2}}_{1}}=\mathbf{C^{T; L_{1},
L_{2}}_{1--10}}\cup \mathbf{C^{T; L_{1}, L_{2}}_{1--12}}\cup
\mathbf{C^{T; L_{1}, L_{2}}_{1--13}}\cup \mathbf{C^{T; L_{1},
L_{2}}_{1--14}}\cup \mathbf{C^{T; L_{1}, L_{2}}_{1--16}}$. Here
\begin{center} $\mathbf{C^{T; L_{1}, L_{2}}_{1--10}}=\{$ (1)
1--10--19--28--37--46--55--64. (2) 1--10--19--28--40--47--54--61.\\
(3) 1--10--20--27--37--46--56--63. (7)
1--10--22--29--36--43--55--64.$\}$,

$\mathbf{C^{T; L_{1}, L_{2}}_{1--12}}=\{$ (1)
1--12--23--30--34--43--56--61. (8)
1--12--19--26--39--46--53--64.$\}$,

$\mathbf{C^{T; L_{1}, L_{2}}_{1--13}}=\{$ (1)
1--13--22--26--40--44--51--63. (6)
1--13--24--28--34--46--55--59.$\}$,

$\mathbf{C^{T; L_{1}, L_{2}}_{1--14}}=\{$ (1)
1--14--24--27--36--47--53--58. (5)
1--14--23--28--37--42--51--64.$\}$,

$\mathbf{C^{T; L_{1}, L_{2}}_{1--16}}=\{$ (1)
1--16--20--29--39--42--54--59. (4)
1--16--19--30--37--44--55--58.$\}$.
\end{center}

Carefully filtrate, it is grasped that these thirty-five treatments can be
equally divided into five pairwise disjoint sets, so that arbitrary two different
treatments inside each of these five sets are first associates
of the pseudo-$L^{*}_{[8-1-2]}(8)$ association scheme. Such as
\begin{center}Type $1 \left\{\begin{array}{l}
\hbox{(1) of $\mathbf{C^{T; L_{1}, L_{2}}_{1--10}}$}\\
\hbox{(1) of $\mathbf{C^{T; L_{1}, L_{2}}_{1--12}}$}\\
\hbox{(1) of $\mathbf{C^{T; L_{1}, L_{2}}_{1--13}}$}\\
\hbox{(1) of $\mathbf{C^{T; L_{1}, L_{2}}_{1--14}}$}\\
\hbox{(1) of $\mathbf{C^{T; L_{1}, L_{2}}_{1--16}}$}\\
\end{array}\right.$. \end{center}

Similarly to $1$, single the remainders out from the pseudo-$L^{*}_{[8-1-2]}(8)$
association scheme, those similar consequences are apparently deduced.

It is known that there exist two treatments $1$ and $10$ of
$\mathbf{C^{T; L_{1}, L_{2}}_{1--10}}$ such that the common transversal
of $\mathbf{L_{1}, L_{2}}$ containing them is non-unique.

Percolate five pairwise disjoint sets concerning each treatment,
the next step is to examine whether the inducing pseudo-$L^{*}_{[8-1-2]}(8)$
association scheme is an $L_{5}(8)$ association scheme.

(a) When the following confining conditions are satisfied, its sixty-four
treatments can form one $L_{5}(8)$ association scheme.

After calculating $\mathbf{C^{T; L_{1}, L_{2}}_{t}}$, we are able to filtrate
forty distinct common transversals of $\mathbf{L_{1}, L_{2}}$ from more than
these ones of the pseudo-$L^{*}_{[8-1-2]}(8)$ association scheme,
so that they fall into five paralleling classifications of eight transversals each.
In this case, eight pairwise parallel transversals $B^{*i}_1, \cdots, B^{*i}_8$
are identical to those of (a) within {\bf Case 1}, $i=1,2,3,4,5$.

In short, $8^2$ treatments of the pseudo-$L^{*}_{[8-1-2]}(8)$
association scheme can be exactly divided into eight pairwise
parallel transversals $B^{*i}_1, \cdots, B^{*i}_8$ inside the $i$-th
classification, $i=1,2,3,4,5$. Furthermore, two transversals of different classifications
have exactly one common treatment. By Consequence I of Theorem 3.6,
it is elucidated that its $8^2$ treatments can form one $L_{5}(8)$ association scheme.

(b) When $\mathbf{(V)}$ of Theorem 3.5 is established,
the pseudo-$L^{*}_{[8-1-2]}(8)$ association scheme is not an $L_{5}(8)$ association scheme.

Still take $t^{*}_1=10$ and $t^{*}_2=1$, since the common transversal
of $\mathbf{L_{1}, L_{2}}$ containing $10$ and $1$ is non-unique;
$\{10\}\cup\{$7--22--27--33--48--52--61$\}, \{1\}\cup\{$12--19--26--39--46--53--64$\}$ and $\{1\}\cup\{$13--24--28--34--46--55--59$\}$ are also common transversals of
$\mathbf{L_{1}, L_{2}}$. It follows that $\mathbf{(V)}$ still holds.
Because of Theorem 3.5, the inducing pseudo-$L^{*}_{[8-1-2]}(8)$ association
scheme is not an $L_{5}(8)$ association scheme.

{\bf Case 3: w=3. } Utilize $\mathbf{L_{1}, L_{2}, L_{3}}$, one $\mathbf{L_{5}(8)}$
association scheme can be constructed, thus we obtain one inducing
pseudo-$L^{*}_{[8-1-3]}(8)$ association scheme.

Taking any treatment $t$, we firstly calculate $\mathbf{C^{T; L_{1},
L_{2}, L_{3}}_{t}}$. We single $t_1=1$ out as an instance. There are
twenty-eight treatments which are first associates of $1$ in the
pseudo-$L^{*}_{[8-1-3]}(8)$ association scheme, such that they and
$1$ are together located in $\mathbf{\Delta^{T; L_{1}, L_{2}, L_{3}}_{1}}$.
Next we show the collection of all the
common transversals of $\mathbf{L_{1}, L_{2}, L_{3}}$ containing
$1$, i.e., $\mathbf{C^{T; L_{1}, L_{2},
L_{3}}_{1}}=\mathbf{C^{T; L_{1}, L_{2}, L_{3}}_{1--10}}\cup
\mathbf{C^{T; L_{1}, L_{2}, L_{3}}_{1--13}}\cup \mathbf{C^{T; L_{1},
L_{2}, L_{3}}_{1--14}}\cup \mathbf{C^{T; L_{1}, L_{2},
L_{3}}_{1--16}}$. Here \begin{center}$\mathbf{C^{T; L_{1}, L_{2},
L_{3}}_{1--10}}=\{$(1) 1--10--19--28--37--46--55--64.$\}$,

$\mathbf{C^{T; L_{1}, L_{2}, L_{3}}_{1--13}}=\{$(1)
1--13--22--26--40--44--51--63.$\}$,

$\mathbf{C^{T; L_{1}, L_{2}, L_{3}}_{1--14}}=\{$(1)
1--14--24--27--36--47--53--58.$\}$,

$\mathbf{C^{T; L_{1}, L_{2}, L_{3}}_{1--16}}=\{$(1)
1--16--20--29--39--42--54--59.$\}$.
\end{center}

We set
$$\begin{array}{ccc}
&\nearrow&A^{1}_1=\{10,19,28,37,46,55,64\}\\
t_1=1&\rightarrow&A^{1}_2=\{13,22,26,40,44,51,63\}\\
&\searrow&A^{1}_3=\{14,24,27,36,47,53,58\}\\
&\searrow&A^{1}_4=\{16,20,29,39,42,54,59\}
\end{array}$$
in Consequence III of Theorem 3.6.

It is stressed that these twenty-eight treatments can be
equally divided into four pairwise disjoint sets, so that arbitrary two different
treatments inside each of these four sets are first associates
of the pseudo-$L^{*}_{[8-1-3]}(8)$ association scheme.

Choosing the remainders except $1$ from the pseudo-$L^{*}_{[8-1-3]}(8)$
association scheme, we detect those similar conclusions are also inferred.

Calculating $\mathbf{C^{T; L_{1}, L_{2}, L_{3}}_{t}}$, there are totally
thirty-two distinct common transversals of $\mathbf{L_{1}, L_{2}, L_{3}}$,
such that they fall into four paralleling classifications of eight transversals each.
In the circumstances, Classifications 1 to 4 are identical to those of {\bf Case 2}.

As a whole, $8^2$ treatments of the pseudo-$L^{*}_{[8-1-3]}(8)$
association scheme can be exactly divided into eight pairwise
parallel transversals $B^{*i}_1, \cdots, B^{*i}_8$ inside the $i$-th
classification, $i=1,2,3,4$. Moreover, two transversals of different classifications
have exactly one common treatment. By Consequence I of Theorem 3.6,
it is deduced that its $8^2$ treatments can form one $L_{4}(8)$ association scheme.

With respect to arbitrary two treatments, which are first
associates of the pseudo-$L^{*}_{[8-1-3]}(8)$ association scheme, via computing, it
is exhibited that the common transversal of $\mathbf{L_{1}, L_{2}, L_{3}}$
containing them is sole one, that is, the set of cardinality $8$ containing
them is unique. Further, it may be clarified that
$\mathbf{(v)}$ in Consequence III of Theorem 3.6 can be satisfied.
By Corollary 3.7, it goes without saying that the inducing
pseudo-$L^{*}_{[8-1-3]}(8)$ association scheme must be
an $L_{4}(8)$ association scheme.

{\bf Case 4: w=4. } Utilize $\mathbf{L_{1}, L_{2}, L_{3}, L_{4}}$,
one $\mathbf{L_{6}(8)}$ association scheme can be constructed,
thus we have one inducing pseudo-$L^{*}_{[8-1-4]}(8)$ association scheme.

Choosing $x=1$ in Theorem 3.8, we still count $\mathbf{C^{T; L_{1},
L_{2}, L_{3}, L_{4}}_{1}}$. Now we will exhibit the
collection of all the common transversals of $\mathbf{L_{1}, L_{2}, L_{3}, L_{4}}$
containing the treatment $1$, i.e., $\mathbf{C^{T; L_{1}, L_{2}, L_{3},
L_{4}}_{1}}\\=\mathbf{C^{T; L_{1}, L_{2}, L_{3}, L_{4}}_{1--10}}\cup
\mathbf{C^{T; L_{1}, L_{2}, L_{3}, L_{4}}_{1--14}}\cup \mathbf{C^{T;
L_{1}, L_{2}, L_{3}, L_{4}}_{1--16}}$. Here \begin{center}
$\mathbf{C^{T; L_{1}, L_{2}, L_{3}, L_{4}}_{1--10}}=\{$(1)
1--10--19--28--37--46--55--64.$\}$,

$\mathbf{C^{T; L_{1}, L_{2}, L_{3}, L_{4}}_{1--14}}=\{$(1)
1--14--24--27--36--47--53--58.$\}$,

$\mathbf{C^{T; L_{1}, L_{2}, L_{3}, L_{4}}_{1--16}}=\{$(1)
1--16--20--29--39--42--54--59.$\}$.

Set $\begin{array}{ccc}
&\nearrow&Y=\{10,19,28,37,46,55,64\}\\
x=1&\rightarrow&Z=\{16,20,29,39,42,54,59\} \\
&\searrow&W=\{14,24,27,36,47,53,58\}
\end{array}$ in Theorem 3.8.\end{center}

It is highlighted that there are twenty-one treatments which are first associates of $1$ in the
pseudo-$L^{*}_{[8-1-4]}(8)$ association scheme, such that they can be
equally divided into three pairwise disjoint sets.
Take the remainders except $1$ from the pseudo-$L^{*}_{[8-1-4]}(8)$
association scheme, it may be shown that those similar results are also acquired.

Through calculations, there are entirely twenty-four distinct common transversals of
$\mathbf{L_{1}, L_{2}, L_{3}}$, $\mathbf{L_{4}}$, such that they fall into three
paralleling classifications of eight transversals each.
In this case, Classifications 1 to 3 are identical to those of {\bf Case 3}.

In sum, $8^2$ treatments of the pseudo-$L^{*}_{[8-1-4]}(8)$ association
scheme can be exactly divided into eight pairwise parallel transversals
$B^{*i}_1, \cdots, B^{*i}_8$ inside the $i$-th classification, $i=1,2,3$.
Moreover, two transversals of different classifications
have exactly one common treatment. By Consequence I of Theorem 3.6,
it is judged that its $8^2$ treatments can form one $L_{3}(8)$ association scheme.

Due to $8>(8-1-4-1)^2$, calculating $\mathbf{C^{T; L_{1}, L_{2}, L_{3}, L_{4}}_{t}}$,
we observe that twenty-one treatments, which are first associates of $t$ in the
pseudo-$L^{*}_{[8-1-4]}(8)$ association scheme, can be equally divided
into three pairwise disjoint sets, and pay attention to Theorem 3.9 again.
With respect to arbitrary two distinct treatments, which are first associates of the
pseudo-$L^{*}_{[8-1-4]}(8)$ association scheme, according to Consequence IV of Theorem 3.6 or Theorem 3.8,
it is known that the common transversal of $\mathbf{L_{1}, L_{2}, L_{3}, L_{4}}$
containing them is sole one, that is, the set of cardinality $8$ containing
them is unique. Further, we examine $\mathbf{(v)}$ in Theorem 3.12 can be satisfied.
By Corollary 3.7, thus it is deduced that the inducing pseudo-$L^{*}_{[8-1-4]}(8)$ association
scheme must be an $L_{3}(8)$ association scheme.

By Theorem 3.12, $\mathbf{L_5, L_6, L_7}$ can be added to the
$POL(8,4)=\{\mathbf{L_1, L_2, L_3, L_4}\}$ to obtain a $POL(8,7)$. Here
\begin{center}
$\mathbf{L_5}=\left(
\begin{array}{cccccccc}
a&b&c&d&e&f&g&h\\
b&a&d&c&f&e&h&g\\
c&d&a&b&g&h&e&f\\
d&c&b&a&h&g&f&e\\
e&f&g&h&a&b&c&d\\
f&e&h&g&b&a&d&c\\
g&h&e&f&c&d&a&b\\
h&g&f&e&d&c&b&a\\
\end{array}
\right),\ \ \mathbf{L_6}=\left(
\begin{array}{cccccccc}
a&d&g&f&b&c&h&e\\
c&b&e&h&d&a&f&g\\
e&h&c&b&f&g&d&a\\
g&f&a&d&h&e&b&c\\
f&g&d&a&e&h&c&b\\
h&e&b&c&g&f&a&d\\
b&c&h&e&a&d&g&f\\
d&a&f&g&c&b&e&h\\
\end{array}
\right)$,

$\mathbf{L_7}=\left(
\begin{array}{cccccccc}
a&c&e&g&f&h&b&d\\
d&b&h&f&g&e&c&a\\
g&e&c&a&d&b&h&f\\
f&h&b&d&a&c&e&g\\
b&d&f&h&e&g&a&c\\
c&a&g&e&h&f&d&b\\
h&f&d&b&c&a&g&e\\
e&g&a&c&b&d&f&h\\
\end{array}
\right)$.\end{center}

{\bf Example 4.6 } Let \begin{center}
$\mathbf{T}$=$\left(\begin{array}{ccccccccc}
1&2&3&4&5&6&7&8&9\\
10&11&12&13&14&15&16&17&18\\
19&20&21&22&23&24&25&26&27\\
28&29&30&31&32&33&34&35&36\\
37&38&39&40&41&42&43&44&45\\
46&47&48&49&50&51&52&53&54\\
55&56&57&58&59&60&61&62&63\\
64&65&66&67&68&69&70&71&72\\
73&74&75&76&77&78&79&80&81
\end{array}
\right)$. $\mathbf{L_1}$=$\left(
\begin{array}{ccccccccc}
a&f&h&e&g&c&i&b&d\\
i&b&d&a&f&h&e&g&c\\
e&g&c&i&b&d&a&f&h\\
c&e&g&d&i&b&h&a&f\\
h&a&f&c&e&g&d&i&b\\
d&i&b&h&a&f&c&e&g\\
b&d&i&f&h&a&g&c&e\\
g&c&e&b&d&i&f&h&a\\
f&h&a&g&c&e&b&d&i
\end{array}
\right)$,

$\mathbf{L_2}=\left(
\begin{array}{ccccccccc}
a&i&e&c&h&d&b&g&f\\
f&b&g&e&a&i&d&c&h\\
h&d&c&g&f&b&i&e&a\\
e&a&i&d&c&h&f&b&g\\
g&f&b&i&e&a&h&d&c\\
c&h&d&b&g&f&a&i&e\\
i&e&a&h&d&c&g&f&b\\
b&g&f&a&i&e&c&h&d\\
d&c&h&f&b&g&e&a&i
\end{array}
\right),\ \ \mathbf{L_3}=\left(
\begin{array}{ccccccccc}
a&g&d&f&c&i&h&e&b\\
e&b&h&g&d&a&c&i&f\\
i&f&c&b&h&e&d&a&g\\
b&h&e&d&a&g&i&f&c\\
f&c&i&h&e&b&a&g&d\\
g&d&a&c&i&f&e&b&h\\
c&i&f&e&b&h&g&d&a\\
d&a&g&i&f&c&b&h&e\\
h&e&b&a&g&d&f&c&i
\end{array}
\right)$,

$\mathbf{L_4}=\left(
\begin{array}{ccccccccc}
a&e&i&b&f&g&c&d&h\\
g&b&f&h&c&d&i&a&e\\
d&h&c&e&i&a&f&g&b\\
f&g&b&d&h&c&e&i&a\\
c&d&h&a&e&i&b&f&g\\
i&a&e&g&b&f&h&c&d\\
h&c&d&i&a&e&g&b&f\\
e&i&a&f&g&b&d&h&c\\
b&f&g&c&d&h&a&e&i
\end{array}
\right),\ \ \mathbf{L_5}=\left(
\begin{array}{ccccccccc}
a&h&f&i&d&b&e&c&g\\
d&b&i&c&g&e&h&f&a\\
g&e&c&f&a&h&b&i&d\\
h&f&a&d&b&i&c&g&e\\
b&i&d&g&e&c&f&a&h\\
e&c&g&a&h&f&i&d&b\\
f&a&h&b&i&d&g&e&c\\
i&d&b&e&c&g&a&h&f\\
c&g&e&h&f&a&d&b&i
\end{array}
\right)$.
\end{center}
Here $POL(9,5)=\{\mathbf{L_{1}, L_{2}, L_{3}, L_{4}, L_{5}}\}$.

Taking $s=9, w=1,2,3,4,5$ in Theorem 3.12, we utilize eighty-one treatments of $\mathbf{T}$ and the
above $POL(9,5)$ to construct some $\mathbf{L_{w+2}(9)}$
association schemes. It has already been known that some
pseudo-$L^{*}_{[9-1-w]}(9)$ association schemes can be
induced by some $\mathbf{L_{w+2}(9)}$ association schemes, moreover,
they are some pseudo-$L_{9-1-w}(9)$ association schemes. Now we
probe how $9^2$ treatments of some pseudo-$L^{*}_{[9-1-w]}(9)$
association schemes, which are also those of some $\mathbf{L_{w+2}(9)}$
association schemes, can be distributed onto a $9\times 9$ matrix.

{\bf Case 1: w=5. } Applying $\mathbf{L_{1}, L_{2}, L_{3}, L_{4}, L_{5}}$, we
can construct one $\mathbf{L_{7}(9)}$ association scheme and obtain one
inducing pseudo-$L^{*}_{[9-1-5]}(9)$ association scheme.

Choosing $t_1=1$ in Consequence III of Theorem 3.6, we calculate $\mathbf{C^{T; L_{1},
L_{2}, L_{3}, L_{4}, L_{5}}_{1}}$. By Definition 2.7, there are
twenty-four treatments which are first
associates of $1$ in the pseudo-$L^{*}_{[9-1-5]}(9)$ association
scheme, that is, they are exactly second associates of $1$
in the $\mathbf{L_{7}(9)}$ association scheme, such that they and $1$ are together located in
$$\mathbf{\Delta^{T; L_{1}, L_{2}, L_{3}, L_{4}, L_{5}}_{1}}=\left(
\begin{array}{ccccccccc}
1&  &  &  &  &  &  &  &  \\
 &11&12&  &  &  &16&  &  \\
 &20&21&22&  &  &  &  &  \\
 &  &  &31&  &33&34&  &  \\
 &  &39&  &41&  &  &  &45\\
 &  &  &  &  &51&  &53&54\\
 &  &  &58&  &  &61&62&  \\
 &  &  &  &68&69&  &71&  \\
 &74&  &  &77&  &  &  &81
\end{array}
\right).$$
By Definitions 2.8 and 2.9, we shall write the
collection of all the common transversals of $\mathbf{L_{1}, L_{2}, L_{3}},\\ \mathbf{L_{4}, L_{5}}$
which contain the treatment $1$, i.e., $\mathbf{C^{T; L_{1},
L_{2}, L_{3}, L_{4}, L_{5}}_{1}}=\mathbf{C^{T; L_{1}, L_{2}, L_{3},
L_{4}, L_{5}}_{1--11}}\cup \\\mathbf{C^{T; L_{1}, L_{2}, L_{3}, L_{4},
L_{5}}_{1--12}}\cup \mathbf{C^{T; L_{1}, L_{2}, L_{3}, L_{4},
L_{5}}_{1--16}}$. Here \begin{center}$\mathbf{C^{T; L_{1}, L_{2},
L_{3}, L_{4}, L_{5}}_{1--11}}=\{$(1)
1--11--21--31--41--51--61--71--81.$\}$,

$\mathbf{C^{T; L_{1}, L_{2}, L_{3}, L_{4}, L_{5}}_{1--12}}=\{$(1)
1--12--20--34--45--53--58--69--77.$\}$,

$\mathbf{C^{T; L_{1}, L_{2}, L_{3}, L_{4}, L_{5}}_{1--16}}=\{$(1)
1--16--22--33--39--54--62--68--74.$\}$.

Let $\begin{array}{ccc}
&\nearrow&A^{1}_1=\{11,21,31,41,51,61,71,81\}\\
t_1=1&\rightarrow&A^{1}_2=\{12,20,34,45,53,58,69,77\}\\
&\searrow&A^{1}_3=\{16,22,33,39,54,62,68,74\}\\
\end{array}$ in Consequence III of Theorem 3.6.\end{center}

It is underlined that these twenty-four treatments can be equally divided
into three pairwise disjoint sets, so that arbitrary two different
treatments inside each of these three sets are first associates of
the pseudo-$L^{*}_{[9-1-5]}(9)$ association scheme.

Select the remaining treatments except $1$ from the pseudo-$L^{*}_{[9-1-5]}(9)$
association scheme, those corresponding consequences are apparently deduced in a similar manner.

Through calculations, there are totally twenty-seven distinct common
transversals of $\mathbf{L_{1}, L_{2}, L_{3}}$, $\mathbf{L_{4}, L_{5}}$,
such that they fall into three paralleling classifications of
nine transversals each. In this case, each transversal $B^{*1}_m$ and
each transversal $B^{*2}_n$ are every row and every column of $\mathbf{T_{1}^{*}}$, respectively,
$m, n=1,\cdots,9$, these obtain Classifications 1 and 2; each transversal $B^{*3}_m$
is the set of nine treatments of $\mathbf{T_{1}^{*}}$ altogether corresponding to
the same symbol of $\mathbf{L^{*}_1}$, this obtains Classification 3. Where
\begin{center}$\mathbf{T_{1}^{*}}=\left(
\begin{array}{ccccccccc}
1 &11&21&31&41&51&61&71&81\\
20&3 &10&50&33&40&80&63&70\\
12&19&2 &42&49&32&72&79&62\\
58&68&78&7 &17&27&28&38&48\\
77&60&67&26&9 &16&47&30&37\\
69&76&59&18&25&8 &39&46&29\\
34&44&54&55&65&75&4 &14&24\\
53&36&43&74&57&64&23&6 &13\\
45&52&35&66&73&56&15&22&5
\end{array}
\right)$,

$\mathbf{L^{*}_1}=\left(
\begin{array}{ccccccccc}
a^{*}&b^{*}&c^{*}&d^{*}&e^{*}&f^{*}&g^{*}&h^{*}&i^{*}\\
i^{*}&g^{*}&h^{*}&c^{*}&a^{*}&b^{*}&f^{*}&d^{*}&e^{*}\\
e^{*}&f^{*}&d^{*}&h^{*}&i^{*}&g^{*}&b^{*}&c^{*}&a^{*}\\
c^{*}&a^{*}&b^{*}&f^{*}&d^{*}&e^{*}&i^{*}&g^{*}&h^{*}\\
h^{*}&i^{*}&g^{*}&b^{*}&c^{*}&a^{*}&e^{*}&f^{*}&d^{*}\\
d^{*}&e^{*}&f^{*}&g^{*}&h^{*}&i^{*}&a^{*}&b^{*}&c^{*}\\
b^{*}&c^{*}&a^{*}&e^{*}&f^{*}&d^{*}&h^{*}&i^{*}&g^{*}\\
g^{*}&h^{*}&i^{*}&a^{*}&b^{*}&c^{*}&d^{*}&e^{*}&f^{*}\\
f^{*}&d^{*}&e^{*}&i^{*}&g^{*}&h^{*}&c^{*}&a^{*}&b^{*}
\end{array}
\right)$.\end{center}

In short, $9^2$ treatments of the pseudo-$L^{*}_{[9-1-5]}(9)$ association
scheme can be exactly divided into nine pairwise parallel transversals
$B^{*i}_1, \cdots, B^{*i}_9$ inside the $i$-th classification, $i=1,2,3$.
Moreover, two transversals of different classifications
have exactly one common treatment. By Consequence I of Theorem 3.6,
it is asserted that its $9^2$ treatments can form one $L_{3}(9)$ association scheme.

Because of $9>(9-1-5-1)^2$, calculating $\mathbf{C^{T; L_{1}, L_{2}, L_{3}, L_{4}, L_{5}}_{t}}$,
we know that twenty-four treatments, which are first associates of $t$ in the
pseudo-$L^{*}_{[9-1-5]}(9)$ association scheme, can be equally divided
into three pairwise disjoint sets, and pay attention to Theorem 3.9 again.
With respect to arbitrary two distinct treatments, which are first associates of the
pseudo-$L^{*}_{[9-1-5]}(9)$ association scheme, due to Consequence IV of Theorem 3.6,
it is shown that the common transversal of $\mathbf{L_{1}, L_{2}, L_{3}, L_{4}, L_{5}}$
containing them is sole one, that is, the set of cardinality $9$ containing
them is unique. Further, we check $\mathbf{(v)}$ in Theorem 3.12 can be satisfied.
By Corollary 3.7, thus it is inferred that the inducing pseudo-$L^{*}_{[9-1-5]}(9)$ association
scheme must be an $L_{3}(9)$ association scheme.

By Theorem 3.12, $\mathbf{L_6, L_7, L_8}$ can be added to the
$POL(9,5)=\{\mathbf{L_1, L_2, L_3, L_4, L_5}\}$ to obtain a $POL(9,8)$. Here
\begin{center}
$\mathbf{L_6}=\left(
\begin{array}{ccccccccc}
a&c&b&g&i&h&d&f&e\\
b&a&c&h&g&i&e&d&f\\
c&b&a&i&h&g&f&e&d\\
d&f&e&a&c&b&g&i&h\\
e&d&f&b&a&c&h&g&i\\
f&e&d&c&b&a&i&h&g\\
g&i&h&d&f&e&a&c&b\\
h&g&i&e&d&f&b&a&c\\
i&h&g&f&e&d&c&b&a
\end{array}
\right),\ \ \mathbf{L_7}=\left(
\begin{array}{ccccccccc}
a&c&b&g&i&h&d&f&e\\
c&b&a&i&h&g&f&e&d\\
b&a&c&h&g&i&e&d&f\\
g&i&h&d&f&e&a&c&b\\
i&h&g&f&e&d&c&b&a\\
h&g&i&e&d&f&b&a&c\\
d&f&e&a&c&b&g&i&h\\
f&e&d&c&b&a&i&h&g\\
e&d&f&b&a&c&h&g&i
\end{array}
\right)$,

$\mathbf{L_8}=\left(
\begin{array}{ccccccccc}
a&d&g&h&b&e&f&i&c\\
h&b&e&f&i&c&a&d&g\\
f&i&c&a&d&g&h&b&e\\
i&c&f&d&g&a&b&e&h\\
d&g&a&b&e&h&i&c&f\\
b&e&h&i&c&f&d&g&a\\
e&h&b&c&f&i&g&a&d\\
c&f&i&g&a&d&e&h&b\\
g&a&d&e&h&b&c&f&i
\end{array}
\right)$.\end{center}

{\bf Case 2: w=4. } Applying $\mathbf{L_{1}, L_{2}, L_{3}, L_{4}}$, we
can construct one $\mathbf{L_{6}(9)}$ association scheme and acquire one
inducing pseudo-$L^{*}_{[9-1-4]}(9)$ association scheme.

Taking any treatment $t$, we may compute $\mathbf{C^{T; L_{1}, L_{2}, L_{3}, L_{4}}_{t}}$.
We single the treatment $1$ out as an instance. There are thirty-two treatments
which are first associates of $1$ in the pseudo-$L^{*}_{[9-1-4]}(9)$ association
scheme, such that they and $1$ are together located in
$$\mathbf{\Delta^{T; L_{1}, L_{2}, L_{3}, L_{4}}_{1}}=\left(
\begin{array}{ccccccccc}
1&  &  &  &  &  &  &  &  \\
 &11&12&  &  &  &16&  &18\\
 &20&21&22&23&  &  &  &  \\
 &  &30&31&  &33&34&  &  \\
 &  &39&  &41&  &  &44&45\\
 &  &  &49&  &51&  &53&54\\
 &56&  &58&  &  &61&62&  \\
 &  &  &  &68&69&70&71&  \\
 &74&  &  &77&78&  &  &81
\end{array}
\right).$$

Next we will write the collection of all the common transversals of
$\mathbf{L_{1}, L_{2}, L_{3}, L_{4}}$ containing $1$, i.e.,
$\mathbf{C^{T; L_{1}, L_{2}, L_{3}, L_{4}}_{1}}=
\mathbf{C^{T; L_{1}, L_{2}, L_{3}, L_{4}}_{1--11}}\cup
\mathbf{C^{T; L_{1}, L_{2}, L_{3}, L_{4}}_{1--12}}\cup \mathbf{C^{T;
L_{1}, L_{2}, L_{3}, L_{4}}_{1--16}}\cup \mathbf{C^{T; L_{1}, L_{2},
L_{3}, L_{4}}_{1--18}}$. Here\begin{center}
$\mathbf{C^{T; L_{1}, L_{2}, L_{3}, L_{4}}_{1--11}}=\{$(1)
1--11--21--31--41--51--61--71--81.\;\;(2)
1--11--21--34--44--54--58--68--78.$\}$,

$\mathbf{C^{T; L_{1}, L_{2}, L_{3}, L_{4}}_{1--12}}=\{$(1)
1--12--20--34--45--53--58--69--77.\;\;(2)
1--12--20--33--41--49--62--70--81.$\}$,

$\mathbf{C^{T; L_{1}, L_{2}, L_{3}, L_{4}}_{1--16}}=\{$(1)
1--16--22--33--39--54--62--68--74.\;\;(2)
1--16--22--30--45--51--56--71--77.$\}$,

$\mathbf{C^{T; L_{1}, L_{2}, L_{3}, L_{4}}_{1--18}}=\{$(1)
1--18--23--30--44--49--56--70--78.\;\;(2)
1--18--23--31--39--53--61--69--74.$\}$.\end{center}

It is highlighted that these thirty-two treatments can be equally divided
into four pairwise disjoint sets, so that arbitrary two different
treatments inside each of these four sets are first associates of
the pseudo-$L^{*}_{[9-1-4]}(9)$ association scheme. Moreover,
there are two distinct types of these four sets, shown
in \begin{center}Type $u \left\{\begin{array}{l}
\hbox{(u) of $\mathbf{C^{T; L_{1}, L_{2}, L_{3}, L_{4}}_{1--11}}$}\\
\hbox{(u) of $\mathbf{C^{T; L_{1}, L_{2}, L_{3}, L_{4}}_{1--12}}$}\\
\hbox{(u) of $\mathbf{C^{T; L_{1}, L_{2}, L_{3}, L_{4}}_{1--16}}$}\\
\hbox{(u) of $\mathbf{C^{T; L_{1}, L_{2}, L_{3}, L_{4}}_{1--18}}$}\\
\end{array}\right.$, here $u=1,2$. \end{center}

Similarly to $1$, we choose the remainders from the pseudo-$L^{*}_{[9-1-4]}(9)$
association scheme, those corresponding conclusions are obviously inferred.

It is known from $\mathbf{C^{T; L_{1}, L_{2}, L_{3}, L_{4}}_{1--11}}$:
there exist two distinct treatments $1$ and $11$, which are first
associates of the pseudo-$L^{*}_{[9-1-4]}(9)$ association scheme,
such that the common transversal of $\mathbf{L_{1}, L_{2}, L_{3}, L_{4}}$
containing them is non-unique. Next we shall check whether the inducing pseudo-$L^{*}_{[9-1-4]}(9)$
association scheme is an $L_{4}(9)$ association scheme.

(a) When the following confining conditions are satisfied, its eighty-one
treatments can form two distinct $L_{4}(9)$ association schemes.

After calculating $\mathbf{C^{T; L_{1}, L_{2}, L_{3}, L_{4}}_{t}}$, we are able to
find out seventy-two distinct common transversals of $\mathbf{L_{1}, L_{2}, L_{3}, L_{4}}$
from the pseudo-$L^{*}_{[9-1-4]}(9)$ association scheme,
so that they fall into eight paralleling classifications of nine transversals each.

In the circumstances, nine pairwise parallel transversals $B^{*i}_1, \cdots, B^{*i}_9$
are identical to those of {\bf Case 1}, $i=1,2,3$;
each transversal $B^{*4}_m$ is the set of nine
treatments of $\mathbf{T_{1}^{*}}$ altogether corresponding to the same symbol
of $\mathbf{L^{*}_2}$, this obtains Classification 4, $m=1,\cdots,9$.
Moreover, two transversals of different classifications
have exactly one common treatment. Afterwards,
each transversal $B^{*5}_m$ and each transversal $B^{*6}_n$ are every row and every
column of $\mathbf{T_{2}^{*}}$, respectively, $m, n=1,\cdots,9$, these obtain
Classifications 5 and 6; each transversal $B^{*(k+6)}_m$ is the set of nine
treatments of $\mathbf{T_{2}^{*}}$ altogether corresponding to the same symbol
of $\mathbf{L^{*}_k}$, this obtains Classification ($k+6$), $k=1,2$.
Meanwhile, two transversals of distinct classifications
have exactly one common treatment. Where
$\mathbf{T_{1}^{*}}$ and $\mathbf{L^{*}_1}$ are as above,
\begin{center} $\mathbf{T_{2}^{*}}=\left(
\begin{array}{ccccccccc}
1 &11&21&58&68&78&34&44&54\\
74&57&64&50&33&40&26&9 &16\\
39&46&29&15&22&5 &72&79&62\\
31&41&51&7 &17&27&55&65&75\\
23&6 &13&80&63&70&47&30&37\\
69&76&59&45&52&35&12&19&2 \\
61&71&81&28&38&48&4 &14&24\\
53&36&43&20&3 &10&77&60&67\\
18&25&8 &66&73&56&42&49&32
\end{array}
\right)$.

$\mathbf{L^{*}_2}=\left(
\begin{array}{ccccccccc}
a^{*}&b^{*}&c^{*}&d^{*}&e^{*}&f^{*}&g^{*}&h^{*}&i^{*}\\
e^{*}&f^{*}&d^{*}&h^{*}&i^{*}&g^{*}&b^{*}&c^{*}&a^{*}\\
i^{*}&g^{*}&h^{*}&c^{*}&a^{*}&b^{*}&f^{*}&d^{*}&e^{*}\\
b^{*}&c^{*}&a^{*}&e^{*}&f^{*}&d^{*}&h^{*}&i^{*}&g^{*}\\
f^{*}&d^{*}&e^{*}&i^{*}&g^{*}&h^{*}&c^{*}&a^{*}&b^{*}\\
g^{*}&h^{*}&i^{*}&a^{*}&b^{*}&c^{*}&d^{*}&e^{*}&f^{*}\\
c^{*}&a^{*}&b^{*}&f^{*}&d^{*}&e^{*}&i^{*}&g^{*}&h^{*}\\
d^{*}&e^{*}&f^{*}&g^{*}&h^{*}&i^{*}&a^{*}&b^{*}&c^{*}\\
h^{*}&i^{*}&g^{*}&b^{*}&c^{*}&a^{*}&e^{*}&f^{*}&d^{*}
\end{array}
\right)$.
\end{center}
It is easily seen that $\mathbf{L^{*}_1}$ and $\mathbf{L^{*}_2}$ are orthogonal.

As a whole, $9^2$ treatments of the pseudo-$L^{*}_{[9-1-4]}(9)$ association
scheme can be exactly divided into nine pairwise parallel transversals
$B^{*i}_1, \cdots, B^{*i}_9$ inside the $i$-th classification, here
$i=4v-3,4v-2,4v-1,4v; v=1,2$. By Consequence I of Theorem 3.6, it can be illustrated
that its $9^2$ treatments can form two distinct $L_{4}(9)$ association schemes.

(b) When anyone of five conditions of Theorem 3.5 is established,
the pseudo-$L^{*}_{[9-1-4]}(9)$ association scheme is not an $L_{4}(9)$ association scheme.

It is obviously conscious that $\mathbf{(I)}$ and $\mathbf{(III)}$ have already been
established in Theorem 3.5. We notice $B^{*1}_1\cap B^{*5}_1=\{$1--11--21--31--41--51--61--71--81.$\}\cap \{$1--11--21--34--44--54--58--68--78.$\}=\{1, 11, 21\}$.
Let $t^{*}_2=1$, none of (1) in $\mathbf{C^{T; L_{1}, L_{2}, L_{3}, L_{4}}_{1--11}}$, (1) in
$\mathbf{C^{T; L_{1}, L_{2}, L_{3}, L_{4}}_{1--12}}$, (1) in
$\mathbf{C^{T; L_{1}, L_{2}, L_{3}, L_{4}}_{1--16}}$, (1) in
$\mathbf{C^{T; L_{1}, L_{2}, L_{3}, L_{4}}_{1--18}}$ is equal
to $B^{*5}_1$. Thus these elucidate Conditions $\mathbf{(II)}$ and $\mathbf{(IV)}$
of Theorem 3.5 are established. Set $t^{*}_1=11$ and $t^{*}_2=1$, it is seen that
$(\{11\}\cup\{$3--19--36--44--52--60--68--76$\})\cap
(\{1\}\cup\{$12--20--33--41--49--62--70--81$\})=\emptyset$
and $(\{11\}\cup\{$3--19--36--44--52--60--68--76$\})\cap
(\{1\}\cup\{$16--22--30--45--51--56--71--77\}$)=\emptyset$
hold; and it is detected from $\mathbf{C^{T; L_{1}, L_{2}, L_{3}, L_{4}}_{1--11}}$ that
the common transversal of $\mathbf{L_{1}, L_{2}, L_{3}, L_{4}}$ containing $11$ and $1$ is
not sole one, that is, the set of cardinality $9$ containing them is non-unique.
It follows that $\mathbf{(V)}$ holds. Due to Theorem 3.5, but the
inducing pseudo-$L^{*}_{[9-1-4]}(9)$ association scheme is not an $L_{4}(9)$ association scheme.

Even if two treatments $31$ and $34$ are second associates of the
pseudo-$L^{*}_{[9-1-4]}(9)$ association scheme, we apparently have
$(\{31\}\cup\{$1--11--21--41--51--61--71--81$\})\cap (\{34\}\cup\{$1--11--21--44--54--58--68--78\}$)
=\{1, 11, 21\}$.

{\bf Case 3: w=3. } Applying $\mathbf{L_{1}, L_{2}, L_{3}}$, we
can construct one $\mathbf{L_{5}(9)}$ association scheme and obtain one
inducing pseudo-$L^{*}_{[9-1-3]}(9)$ association scheme.

Since it is very certain that all the common transversals of
$\mathbf{L_{1}, L_{2}, L_{3}, L_{4}}$ are the common transversals of $\mathbf{L_{1}, L_{2}, L_{3}}$.
Through calculations, it is advertent that forty treatments are
first associates of any treatment in the
pseudo-$L^{*}_{[9-1-3]}(9)$ association scheme, they can be
equally divided into five pairwise disjoint sets, so that arbitrary two different
treatments inside each of these five sets are first associates
of the pseudo-$L^{*}_{[9-1-3]}(9)$ association scheme.
For example, let (1) of $\mathbf{C^{T; L_{1}, L_{2}, L_{3}}_{1--17}}$ be
$\{$1--17--24--36--40--47--59--66--79.$\}$,
there are at least two distinct types of these five sets concerning the treatment $1$, exhibited
in \begin{center}Type $u \left\{\begin{array}{l}
\hbox{(u) of $\mathbf{C^{T; L_{1}, L_{2}, L_{3}, L_{4}}_{1--11}}$}\\
\hbox{(u) of $\mathbf{C^{T; L_{1}, L_{2}, L_{3}, L_{4}}_{1--12}}$}\\
\hbox{(u) of $\mathbf{C^{T; L_{1}, L_{2}, L_{3}, L_{4}}_{1--16}}$}\\
\hbox{(u) of $\mathbf{C^{T; L_{1}, L_{2}, L_{3}, L_{4}}_{1--18}}$}\\
\hbox{(1) of $\mathbf{C^{T; L_{1}, L_{2}, L_{3}}_{1--17}}$}\\
\end{array}\right.$, here $u=1,2$. \end{center}

Next we shall discuss if the inducing pseudo-$L^{*}_{[9-1-3]}(9)$
association scheme is an $L_{5}(9)$ association scheme.

(a) When the following confining conditions are satisfied, its eighty-one
treatments can form two distinct $L_{5}(9)$ association schemes.

Also since $9^2$ treatments of the pseudo-$L^{*}_{[9-1-4]}(9)$ association
scheme can be exactly divided into nine pairwise parallel transversals
$B^{*i}_1, \cdots, B^{*i}_9$ inside the $i$-th classification,
of course, $9^2$ treatments of the pseudo-$L^{*}_{[9-1-3]}(9)$ association
scheme can also be divided into nine pairwise parallel transversals
$B^{*i}_1, \cdots, B^{*i}_9$ inside that above $i$-th classification.
Calculating $\mathbf{C^{T; L_{1}, L_{2}, L_{3}}_{t}}$,
we are able to find nine pairwise parallel common transversals of $\mathbf{L_{1}, L_{2}, L_{3}}$,
which constitute another paralleling classification.

In this case, nine pairwise parallel transversals $B^{*i}_1, \cdots, B^{*i}_9$
are identical to those of (a) within {\bf Case 2}, $i=1,2,3,4$;
each transversal $B^{*5}_m$ is the set of nine treatments of $\mathbf{T_{1}^{*}}$
altogether corresponding to the same symbol of $\mathbf{L^{*}_3}$,
this obtains Classification 5, $m=1,\cdots,9$.
Moreover, two transversals of different classifications
have exactly one common treatment. Subsequently,
Classifications 6 to 9 are correspondingly identical to Classifications 5 to 8 of (a)
within {\bf Case 2}; each transversal $B^{*10}_m$ is the set of nine
treatments of $\mathbf{T_{2}^{*}}$ altogether corresponding to the same symbol
of $\mathbf{L^{*}_3}$, this obtains Classification 10, $m=1,\cdots,9$.
Meanwhile, two transversals of distinct classifications
have exactly one common treatment. Now we must point out that Classification 5
and Classification 10 are exactly of the same class. Where
$\mathbf{T_{1}^{*}, T_{2}^{*}}$ and $\mathbf{L^{*}_1, L^{*}_2}$ are as above,
\begin{center}$\mathbf{L^{*}_3}=\left(
\begin{array}{ccccccccc}
a^{*}&b^{*}&c^{*}&d^{*}&e^{*}&f^{*}&g^{*}&h^{*}&i^{*}\\
h^{*}&i^{*}&g^{*}&b^{*}&c^{*}&a^{*}&e^{*}&f^{*}&d^{*}\\
f^{*}&d^{*}&e^{*}&i^{*}&g^{*}&h^{*}&c^{*}&a^{*}&b^{*}\\
i^{*}&g^{*}&h^{*}&c^{*}&a^{*}&b^{*}&f^{*}&d^{*}&e^{*}\\
d^{*}&e^{*}&f^{*}&g^{*}&h^{*}&i^{*}&a^{*}&b^{*}&c^{*}\\
b^{*}&c^{*}&a^{*}&e^{*}&f^{*}&d^{*}&h^{*}&i^{*}&g^{*}\\
e^{*}&f^{*}&d^{*}&h^{*}&i^{*}&g^{*}&b^{*}&c^{*}&a^{*}\\
c^{*}&a^{*}&b^{*}&f^{*}&d^{*}&e^{*}&i^{*}&g^{*}&h^{*}\\
g^{*}&h^{*}&i^{*}&a^{*}&b^{*}&c^{*}&d^{*}&e^{*}&f^{*}
\end{array}
\right)$.
\end{center}
It is apparently checked that $\mathbf{L^{*}_1, L^{*}_2, L^{*}_3}$
are mutually orthogonal.

In sum, $9^2$ treatments of the pseudo-$L^{*}_{[9-1-3]}(9)$ association
scheme can be exactly divided into nine pairwise parallel transversals
$B^{*i}_1, \cdots, B^{*i}_9$ inside the $i$-th classification, where
$i=5v-4,5v-3,5v-2,5v-1,5v; v=1,2$. By Consequence I of Theorem 3.6, it may be verified
that its $9^2$ treatments can form two distinct $L_{5}(9)$ association schemes.

(b) When anyone of five conditions of Theorem 3.5 is established,
the pseudo-$L^{*}_{[9-1-3]}(9)$ association scheme is not an $L_{5}(9)$ association scheme.

Still let $t^{*}_2=1$, none of (1) in $\mathbf{C^{T; L_{1}, L_{2}, L_{3}, L_{4}}_{1--11}}$, (1) in
$\mathbf{C^{T; L_{1}, L_{2}, L_{3}, L_{4}}_{1--12}}$, (1) in
$\mathbf{C^{T; L_{1}, L_{2}, L_{3}, L_{4}}_{1--16}}$, (1) in
$\mathbf{C^{T; L_{1}, L_{2}, L_{3}, L_{4}}_{1--18}}$, (1) in
$\mathbf{C^{T; L_{1}, L_{2}, L_{3}}_{1--17}}$ is equal
to $\{$1--11--21--34--44--54--58--68--78.$\}$. Being similar to (b) of {\bf Case 2},
Conditions $\mathbf{(I)}$ to $\mathbf{(V)}$ of Theorem 3.5 are satisfied, thus the
inducing pseudo-$L^{*}_{[9-1-3]}(9)$ association scheme is not an $L_{5}(9)$ association scheme.

{\bf Case 4: w=2. } Applying $\mathbf{L_{1}, L_{2}}$, we
can construct one $\mathbf{L_{4}(9)}$ association scheme and acquire one
inducing pseudo-$L^{*}_{[9-1-2]}(9)$ association scheme.

It is sure that the whole common transversals of
$\mathbf{L_{1}, L_{2}, L_{3}}$ are the common transversals of $\mathbf{L_{1}, L_{2}}$.
We emphasize that forty-eight treatments, which are first associates of any treatment in the
pseudo-$L^{*}_{[9-1-2]}(9)$ association scheme, can be equally divided into
six pairwise disjoint sets, so that arbitrary two different
treatments inside each of these six sets are first associates
of the pseudo-$L^{*}_{[9-1-2]}(9)$ association scheme.
For instance, label $\{$1--15--26--32--43--48--63--65--76.$\}$ by (1)
of $\mathbf{C^{T; L_{1}, L_{2}}_{1--15}}$,
there are at least two distinct types of these six sets concerning the treatment $1$, displayed
in \begin{center}Type $u \left\{\begin{array}{l}
\hbox{(u) of $\mathbf{C^{T; L_{1}, L_{2}, L_{3}, L_{4}}_{1--11}}$}\\
\hbox{(u) of $\mathbf{C^{T; L_{1}, L_{2}, L_{3}, L_{4}}_{1--12}}$}\\
\hbox{(u) of $\mathbf{C^{T; L_{1}, L_{2}, L_{3}, L_{4}}_{1--16}}$}\\
\hbox{(u) of $\mathbf{C^{T; L_{1}, L_{2}, L_{3}, L_{4}}_{1--18}}$}\\
\hbox{(1) of $\mathbf{C^{T; L_{1}, L_{2}, L_{3}}_{1--17}}$}\\
\hbox{(1) of $\mathbf{C^{T; L_{1}, L_{2}}_{1--15}}$}\\
\end{array}\right.$, here $u=1,2$. \end{center}

Next we shall determine whether the inducing pseudo-$L^{*}_{[9-1-2]}(9)$
association scheme is an $L_{6}(9)$ association scheme.

(a) When the following confining conditions are satisfied, its eighty-one
treatments can form one $L_{6}(9)$ association scheme.

It is noted that $9^2$ treatments of the pseudo-$L^{*}_{[9-1-3]}(9)$ association
scheme can be exactly divided into nine pairwise parallel transversals
$B^{*i}_1, \cdots, B^{*i}_9$ inside the $i$-th classification, $i=1,2,3,4,5$.
Via computing, we may look for nine pairwise parallel common transversals
of $\mathbf{L_{1}, L_{2}}$, which constitute the sixth paralleling classification.

In the circumstances, nine pairwise parallel transversals $B^{*i}_1, \cdots, B^{*i}_9$
are identical to those of (a) within {\bf Case 3}, $i=1,2,3,4,5$;
each transversal $B^{*6}_m$ is the set of nine treatments of $\mathbf{T_{1}^{*}}$
altogether corresponding to the same symbol of $\mathbf{L^{*}_4}$,
this obtains Classification 6, $m=1,\cdots,9$. Where
$\mathbf{T_{1}^{*}}$ and $\mathbf{L^{*}_1, L^{*}_2, L^{*}_3}$ are as above,
\begin{center}$\mathbf{L^{*}_4}=\left(
\begin{array}{ccccccccc}
a^{*}&b^{*}&c^{*}&d^{*}&e^{*}&f^{*}&g^{*}&h^{*}&i^{*}\\
f^{*}&d^{*}&e^{*}&i^{*}&g^{*}&h^{*}&c^{*}&a^{*}&b^{*}\\
h^{*}&i^{*}&g^{*}&b^{*}&c^{*}&a^{*}&e^{*}&f^{*}&d^{*}\\
e^{*}&f^{*}&d^{*}&h^{*}&i^{*}&g^{*}&b^{*}&c^{*}&a^{*}\\
g^{*}&h^{*}&i^{*}&a^{*}&b^{*}&c^{*}&d^{*}&e^{*}&f^{*}\\
c^{*}&a^{*}&b^{*}&f^{*}&d^{*}&e^{*}&i^{*}&g^{*}&h^{*}\\
i^{*}&g^{*}&h^{*}&c^{*}&a^{*}&b^{*}&f^{*}&d^{*}&e^{*}\\
b^{*}&c^{*}&a^{*}&e^{*}&f^{*}&d^{*}&h^{*}&i^{*}&g^{*}\\
d^{*}&e^{*}&f^{*}&g^{*}&h^{*}&i^{*}&a^{*}&b^{*}&c^{*}
\end{array}
\right)$.
\end{center}
It is obviously verified that $\mathbf{L^{*}_1, L^{*}_2, L^{*}_3, L^{*}_4}$
are mutually orthogonal.

On balance, $9^2$ treatments of the pseudo-$L^{*}_{[9-1-2]}(9)$ association
scheme can be exactly divided into nine pairwise parallel transversals
$B^{*i}_1, \cdots, B^{*i}_9$ inside the $i$-th classification, here
$i=1,2,3,4,5,6$. Moreover, two transversals of different classifications
have exactly one common treatment. By Consequence I of Theorem 3.6, it is illustrated
that its $9^2$ treatments can form one $L_{6}(9)$ association scheme.

(b) When anyone of five conditions of Theorem 3.5 is established,
the pseudo-$L^{*}_{[9-1-2]}(9)$ association scheme is not an $L_{6}(9)$ association scheme.

Still put $t^{*}_2=1$, none of (1) in $\mathbf{C^{T; L_{1}, L_{2}, L_{3}, L_{4}}_{1--11}}$, (1) in
$\mathbf{C^{T; L_{1}, L_{2}, L_{3}, L_{4}}_{1--12}}$, (1) in
$\mathbf{C^{T; L_{1}, L_{2}, L_{3}, L_{4}}_{1--16}}$, (1) in
$\mathbf{C^{T; L_{1}, L_{2}, L_{3}, L_{4}}_{1--18}}$, (1) in
$\mathbf{C^{T; L_{1}, L_{2}, L_{3}}_{1--17}}$, (1) in
$\mathbf{C^{T; L_{1}, L_{2}}_{1--15}}$ is equal
to $\{$1--11--21--34--44--54--58--68--78.$\}$. Similarly to (b) of {\bf Case 3},
Conditions $\mathbf{(I)}$ to $\mathbf{(V)}$ are satisfied in Theorem 3.5, thus the
inducing pseudo-$L^{*}_{[9-1-2]}(9)$ association scheme is not an $L_{6}(9)$ association scheme.

{\bf Case 5: w=1. } Applying $\mathbf{L_{1}}$, we
can construct one $\mathbf{L_{3}(9)}$ association scheme and obtain one
inducing pseudo-$L^{*}_{[9-1-1]}(9)$ association scheme.

It is of course true that all the common transversals of
$\mathbf{L_{1}, L_{2}}$ are the transversals of $\mathbf{L_{1}}$.
We stress that fifty-six treatments, which are first associates of any treatment in the
pseudo-$L^{*}_{[9-1-1]}(9)$ association scheme, can be equally divided into
seven pairwise disjoint sets. Such as, we denote $\{$1--14--27--29--42--52--57--67--80.$\}$ by (1)
of $\mathbf{C^{T; L_{1}}_{1--14}}$, there are at least two distinct types of these
seven sets concerning the treatment $1$, shown
in \begin{center}Type $u \left\{\begin{array}{l}
\hbox{(u) of $\mathbf{C^{T; L_{1}, L_{2}, L_{3}, L_{4}}_{1--11}}$}\\
\hbox{(u) of $\mathbf{C^{T; L_{1}, L_{2}, L_{3}, L_{4}}_{1--12}}$}\\
\hbox{(u) of $\mathbf{C^{T; L_{1}, L_{2}, L_{3}, L_{4}}_{1--16}}$}\\
\hbox{(u) of $\mathbf{C^{T; L_{1}, L_{2}, L_{3}, L_{4}}_{1--18}}$}\\
\hbox{(1) of $\mathbf{C^{T; L_{1}, L_{2}, L_{3}}_{1--17}}$}\\
\hbox{(1) of $\mathbf{C^{T; L_{1}, L_{2}}_{1--15}}$}\\
\hbox{(1) of $\mathbf{C^{T; L_{1}}_{1--14}}$}\\
\end{array}\right.$, here $u=1,2$. \end{center}

Next we shall judge if the inducing pseudo-$L^{*}_{[9-1-1]}(9)$
association scheme is an $L_{7}(9)$ association scheme.

(a) When the following confining conditions are satisfied, its eighty-one
treatments can form one $L_{7}(9)$ association scheme.

It is noticed that $9^2$ treatments of the pseudo-$L^{*}_{[9-1-2]}(9)$ association
scheme can be exactly divided into nine pairwise parallel transversals
$B^{*i}_1, \cdots, B^{*i}_9$ inside the $i$-th classification, $i=1,2,3,4,5,6$.
Through calculations, we may search for nine pairwise parallel transversals
of $\mathbf{L_{1}}$, which make up the seventh paralleling classification.

In this case, nine pairwise parallel transversals $B^{*i}_1, \cdots, B^{*i}_9$
are identical to those of (a) within {\bf Case 4}, $i=1,2,3,4,5,6$;
each transversal $B^{*7}_m$ is the set of nine treatments of $\mathbf{T_{1}^{*}}$
altogether corresponding to the same symbol of $\mathbf{L^{*}_5}$,
this obtains Classification 7, $m=1,\cdots,9$. Where
$\mathbf{T_{1}^{*}}$ and $\mathbf{L^{*}_1, L^{*}_2, L^{*}_3, L^{*}_4}$ are as above,
\begin{center}$\mathbf{L^{*}_5}=\left(
\begin{array}{ccccccccc}
a^{*}&b^{*}&c^{*}&d^{*}&e^{*}&f^{*}&g^{*}&h^{*}&i^{*}\\
d^{*}&e^{*}&f^{*}&g^{*}&h^{*}&i^{*}&a^{*}&b^{*}&c^{*}\\
g^{*}&h^{*}&i^{*}&a^{*}&b^{*}&c^{*}&d^{*}&e^{*}&f^{*}\\
h^{*}&i^{*}&g^{*}&b^{*}&c^{*}&a^{*}&e^{*}&f^{*}&d^{*}\\
b^{*}&c^{*}&a^{*}&e^{*}&f^{*}&d^{*}&h^{*}&i^{*}&g^{*}\\
e^{*}&f^{*}&d^{*}&h^{*}&i^{*}&g^{*}&b^{*}&c^{*}&a^{*}\\
f^{*}&d^{*}&e^{*}&i^{*}&g^{*}&h^{*}&c^{*}&a^{*}&b^{*}\\
i^{*}&g^{*}&h^{*}&c^{*}&a^{*}&b^{*}&f^{*}&d^{*}&e^{*}\\
c^{*}&a^{*}&b^{*}&f^{*}&d^{*}&e^{*}&i^{*}&g^{*}&h^{*}
\end{array}
\right)$.
\end{center}
It can be tested that $\mathbf{L^{*}_1, L^{*}_2, L^{*}_3, L^{*}_4, L^{*}_5}$
are mutually orthogonal.

In sum, $9^2$ treatments of the pseudo-$L^{*}_{[9-1-1]}(9)$ association
scheme can be exactly divided into nine pairwise parallel transversals
$B^{*i}_1, \cdots, B^{*i}_9$ inside the $i$-th classification, here
$i=1,2,3,4,5,6,7$. Moreover, two transversals of different classifications
have exactly one common treatment. By Consequence I of Theorem 3.6, it is elucidated
that its $9^2$ treatments can form one $L_{7}(9)$ association scheme.

(b) When anyone of five conditions of Theorem 3.5 is established,
the pseudo-$L^{*}_{[9-1-1]}(9)$ association scheme is not an $L_{7}(9)$ association scheme.

Similarly to (b) of {\bf Case 4}, Conditions $\mathbf{(I)}$ to $\mathbf{(V)}$ are
satisfied in Theorem 3.5, thus the inducing pseudo-$L^{*}_{[9-1-1]}(9)$ association
scheme is not an $L_{7}(9)$ association scheme.

\section{ Some Research Problems }

Set $POL(s,w)=\{\mathbf{L_1, \cdots, L_{w}}\}, s\geq 10$,
as long as the $POL(s,w)$ can be extended to a
$POL(s,s-1)$. Whether the order $s$ is a
prime power or not, theoretically speaking, the remaining $s-1-w$ pairwise
orthogonal Latin squares of order $s$ could be added to the $POL(s,w)$ to
obtain a $POL(s,s-1)$. However, it is possible that there are most tremendous
common transversals of $\mathbf{L_1, \cdots, L_{w}}$. The most important open problem is
how to select $s(s-1-w)$ distinct common transversals
with efficiency from much more than these $s(s-1-w)$ ones of the inducing
pseudo-$L^{*}_{[s-1-w]}(s)$ association scheme, so that they fall into
$(s-1-w)$ paralleling classifications of $s$ transversals each,
two transversals of different classifications have exactly one common treatment.
It follows that $s^2$ treatments of the pseudo-$L^{*}_{[s-1-w]}(s)$ association
scheme can form an $L_{s-1-w}(s)$ association scheme.\\

\noindent\textbf{Acknowledgements}

This work was supported by the National Natural Science Foundation of China (Grant
Nos. 11971104 and 11871417).

\end{document}